\documentclass[10pt]{article}

\usepackage{jmlr2e}
\usepackage[utf8]{inputenc} 
\usepackage[T1]{fontenc}    
\usepackage{lmodern}
\usepackage{hyperref}       
\usepackage{url}            
\usepackage{booktabs}       
\usepackage{amsfonts}       
\usepackage{nicefrac}       
\usepackage{microtype}      
\usepackage{amsmath,dsfont,empheq}
\usepackage{algorithm}
\usepackage{algpseudocode}
\usepackage{syntonly,bm,wrapfig,epsfig}
\usepackage{xcolor,subcaption}

\usepackage{multibib}
\usepackage{graphicx}
\newcites{apndx}{References}

\usepackage{balance}

\usepackage{multirow}
\def\1{\bm{1}}

\DeclareMathAlphabet{\mathsfit}{\encodingdefault}{\sfdefault}{m}{sl}
\SetMathAlphabet{\mathsfit}{bold}{\encodingdefault}{\sfdefault}{bx}{n}

\newcommand{\RR}{\mathbb{R}}
\newcommand{\EE}{\mathbb{E}}
\newcommand{\cL}{\mathcal{L}}

\newcommand{\ttheta}{\tilde{\theta}}
\newcommand{\tdelta}{\tilde{\delta}}

\newcommand{\htau}{\hat{\tau}}

\DeclarePairedDelimiter{\dotp}{\langle}{\rangle}

\newtheorem{assumption}{Assumption}
\renewcommand\footnotemark{}

\algrenewcommand{\algorithmiccomment}[1]{\hskip0em$\triangleright$ #1}
\algblock{ParFor}{EndParFor}
\algnewcommand\algorithmicparfor{\textbf{for}}
\algnewcommand\algorithmicpardo{\textbf{do in parallel}}
\algnewcommand\algorithmicendparfor{\textbf{end for}}
\algrenewtext{ParFor}[1]{\algorithmicparfor\ #1\ \algorithmicpardo}
\algrenewtext{EndParFor}{\algorithmicendparfor}

\title{LASG: Lazily Aggregated Stochastic Gradients\\for Communication-Efficient Distributed Learning}

\author{\qquad\qquad Tianyi Chen$^{\star}$ \qquad Yuejiao Sun$^{\dag}$ \qquad Wotao Yin$^{\dag}$\\\\
{\centering
\addr $^{\star}$\textit{Rensselaer Polytechnic Institute - Troy, NY 12180, USA}\\
$^{\dag}$\textit{University of California - Los Angeles, Los Angeles, CA 90095, USA}\\
\email \texttt{chent18@rpi.edu}~~~\{\texttt{sunyj,wotaoyin\}@math.ucla.edu}\qquad\qquad\qquad\qquad\qquad}
\thanks{\hspace{-0.5cm}Authors are listed in alphabetical order.}
}

\ShortHeadings{LASG: Lazily Aggregated Stochastic Gradients for Communication-Efficient Distributed Learning}{Chen and Sun and Yin}
\firstpageno{1}

\begin{document}

\maketitle

\thispagestyle{empty}

\begin{abstract}
This paper targets solving distributed machine learning problems such as federated learning in a communication-efficient fashion. A class of new stochastic gradient descent (SGD) approaches have been developed, which can be viewed as the stochastic generalization to the recently developed lazily aggregated gradient (LAG) method --- justifying the name LASG. LAG adaptively predicts the contribution of each round of communication and chooses only the significant ones to perform. It saves communication while also maintains the rate of convergence. However, LAG only works with deterministic gradients, and applying it to stochastic gradients yields poor performance.
The key components of LASG are a set of new rules tailored for stochastic gradients that can be implemented either to save download, upload, or both.
The new algorithms adaptively choose between fresh and stale stochastic gradients and have convergence rates comparable to the original SGD.
LASG achieves impressive empirical performance --- it typically saves total communication by an order of magnitude.

\end{abstract} 

\section{Introduction}
We consider a cloud server and a set of $M$ mobile devices (workers) collected in ${\cal M}:=\{1,\ldots,M\}$.
Each device $m$ has its local dataset $\{\xi_n,\,n\in {\cal N}_m\}$, which defines the loss function of device $m$ as
\begin{align}
{\cal L}_m(\theta)\!:=\!\sum_{n\in{\cal N}_m}\!\ell(\theta;\xi_n),\quad m\in{\cal M}
\end{align}
where $\theta\in \mathbb{R}^p$ is the sought vector (e.g., parameters of a prediction model) and $\xi_n$ is a data sample.
For example, in linear regression, $\ell(\theta;\xi_n)$ is the square loss; 
and, in deep learning, $\ell(\theta;\xi_n)$ is the loss function of a neural network, and $\theta$ concatenates the weights.
The goal is to solve %
\begin{align}\label{eqn: problem}
	\min_{\theta\in \mathbb{R}^p}~{\cal L}(\theta)~~~{\rm with}~~~{\cal L}(\theta):=\sum_{m\in{\cal M}}{\cal L}_m(\theta).
\end{align}
Problem \eqref{eqn: problem} also arises in a number of areas, such as
multi-agent optimization \citep{nedic2009}, distributed signal processing \citep{msechu2011tsp}, and distributed machine learning \citep{dean2012}.
While our algorithms can be applied to other settings, we focus on the federated learning setting.
In this case, for bandwidth and privacy concerns, local data $\{\xi_n, n\in{\cal N}_m\}$ at each worker $m$ are not uploaded to the server, and collaboration is needed through communication between the server and workers.

\begin{figure}[t]
\centering
\includegraphics[width=8cm]{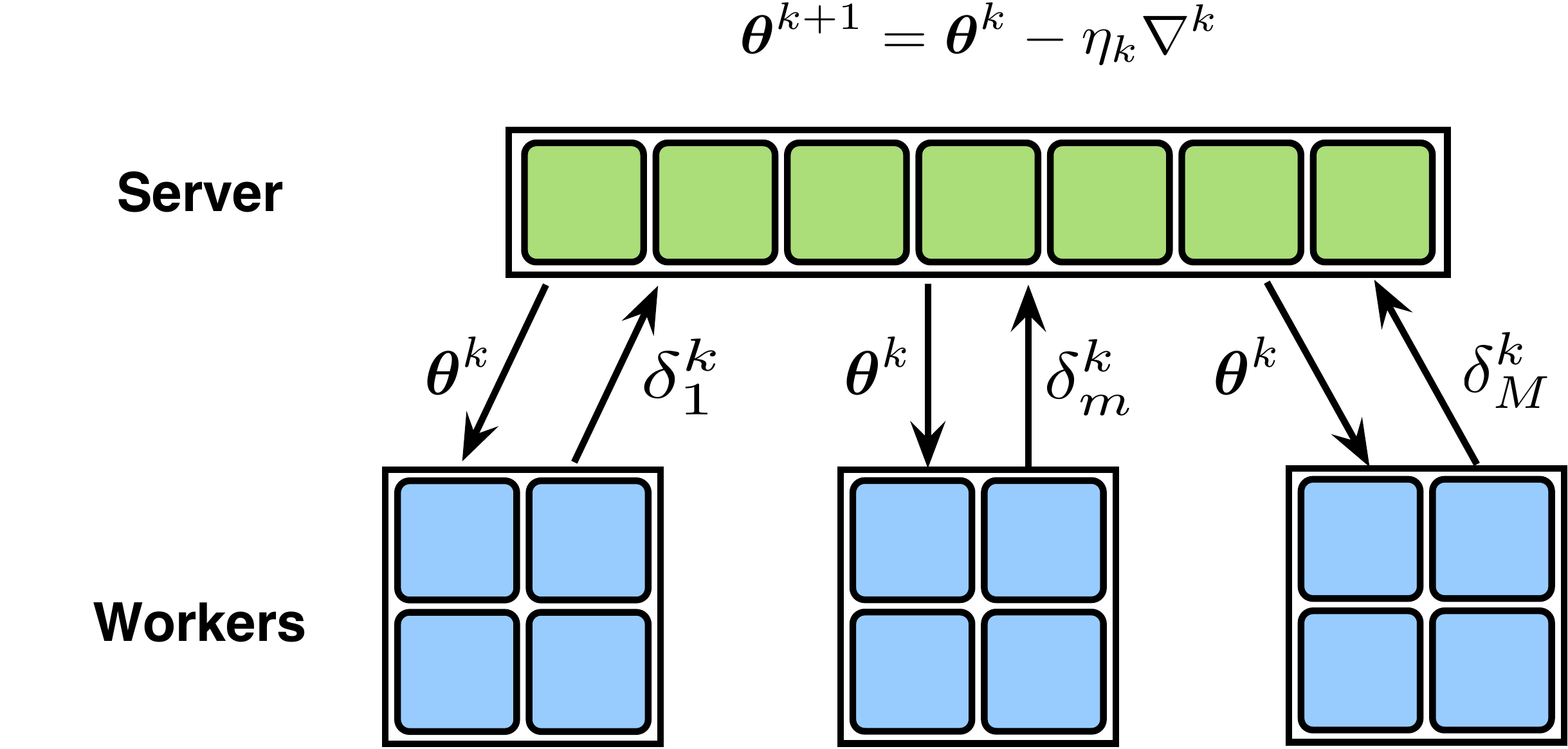}
  \caption{Generic LASG implementation.}
\label{fig:lasg-diag}
\end{figure}

To solve \eqref{eqn: problem}, we can in principle apply the distributed version of stochastic gradient descent (SGD) method.
In this case, at iteration $k$, the server broadcasts the current model $\theta^k$ to \emph{all} the workers; each worker $m$ computes $\nabla  \ell(\theta^k; \xi_m^k)$ using a randomly selected sample or a minibatch of samples $\{\xi_m^k\}\subseteq\{\xi_n, n\in{\cal N}_m\}$, and then uploads it to the server; and once receiving stochastic gradients from all workers, the server updates the model parameters via 
\begin{flalign}	\label{eqn:SGDiteration}
& {\rm \textbf{SGD}}\qquad\qquad\qquad\qquad\qquad\qquad\qquad \theta^{k+1}=\theta^k-\eta_k\sum\limits_{m\in{\cal M}} \nabla  \ell(\theta^k; \xi_m^k) &
\end{flalign}
where $\eta_k>0$ is the (possibly) time-varying stepsize used at iteration $k$.
When $\nabla  \ell(\theta^k; \xi_m^k)$ is an unbiased gradient estimator of $\cL_m(\theta)$, the convergence of SGD update \eqref{eqn:SGDiteration} is guaranteed \citep{bottou2016}.
To implement \eqref{eqn:SGDiteration}, however, the server has to communicate with \emph{all} workers to obtain fresh $\{\nabla \ell(\theta^k; \xi_m^k) \}$.
This prevents the efficient implementation of SGD in scenarios where communication between the server and the workers is costly \citep{mcmahan2017}, since latency will degrade the overall performance.
Therefore, \emph{our goal} is to find the parameter $\theta$ that minimizes \eqref{eqn: problem} with minimal communication overhead.

\subsection{Our approach}
This paper puts forward a class of new stochastic optimization methods that can considerably reduce the redundant or less informative communication of SGD.
The key motivation is that during the distributed learning process, not all communication rounds between the server and the worker are equally important. So a natural solution is to use a condition that decides whether to communicate or not. In this case, if some workers are not communicating, the server uses their stale gradients so that stale gradients replace skipped fresh gradients.

Analogous to the distributed implementation of SGD \eqref{eqn:SGDiteration}, our new approaches also aggregate the stochastic gradients from all workers, but in a fairly lazy manner. Hence, we term our algorithms as \textbf{L}azily \textbf{A}ggregated \textbf{S}tochastic \textbf{G}radient (\textbf{LASG}). LASG has the following generic update
\begin{align}\label{eqn:LASGiteration}
\theta^{k+1} = \theta^k-\eta_k\!\!\!\sum\limits_{m\in{\cal M}\backslash{\cal M}^k} \!\!\!\nabla \ell(\theta^{k-\tau_m^k};\xi_m^{k-\tau_m^k})
-\eta_k\!\!\!\sum\limits_{m\in{\cal M}^k}\!\!\! \nabla\ell(\theta^k;\xi_m^k)
\end{align}
or equivalently (see also Figure \ref{fig:lasg-diag})
 \begin{flalign}\label{eqn:LASGiteration2}
&{\rm \textbf{Generic LASG}}\qquad\qquad\qquad\qquad  \theta^{k+1} = \theta^k-\eta_k\nabla^k\quad{\rm with}\quad\nabla^k\!=\!\nabla^{k-1}+\sum\limits_{m\in{\cal M}^k}\delta_m^k&
\end{flalign}
where $\delta_m^k := \nabla \ell(\theta^k;\xi_m^k)-\nabla \ell(\theta^{k-\tau_m^k};\xi_m^{k-\tau_m^k})$ is the stochastic gradient innovation, $\tau_m^k\geq 0$ is the staleness of the gradient from worker $m$ used by $k$-iteration, and ${\cal M}^k$ is the subset of workers uploading $\delta_m^k$ at iteration $k$. The stalenesses $\{\tau_m^k\}$ are controlled by the selection of subset ${\cal M}^k$: at iteration $k$, if worker $m\notin {\cal M}^k$, the server increases the staleness by $\tau_m^{k+1} = \tau_m^{k} + 1$; otherwise, worker $m$ uploads the stochastic gradient, and the server resets $\tau_m^{k+1}=1$.

Clearly, selection of subset ${\cal M}^k$ is critical in LASG. However, the challenges are 1) the importance of each communication round is dynamic, thus a fixed or nonadaptive condition is ineffective; and 2) the condition needs to be checked either at server or at worker locally and efficiently.

To overcome these challenges, we develop two types of adaptive conditions based on message innovations. They can be chosen under different \emph{communication}, \emph{computation} and \emph{memory} requirements.
The first type is adopted by each worker (WK), and the second one by the server (PS).

\noindent\textbf{LASG-WK}: At iteration $k$, the server broadcasts $\theta^k$ to all workers; every worker computes $\nabla \ell(\theta^k;\xi_m^k)$, and checks if it belongs to ${\cal M}^k$; only the workers in ${\cal M}^k$ upload $\delta^k_m$; the server updates via \eqref{eqn:LASGiteration2}.

\noindent\textbf{LASG-PS}: At iteration $k$, the server decides ${\cal M}^k$ and sends $\theta^k$ to workers in ${\cal M}^k$; each worker $m\!\in\!{\cal M}^k$ computes $\nabla \ell(\theta^k;\xi_m^k)$ and uploads $\delta^k_m$; the rest of workers do nothing; the server updates via \eqref{eqn:LASGiteration2};

With detailed description of LASG rules deferred to Section \ref{sec.LASG}, the contributions of this paper are listed as follows.

{\bf 1)} We introduce a class of novel (quantized) stochastic optimization approaches that
reuses stale stochastic gradients to reduce redundant communication.

{\bf 2)} We establish convergence of our proposed algorithms in strongly convex and nonconvex settings even when the datasets are non-i.i.d. across workers. The convergence rates match those of SGD in the respective settings.

{\bf 3)} We confirm performance gains of our novel distributed
algorithms over some alternatives using extensive numerical tests on logistic regression and neural network training.

\subsection{Related work}
Communication-efficient distributed learning methods have gained popularity recently \citep{Nedic2018,jordan2018}. Most popular methods belong to two categories: c1) reduce the number of bits per communication
round; and, c2) save the number of communication
rounds.
For c1), methods are centered around the ideas of \emph{quantization} and \emph{sparsification}.

\noindent\textbf{Quantization} has been successfully applied to several engineering tasks employing wireless sensor networks \citep{msechu2011tsp}. In the context of distributed machine learning, a 1-bit and multi-bits quantization methods have been developed in \citep{seide20141, bernstein2018icml,alistarh2017qsgd,magnusson2019maintaining}. Other variants of quantized gradient schemes include error compensation \citep{wu2018error}, variance-reduced quantization \citep{zhang2017zipml}, and quantization to a ternary vector \citep{wen2017terngrad}.

\noindent\textbf{Sparsification} amounts to transmitting only gradient coordinates with large enough magnitudes exceeding a certain threshold~\citep{strom2015scalable, aji2017sparse}.
To avoid losing information of skipping communication, small gradient components will be accumulated and then transmitted when they are large enough \citep{lin2017deep,stich2018nips,alistarh2018}. Recently, randomized sparsification approaches have also been developed in \citep{wangni2018gradient,wang2018atomo}.

However, both quantization and sparsification aim to resolve c1).
For exchanging messages, e.g., the $p$-dimensional $\theta$ or its gradient, other latencies (initiating communication links, queueing, and propagating the message) are at least comparable to the message size-dependent transmission latency \citep{peterson2007}. This motivates c2) reducing the number of communication rounds.

\noindent\textbf{Periodic communication.}
In contrast to the gradient compression schemes, schemes that reduce the number of communication rounds have also been developed, including the periodic averaging techniques, e.g., local stochastic gradient descent (a.k.a. local SGD) \citep{lin2018don,stich2019local,wang2018coop,yu2019icml,yu2019aaai}.
In local SGD, workers are allowed to perform local model updates independently and the resultant models are averaged periodically. In this way communication frequency is reduced.
The caveat is that most local SGD methods have performance guarantee in the \emph{homogeneous} settings, where the data are independent and identically distributed over all workers. However, this assumption rarely holds in federated learning \citep{mcmahan2017}.

\noindent\textbf{Intermittent communication.}
Different from periodic communication used in local SGD, adaptive uploading techniques have been studied in e.g., lazily aggregated gradient (LAG) approaches \citep{chen2018lag,sun2019}.
LAG is tailored for the \emph{heterogeneous} learning settings, and has provable performance gain when the data distributions vary across workers.
Models in LAG are updated at the server, and workers only adaptively upload information that is determined to be informative enough. Unfortunately, while the original LAG has good performance in the deterministic settings (e.g., with full gradient), its performance is significantly degraded in the stochastic settings. Recent efforts have been made towards adaptive uploading in stochastic settings \citep{li2019csgd}, but the proposed scheme therein requires an exponentially increasing batch size, which is not favorable in practice. In contrast, the LASG approaches can be viewed as the stochastic counterparts of the LAG, and our adaptive communication rules are new and tailored for SGD, which do not require the increasing batch size.

\section{Why LAG does not work well with SGD?}
Our LASG approaches are closely related to the recently developed LAG method \citep{chen2018lag}.
In this section, we revisit LAG and provide insights why it does not work well in stochastic settings.

\begin{figure}[t]
\centering
\includegraphics[width=.5\textwidth]{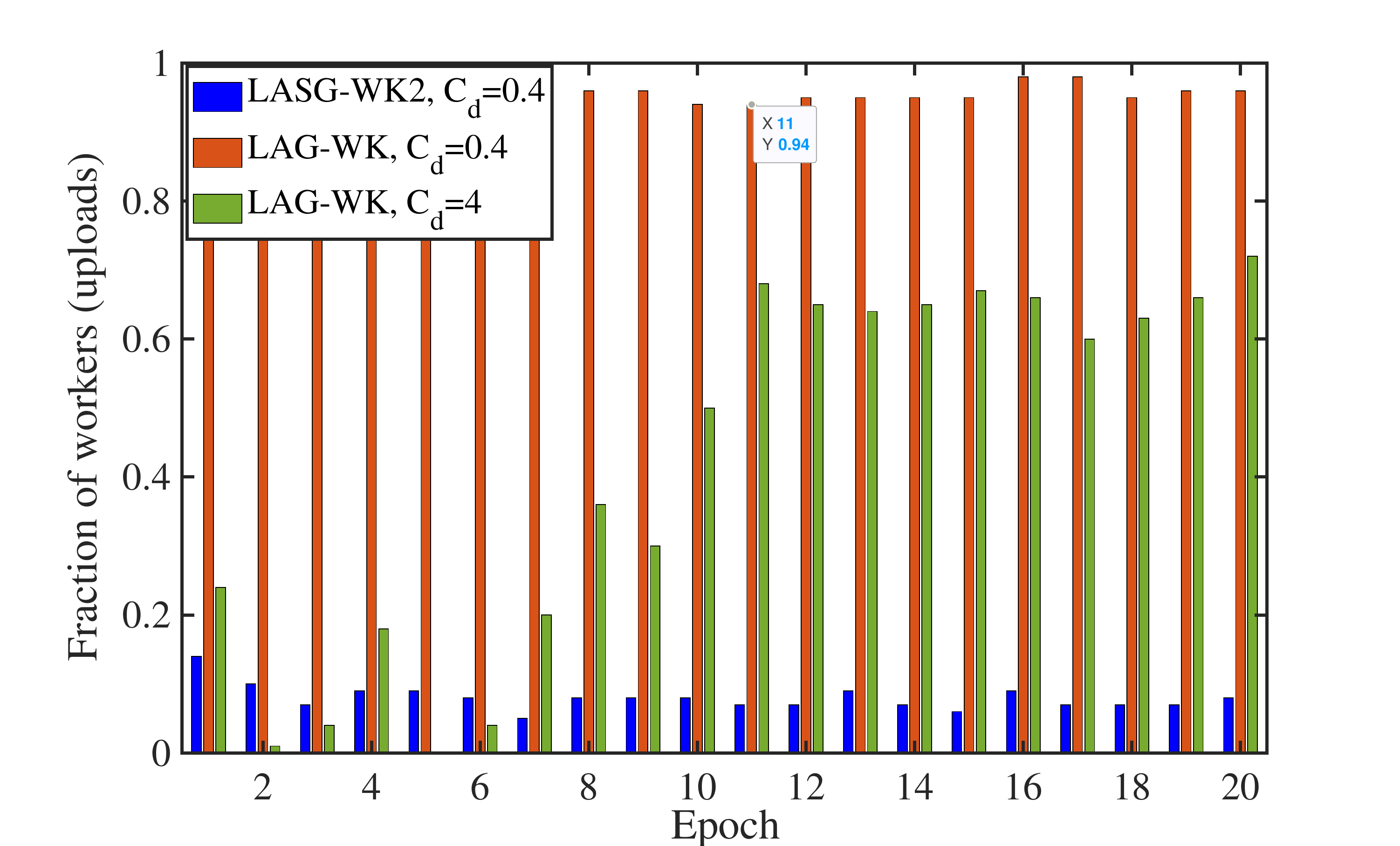}
  \caption{Number of uploads in per epoch (10 iterations) under stochastic LAG-WK and LASG-WK2.}
  \label{toyexample}
\end{figure}


Because not every communication round is equally important during the learning process, LAG only admits useful communication, and otherwise, reuses stale information.
Instead of communicating with all $M$ workers as SGD in \eqref{eqn:SGDiteration}, the direct (or ``naive") stochastic version of LAG (specifically LAG-WK) selects the subset of workers ${\cal M}^k$ to obtain their \emph{fresh} stochastic gradients $\nabla{\cal L}_m\big(\theta^k;\xi_m^k\big)$, $m\in{\cal M}^k$.
The \textbf{direct stochastic LAG} also follows the generic update \eqref{eqn:LASGiteration}, but it selects ${\cal M}^k$ as: if worker $m$ finds the innovation of the fresh stochastic gradient $\nabla \ell(\theta^k;\xi_m^k)$ is small such that it satisfies (with pre-defined $\{c_d\geq0\}$)
{\begin{align}\label{eqn:LAGrule}
\Big\|\nabla \ell(\theta^k;\xi_m^k)-\nabla \ell(\theta^{k-\tau_m^k};\xi_m^{k-\tau_m^k})\Big\|^2
\!\leq\frac{1}{M^2}\sum\limits_{d=1}^Dc_d\big\|\theta^{k+1-d}-\theta^{k-d}\big\|^2\! 
\end{align}}
then we reuse the old gradient, $m\in{\cal M}\backslash{\cal M}^k$, and increase the staleness by $\tau_m^{k+1} = \tau_m^{k} + 1$; otherwise, worker $m$ uploads the stochastic gradient, and resets $\tau_m^{k+1}=1$.

In the deterministic setting, LAG condition \eqref{eqn:LAGrule} is motivated by the elegant \emph{``larger descent per upload"} rationale, and has proved to be effective \citep{chen2018lag}.
Nevertheless, the observation here is that the two stochastic gradients \eqref{eqn:LAGrule} are evaluated on two different iterates ($\theta^k$ and $\theta^{k-\tau_m^k}$) and two different samples ($\xi_m^k$ and $\xi_m^{k-\tau_m^k}$) thus two different loss functions. This is in contrast to the original LAG condition in \citep{chen2018lag} where the gradient innovation is evaluated on the same function. This subtle difference leads to significant degradation in performance.

Figure \ref{toyexample} compares the stochastic LAG and one of our new algorithms LASG-WK2 (introduced later), and demonstrates that the stochastic LAG is not effective in saving communication --- when $c_d$ is set to be small (e.g., 0.4), \eqref{eqn:LAGrule} almost never satisfies; and when $c_d$ is set to be large (e.g., 4), \eqref{eqn:LAGrule} satisfies initially, but stops satisfying later. This can be explained by expanding the left-hand-side (LHS) of \eqref{eqn:LAGrule} by (see the supplemental material for the deduction)
\begin{subequations}\label{eqn:variance-wk}
\begin{align}
\EE\left[\|\nabla\ell(\theta^k;\xi_m^k)-\nabla\ell(\theta^{k-\tau_m^k};\xi_m^{k-\tau_m^k})\|^2\right]
\geq \,&\frac{1}{2}\EE\Big[\big\|\nabla\ell(\theta^k;\xi_m^k)-\nabla{\cal L}_m(\theta^k)\big\|^2\Big]\label{eqn:wk-rq1}\\
+\,&\frac{1}{2}\EE\Big[\big[\big\|\nabla\ell(\theta^{k-\tau_m^k};\xi_m^{k-\tau_m^k})-\nabla{\cal L}_m(\theta^{k-\tau_m^k})\big\|^2\big]\Big]\label{eqn:wk-rq2}\\
-\,&\EE[\|\nabla{\cal L}_m(\theta^k)-\nabla{\cal L}_m(\theta^{k-\tau_m^k})\|^2]. \label{eqn:wk-rq3}
\end{align}
\end{subequations}
Even if the iterate $\theta^k$ converges, e.g., $\theta^k\rightarrow \theta^*$, and thus the right-hand-side (RHS) of \eqref{eqn:LAGrule} $\big\|\theta^{k+1-d}\!-\!\theta^{k-d}\big\|^2\!\rightarrow\! 0$, the LHS of \eqref{eqn:LAGrule} does not, because the gradient variance appearing in \eqref{eqn:wk-rq1} and \eqref{eqn:wk-rq2} does not vanish yet the gradient difference at the same function \eqref{eqn:wk-rq3} diminishes.

Therefore, the key insight here is that the non-diminishing variance of stochastic gradients makes the direct implementation of the LAG rule \eqref{eqn:LAGrule} ineffective eventually.

 \begin{table*}[t]
%
\centering
 \begin{tabular}[t]{ c || c | c ||c |c|| c |c  }
 \hline\hline \textbf{Metric} &\multicolumn{2}{|c||}{\textbf{Communication}}& \multicolumn{2}{|c||}{\textbf{Computation}}& \multicolumn{2}{|c}{\textbf{Memory}}\\
\hline \!\textbf{Algorithm}\!&\textbf{PS$\rightarrow$WK}~$m$\!\!&\!\textbf{WK} $m$ $\rightarrow$\textbf{PS}\!& ~~~\textbf{PS}~~~&\textbf{WK}~$m$& ~~~\textbf{PS}~~~&\textbf{WK}~$m$\\ \hline \hline
\textbf{Sync SGD} & always &  always & \eqref{eqn:SGDiteration} & \eqref{eqn:SGDiteration} &  ${\cal O}(p)$& $/$ \\ \hline
     \textbf{LASG-WK1} & always &\!\!\!only if $m\!\in\!{\cal M}^k$\!\!\!  & \!\!\! \eqref{eqn:LASGiteration2}\!\! & \!\! \eqref{eqn:workerrule1} &  \!\!\!\! ${\cal O}(p)$\!\!\!& \!\!\!\!${\cal O}(p)$\!\!\!\!  \\\hline
    \textbf{LASG-WK2} & always & \!\!\!only if $m\!\in\!{\cal M}^k$\!\!\! & \eqref{eqn:LASGiteration2}  &\!\!\!$\eqref{eqn:workerrule2}$\!\!\!&  ${\cal O}(p)$&\!\!\!\!\! ${\cal O}(p)$\!\!\!\!\\ \hline
    \textbf{LASG-PS} & only if $m\!\in\!{\cal M}^k$ & \!\!\!only if $m\!\in\!{\cal M}^k$\!\!\! & \eqref{eqn:LASGiteration2},\eqref{eqn:serverrule}  & \!\!only if $m\!\in\!{\cal M}^k$ \!\!\!&  ${\cal O}(Mp)$&\!\!\!${\cal O}(p)$\!\!\\ \hline
        \textbf{LASG-PSE} & only if $m\!\in\!{\cal M}^k$ & \!\!\!only if $m\!\in\!{\cal M}^k$\!\!\! & \eqref{eqn:LASGiteration2},\eqref{eqn:eserversule}  & \!\!only if $m\!\in\!{\cal M}^k$ \!\!\!&  ${\cal O}(Mp)$ &\!\!\!${\cal O}(p)$\!\!\!\\ \hline\hline
    \end{tabular}
\caption{A comparison of communication, computation and memory requirements.
\textbf{PS} denotes the parameter server, \textbf{WK} denotes the worker,
\textbf{PS$\rightarrow$WK $m$} is the download from the server to the worker $m$, and \textbf{WK} $m$ $\rightarrow$ \textbf{PS} is the upload from the worker $m$ to the server.}
\label{tab:3algcomp}
\end{table*}

\section{LASG: Lazily Aggregated Stochastic Gradient Approach}\label{sec.LASG}
In this section, we formally develop our LASG method, and present the intuition
behind its design.


While the updates of stochastic LAG and LASG \eqref{eqn:LASGiteration} look identical, the choice of ${\cal M}^k$ in them is very different.
To overcome the limitations of LAG in stochastic settings, the key of the LASG design is to \textbf{reduce the variance of the innovation measure} appeared in the adaptive condition.

Towards this goal, we develop two types of LASG rules to select ${\cal M}^k$.
The first type is adopted by each worker that uses the gradient difference as the innovation measure but a \emph{variance-reduced} gradient difference; and the second one by the parameter server that uses the model difference as the innovation measure, but uses a sequence of \emph{diminishing stepsizes} to control variance.

%

\begin{table*}[t]
\centering
    \begin{tabular}{c c}
    \hspace{-0.2cm}
\begin{minipage}[t]{7.8cm}
  \vspace{0pt}
\begin{algorithm}[H]
\caption{LASG-WK1}\label{alg: LASG-WK1}
    \begin{algorithmic}[1]
    \State{\textbf{Input:} Delay counter $\{\tau_m^0\}$,  stepsizes $\{\eta_k\}$, max delay $D$.}
        \For{$k=0,1,\ldots, K-1$}
            \State{Server broadcasts $\theta^k$ to all workers.}
            \State{All workers save $\ttheta=\theta^k$ if $k \mod D=0$.}
            \ParFor{ Worker $m=1,2,\ldots, M$}
            \State{Compute $\nabla \ell(\theta^k;\xi_m^k)$ and $\nabla \ell(\ttheta;\xi_m^k)$.}
            \State{Check condition \eqref{eqn:workerrule1} with stored $\tdelta_m^{k-\tau_m^{k}}$.}
                \If{\eqref{eqn:workerrule1} is violated, or, $k\mod D=0$}
                \State{Upload $\delta_m^k$.}
               \qquad\quad\Comment{$\tau_m^{k+1}=1$}
                \Else
                \State{Upload nothing.}\quad\Comment{$\tau_m^{k+1} = \tau_m^{k} + 1$}
                \EndIf
            \EndParFor
            \State{Server updates via \eqref{eqn:LASGiteration}.}
        \EndFor
    \end{algorithmic}
  \end{algorithm}
\end{minipage}
\begin{minipage}[t]{7.8cm}
  \vspace{0pt}
  \begin{algorithm}[H]
\caption{LASG-WK2}
\label{alg: LASG-WK2}
    \begin{algorithmic}[1]
        \State{\textbf{Input:} Delay counter $\{\tau_m^0\}$,  stepsizes $\{\eta_k\}$, max delay $D$.}
        \For{$k=0,1,\ldots, K-1$}
            \State{Server broadcasts $\theta^k$ to all workers.}
            \ParFor{ Worker $m=1,2,\ldots, M$}
            \State{Compute $\nabla \ell(\theta^k;\xi_m^k)$ and $\nabla \ell(\theta^{k-\tau_m^k}_m;\xi_m^k)$.}
            \State{Check condition \eqref{eqn:workerrule2}.}
                \If{\eqref{eqn:workerrule2} is violated, or, $\tau_m^k\geq D$}
                \State{Upload $\delta_m^k$. \qquad\qquad~~~~\Comment{$\tau_m^{k+1}=1$}}
                \Else
                \State{Upload nothing.}\qquad~\Comment{$\tau_m^{k+1} = \tau_m^{k} + 1$}
                \EndIf
            \EndParFor
            \State{Server updates via \eqref{eqn:LASGiteration}.}
        \EndFor
    \end{algorithmic}
\end{algorithm}
\end{minipage}
   \end{tabular}
   \caption{A comparison of LASG-WK1 and LASG-WK2.}
\end{table*}

\subsection{Worker LASG: save communication uploads}\label{subsec.wk}
We first introduce two LASG variants that use variance-reduced rules to check gradient innovation at the worker side.
The first one that we term \textbf{LASG-WK1} will reuse the old gradient of worker $m$ at iteration $k$ if it satisfies
\begin{equation}\label{eqn:workerrule1}
\left\|\tdelta_m^k-\tdelta_m^{k-\tau_m^k}\right\|^2\leq\frac{1}{M^2}\sum\limits_{d=1}^Dc_d\left\|\theta^{k+1-d}-\theta^{k-d}\right\|^2
\end{equation}
where $\tdelta_m^k := \nabla \ell(\theta^k;\xi_m^k)-\nabla \ell(\ttheta;\xi_m^k)$ is the stochastic gradient difference at a common sample $\xi_m^k$, $\tdelta_m^{k-\tau_m^k} := \nabla \ell(\theta^{k-\tau_m^k};\xi_m^{k-\tau_m^k})-\nabla \ell(\ttheta;\xi_m^{k-\tau_m^k})$ is the stochastic gradient difference at a common sample $\xi_m^{k-\tau_m^k}$,
and $\ttheta$ is a snapshot of the previous iterate that will be updated every $D$ iterations.
If \eqref{eqn:workerrule1} is satisfied, the staleness increases by $\tau_m^{k+1} = \tau_m^{k} + 1$; otherwise, worker $m$ uploads the fresh stochastic gradient, and resets staleness as $\tau_m^{k+1}=1$.

The rationale of \eqref{eqn:workerrule1} follows next.
In contrast to the non-vanishing variance in the LAG-WK rule (see \eqref{eqn:variance-wk}),
the LASG-WK1 rule \eqref{eqn:workerrule1} reduces its inherent variance.
To see this, we can decompose the LHS of \eqref{eqn:workerrule1} as the difference of two \emph{variance reduced} stochastic gradients at iteration $k$ and $k-\tau_m^k$.
Using the stochastic gradient in SVRG as an example \citep{johnson2013accelerating},
the innovation can be written as
\begin{align*}
\tdelta_m^k-\tdelta_m^{k-\tau_m^k}=\left(\nabla\ell(\theta^k;\xi_m^k)-\nabla\ell(\tilde{\theta};\xi_m^k)+\nabla{\cal L}_m(\tilde{\theta})\right)
-\left(\nabla\ell(\theta^{k-\tau_m^k};\xi_m^{k-\tau_m^k})-\nabla\ell(\tilde{\theta};\xi_m^{k-\tau_m^k})+\nabla{\cal L}_m(\tilde{\theta})\right).
\end{align*}
Define the minimizer of \eqref{eqn: problem} as $\theta^{\star}$ and assume that $\nabla\ell(\theta;\xi_m)$ is $\bar{L}$-Lipschitz continuous for any $\xi_m$.
The expectation of the LHS of \eqref{eqn:workerrule1} can be \emph{upper-bounded} by
\begin{align}\label{eqn:variance-wk1}
\EE\left[\big\|\tdelta_m^k-\tdelta_m^{k-\tau_m^k}\big\|^2\right]\leq8\bar{L}(\EE{\cal L}(\theta^k)-{\cal L}(\theta^{\star}))
+8\bar{L}(\EE{\cal L}(\theta^{k-\tau_m^k})-{\cal L}(\theta^{\star}))\!+\!16\bar{L}(\EE{\cal L}(\ttheta)-{\cal L}(\theta^{\star})).
\end{align}
If the iterate $\theta^k$ converges, e.g., $\theta^k, \theta^{k-\tau_m^k}, \ttheta\rightarrow \theta^*$, the RHS of \eqref{eqn:variance-wk1} diminishes, and thus the LHS of \eqref{eqn:workerrule1} diminishes. This is in contrast to the stochastic LAG-WK rule in \eqref{eqn:variance-wk} that is \emph{lower-bounded} by a non-diminishing value.

In addition to \eqref{eqn:workerrule1}, the second rule that we term \textbf{LASG-WK2} reuses the old gradient of worker $m$ if it satisfies
\begin{equation}\label{eqn:workerrule2}
\!\left\|\nabla \ell(\theta^k;\xi_m^k)-\!\nabla \ell(\theta^{k-\tau_m^k}_m\!;\xi_m^k)\right\|^2\!\leq\!\frac{1}{M^2}\!\sum\limits_{d=1}^D\!c_d\left\|\theta^{k+1-d}\!\!-\theta^{k-d}\right\|^2\!\!.
\end{equation}
If \eqref{eqn:workerrule2} is satisfied, the server will use the stale stochastic gradient $\nabla \ell(\theta^{k-\tau_m^k}_m;\xi_m^{k-\tau_m^k})$ for worker $m$, and the staleness increases by $\tau_m^{k+1} = \tau_m^{k} + 1$; otherwise, worker $m$ uploads the fresh stochastic gradient, and resets the staleness as $\tau_m^{k+1}=1$. Notice that different from the naive LAG-WK \eqref{eqn:LAGrule}, the LASG condition \eqref{eqn:workerrule2} is evaluated at two different iterates but on the same sample $\xi_m^k$.

Similar to LASG-WK1, the LASG-WK2 rule \eqref{eqn:workerrule2} also reduces its inherent variance, since the LHS of \eqref{eqn:workerrule2} can be written as the difference between a \emph{variance reduced} stochastic gradient and a \emph{deterministic} gradient, that is
\begin{align}\label{eqn:variance-wk2}
    \nabla\ell(\theta^k;\xi_m^k)-\nabla\ell(\theta^{k-\tau_m^k};\xi_m^k)=\big(\nabla\ell(\theta^k;\xi_m^k)
    -\nabla\ell(\theta^{k-\tau_m^k};\xi_m^k)+\nabla{\cal L}_m(\theta^{k-\tau_m^k})\big)-\nabla{\cal L}_m(\theta^{k-\tau_m^k}).
\end{align}
With derivations deferred to the supplementary document, similar to \eqref{eqn:variance-wk1} we can also conclude that
$\EE[\|\nabla\ell(\theta^k;\xi_m^k)-\nabla\ell(\theta^{k-\tau_m^k};\xi_m^k)\|^2]\rightarrow 0$ as the iterate $\theta^k\rightarrow \theta^{\star}$.

\subsection{Server LASG: save up/downloads and calculations}\label{subsec.ps}
Besides the worker-side conditions, we next introduce two LASG variants that use variance-reduced rules to check model innovation at the server side, both of which do not even need to broadcast the current models. The rationale is that if the model difference is small, the gradient difference used in Section \ref{subsec.wk} is likely to be small.

The first one that we term \textbf{LASG-PS} will reuse the old gradient of worker $m$, given that the old parameter that worker $m$ used for computing the last stochastic gradient satisfies
\begin{equation}\label{eqn:serverrule}
L_m^2\left\|\theta^{k}-\theta^{k-\tau_m^k}\right\|^2\leq\frac{1}{M^2}\sum\limits_{d=1}^Dc_d\left\|\theta^{k+1-d}-\theta^{k-d}\right\|^2
\end{equation}
where $L_m$ is the smoothness constant of $\cL_m(\theta)$. Condition \eqref{eqn:serverrule} can be checked at the server side without computing new gradients if the server stores $\{\theta^{k-\tau_m^k}\}$ that all workers used for computing the most recent stochastic gradients.

The LHS of \eqref{eqn:serverrule} can be \emph{upper-bounded} in expectation by
\begin{align}\label{eqn:variance-server}
\EE\left[\big\|\theta^{k}-\theta^{k-\tau_m^k}\big\|^2\right]&\leq 2D \sum_{d=1}^{D}\EE\left[\big\|\theta^{k-d}-\theta^{k-d-\tau_m^{k-d}}\big\|^2\right]\eta_{k-D}^2 \nonumber\\
&+2D \sum_{d=1}^{D}\mathbb{E}\left\|\nabla \cL(\theta^{k-d})\right\|^2\eta_{k-D}^2+D^2 \left(\sum_{m\in {\cal M}}\sigma_m^2\right)\eta_{k-D}^2.
\end{align}
If the iterate $\theta^k$ does not diverge so that $\left\|\nabla \cL(\theta^k)\right\|^2$ is bounded, then the diminishing stepsizes $\{\eta_k\}$ ensure that the second and third terms in the RHS of \eqref{eqn:variance-server} vanish. Using mathematical induction, the LHS of \eqref{eqn:serverrule} also diminishes. Therefore, similar to the variance-reduced gradient difference used in LASG-WK, the diminishing stepsizes can also make the LASG-PS condition effective asymptotically.

\begin{table*}[t]
\centering
\vspace{-0.9cm}
\begin{tabular}{c c}
\hspace{-0.5cm}
\begin{minipage}[t]{7.8cm}
\begin{algorithm}[H]
\caption{LASG-PS}
\label{alg: LASG-PS}
    \begin{algorithmic}[1]
  \State{\textbf{Input:} $\theta^0$, delay counter $\{\tau_m^0\}$, smoothness contants $\{L_m\}$, stepsizes $\{\eta_k\}$, maximum delay $D$.} 
        \For{$k=0,1,\ldots, K-1$}
            \ParFor{ Worker $m=1,2,\ldots, M$}
            \State{Server checks condition \eqref{eqn:serverrule}.}
                \If{\eqref{eqn:serverrule} is violated or $\tau_m^k\geq D$} %
                \State{Server sends $\theta^k$ to worker $m$} 
                \State{Worker $m$ computes $\nabla \ell(\theta^k;\xi_m^k)$.}
                \State{Worker $m$ uploads $\delta_m^k$.} \quad~~\Comment{$\tau_m^{k+1}=1$}
                \Else
                \State{No action.} \qquad\qquad~~\Comment{$\tau_m^{k+1} = \tau_m^{k} + 1$}
                \EndIf
            \EndParFor
            \State{Server updates via \eqref{eqn:LASGiteration}.}
        \EndFor
    \end{algorithmic}
\end{algorithm}
\end{minipage}
\begin{minipage}[t]{7.8cm}
\begin{algorithm}[H]
\caption{LASG-PSE}
\label{alg: LASG-PSE}
    \begin{algorithmic}[1]
  \State{\textbf{Input:} $\theta^0$, delay counter $\{\tau_m^0\}$, smoothness estimates $\{\hat{L}_m^0\}$, stepsizes $\{\eta_k\}$, maximum delay $D$.} 
        \For{$k=0,1,\ldots, K-1$}
            \ParFor{ Worker $m=1,2,\ldots, M$}
            \State{Server checks condition \eqref{eqn:eserversule}.}
                \If{\eqref{eqn:eserversule} is violated or $\tau_m^k\geq D$} %
                \State{Server sends $\theta^k$ to worker $m$.} 
                \State{Worker $m$ computes $\nabla \ell(\theta^k;\xi_m^k)$.}
                \State{Worker $m$ uploads $\delta_m^k$.} \quad~~\Comment{$\tau_m^{k+1}=1$}
                \State{Worker $m$ uploads $\hat{L}_m^{k+1}$ in \eqref{eqn:Lmestimate}.}
                \Else
                \State{No action.} \qquad\qquad~~\Comment{$\tau_m^{k+1} = \tau_m^{k} + 1$}
                \EndIf
            \EndParFor
            \State{Server updates via \eqref{eqn:LASGiteration}.}
        \EndFor
    \end{algorithmic}
\end{algorithm}
\end{minipage}
\end{tabular}
\vspace{0.1cm}
\caption{A comparison of LASG-PS and LASG-PSE.}
\end{table*}

In many problems, however, $L_m$ is not always available, or, hard to compute. To resolve this issue, we develop LASG-PSE, a variation of LASG-PS that estimates $L_m$ ``on-the-fly.''
With $\hat{L}_m^k$ denoting the estimate of $L_m$, \textbf{LASG-PSE} will reuse the old gradient of worker $m$ if it satisfies
\begin{equation}\label{eqn:eserversule}
 (\hat{L}_m^k)^2\|\theta^{k}-\theta^{k-\tau_m^k}\|^2\leq\frac{1}{M^2}\sum\limits_{d=1}^Dc_d\|\theta^{k+1-d}-\theta^{k-d}\|^2
\end{equation}
where the estimated constant $\hat{L}_m^k$ is updated iteratively via
\begin{equation}\label{eqn:Lmestimate}
\!\!\hat{L}_m^{k+1}=\max\left\{\hat{L}_m^k, \frac{\|\nabla\ell(\theta^{k};\xi_m^k)-\nabla\ell(\theta^{k-\tau_m^k};\xi_m^k)\|}{\|\theta^{k}-\theta^{k-\tau_m^k}\|}\right\}.
\end{equation}

We summarize LASG-PS and LASG-PSE in Algorithms \ref{alg: LASG-PS} and \ref{alg: LASG-PSE}, and compare our four LASG variants in Table \ref{tab:3algcomp}.

\noindent\textbf{Comparison of all LASG variants.}
Comparing WK conditions with PS conditions in Table \ref{tab:3algcomp}, LASG-PS and LASG-PSE need \emph{extra memory} at the server side but save both \emph{local computation and download communication}, while LASG-WK1 and LASG-WK2 save only upload communication.
Between the two WK conditions, LASG-WK1 is more conservative as LASG-WK1 measures the change of gradients at two model states for both new and old data samples but LASG-WK2 measures only the change of gradient at the new sample.
Between the two PS conditions, LASG-PSE is more flexible since it does not require the knowledge of local smoothness $L_m$.
Given specific communication, computation, and memory requirements, we can flexibly choose different LASG.

%
%
%

\subsection{Quantized LASG: Further save communication bits}
The four LASG variants save the number of communication rounds.
To further reduce communication bits per round, we combine LASG with various quantization mechanisms. With the stochastic gradient $\nabla \ell(\theta;\xi)$, we define the gradient under a quantization operator ${\cal Q}$ as
\begin{equation}
Q(\theta;\xi)	:={\cal Q}\left(\nabla \ell(\theta;\xi)\right).
\end{equation}

We adopt the stochastic quantization scheme in \citep{alistarh2017qsgd} and develop the quantized LASG as
\begin{align*}
\theta^{k+1} \!=\! \theta^k\!-\!\eta_k\!\!\!\!\!\sum\limits_{m\in{\cal M}\backslash{\cal M}^k} \!\!\!\!\! Q(\theta^{k-\tau_m^k};\xi_m^{k-\tau_m^k})-\eta_k\!\!\!\!\sum\limits_{m\in{\cal M}^k}\!\!\!Q(\theta^k;\xi_m^k)
\end{align*}
where ${\cal M}^k$ is determined by one out of four rules \eqref{eqn:workerrule1}-\eqref{eqn:eserversule}.
We term the quantized LASG as LAQSG.

\section{Convergence Analysis}
In this section we present the convergence results of LASG-WK1, LASG-WK2 and LASG-PS in both the nonconvex and strongly convex cases, and the convergence results of LAQSG in the nonconvex case.
Due to the technical reasons, we leave the analysis of LASG-PSE for future work, but it empirically has very impressive performance.

First, we make some basic assumptions, which are standard in analyzing SGD and its variants \citep{ghadimi2013stochastic, johnson2013accelerating,alistarh2017qsgd}.
\begin{assumption}\label{assump:smoothness}
The loss function $\cL$ is smooth with $L>0$.
\end{assumption}
\begin{assumption}\label{assump:gradientestimator}
The samples $\xi_m^1,\xi_m^2,\ldots$ are independent, and the stochastic gradient $\nabla\ell(\theta;\xi_m^k)$ satisfies
\begin{align}
    &\EE_{\xi_m^k} \nabla\ell(\theta;\xi_m^k) = \nabla\cL_m(\theta),\\
    &\EE_{\xi_m^k} \left[\|\nabla\ell(\theta;\xi_m^k)-\nabla\cL_m(\theta)\|^2\right]\leq\sigma_m^2. \label{eqn:boundedvariance}
\end{align}
\end{assumption}
For LASG-PS, we require an extra smoothness assumption.
\begin{assumption}\label{assump:localsmoothness}
The local gradient $\nabla\cL_m$ is $L_m$-Lipschitz continuous, i.e. for any $\theta_1,\theta_2\in\RR^p$, we have
\begin{equation}
    \|\nabla\cL_m(\theta_1)-\nabla\cL_m(\theta_2)\|\leq L_m\|\theta_1-\theta_2\|.
\end{equation}
\end{assumption}

With these assumptions, LASG will yield descent of $\cL(\theta^k)$.
\begin{lemma}\label{lemma:lossdescent}
Under Assumptions \ref{assump:smoothness}, \ref{assump:gradientestimator} and \ref{assump:localsmoothness}, $\{\theta^k\}$ generated by Algorithms \ref{alg: LASG-WK1}, \ref{alg: LASG-WK2} and \ref{alg: LASG-PS} satisfy
\begin{align}\label{eqn:lossdescent}
    \EE\cL(\theta^{k+1})&\leq\EE\cL(\theta^k)-\left(\eta_k-L\eta_k^2\right)\EE\left[\|\nabla\cL(\theta^k)\|^2\right]\nonumber\\
    &+\sum\limits_{d=1}^D\Big(\big(\eta_k+\frac{1}{2L}\big)c_d+\frac{\sqrt{M}L}{12}\Big)\EE\left[\|\theta^{k+1-d}-\theta^{k-d}\|^2\right]+L\eta_k^2\left(\frac{9}{2}+6\sqrt{M}D\right)\sum\limits_{m=1}^M\sigma_m^2.
\end{align}
\end{lemma}
\vspace{-0.15cm}
Note that all the terms on the right hand side of the inequality \eqref{eqn:lossdescent} show up in SGD analysis except $\|\theta^{k+1-d}\!-\!\theta^{k-d}\|^2$, which exists due to stale information.
To deal with this term, we introduce the following Lyapunov function:
\begin{equation}
\label{eqn:Lyapunov}
    V^k:=\EE\cL(\theta^k)-\cL(\theta^{\star})+\sum\limits_{d=1}^{D}\gamma_d\EE\left[\|\theta^{k+1-d}-\theta^{k-d}\|^2\right]
\end{equation}
where $\{\gamma_d\}_{d=1}^D$ are constants to be determined later. The following lemma is a direct application of Lemma \ref{lemma:lossdescent}.

\begin{lemma}\label{lemma:lyapunovdescent}
Under Assumptions \ref{assump:smoothness} and \ref{assump:gradientestimator}, there exist constants $\{A_d^k\}_{d=1}^D $, $B_0^k$ and $B_1^k$ such that 
\begin{align}\label{descent}
    \EE V^{k+1}-\EE V^k\leq -B_0^k\EE\left[\|\nabla\cL(\theta^{k})\|^2\right]+B_1^k\sum\limits_{m=1}^M\sigma_m^2-\sum\limits_{d=1}^DA_d^k\EE\left[\|\theta^{k+1-d}-\theta^{k-d}\|^2\right].
\end{align}
\end{lemma}
The constants $\{A_d^k\}_{d=1}^D$, $B_0^k$ and $B_1^k$ depend on stepsize $\eta_k$, $\gamma$ and $\{c_d\}_{d=1}^D$. Their expressions are specified in the supplementary materials. By choosing proper $\eta_k$ and $c_d$, we are able to ensure the convergence of LASG.

We first present the convergence in nonconvex case.
\begin{theorem}\label{thm:nonconvex}
{\normalfont \textsf{(nonconvex)}} Under Assumptions \ref{assump:smoothness}, \ref{assump:gradientestimator} (for Algorithm \ref{alg: LASG-PS} also Assumption \ref{assump:localsmoothness}), if $\eta_k=\eta={\cal O}(\frac{1}{\sqrt{K}})$ with $c_{\eta}>0$, and $c_d\leq\min\{\frac{1}{12D\eta^2},\frac{\sqrt{M}L^2}{18}\}$, then $\{\theta^k\}$ generated by Algorithms \ref{alg: LASG-WK1}-\ref{alg: LASG-PS} satisfy 
\begin{equation}
    \frac{1}{K}\sum\limits_{k=0}^{K-1}\EE\left[\|\nabla\cL(\theta^k)\|^2\right]=O\left({1}/{\sqrt{K}}\right).
\end{equation}
\end{theorem}

Next we present the convergence results under the following strong convexity assumption on ${\cal L}(\theta)$.
\begin{assumption}\label{assump:strongconvexity}
The overall loss $\cL(\theta)$ is $\mu$-strongly convex.
\end{assumption}
Parallel to the sublinear convergence of SGD in the
strongly convex case, e.g., \citep{rakhlin2011making}, LASG algorithms achieve the $O(1/K)$ order of convergence.
\begin{theorem}\label{thm:stronglyconvex}
{\normalfont \textsf{(strongly convex)}} Under Assumption \ref{assump:smoothness},\ref{assump:gradientestimator},\ref{assump:strongconvexity} (for Algorithm \ref{alg: LASG-PS} only, also Assumption \ref{assump:localsmoothness}), if $\eta_k=\frac{2}{\mu(k+K_0)}$ for a given constant $K_0$, and $c_d\leq\min\{\frac{1}{24D\eta_0^2},\frac{\sqrt{M}L^2}{18}\}$, then $\theta^K$ generated by Algorithms \ref{alg: LASG-WK1}, \ref{alg: LASG-WK2} and \ref{alg: LASG-PS} satisfies
\begin{equation}
    \EE\cL(\theta^K)-\cL(\theta^{\star})=O\left({1}/{K}\right).
\end{equation}
\end{theorem}

For the convergence of LAQSG algorithms, we make the following additional assumption that guarantees the bounded variance of the quantized stochastic gradient.
\begin{assumption}\label{assump:boundedgradient}
The gradient is bounded as $\EE_{\xi_m}\left[\|\nabla\ell(\theta;\!\xi_m)\|^2\right]\!\leq\! B$ for any $\theta \in \RR^p$.
\end{assumption}
Based on this assumption, we have the following result.
\begin{theorem}\label{thm:LAQSG}
{\normalfont \textsf{(LAQSG)}} Under Assumptions \ref{assump:smoothness}, \ref{assump:gradientestimator}, \ref{assump:boundedgradient} (also Assumption \ref{assump:localsmoothness} for Algorithm \ref{alg: LASG-PS}), if $\eta_k\!=\!\eta\!={\cal O}(\frac{1}{\sqrt{K}})$, $c_d\leq\min\{\frac{1}{16D\eta^2},\frac{\sqrt{M}L^2}{24}\}$ where $c_{\eta}>0$ is a constant, then $\{\theta^k\}$ generated by quantized Algorithms \ref{alg: LASG-WK1} - \ref{alg: LASG-PS} satisfy
\begin{equation}
    \frac{1}{K}\sum\limits_{k=0}^{K-1}\EE\left[\|\nabla \cL(\theta^{k})\|^2\right]=O\left({1}/{\sqrt{K}}\right).
\end{equation}
\end{theorem}

\begin{figure*}[t]
\centering
 \hspace*{-3ex}
    \includegraphics[width=.35\textwidth]{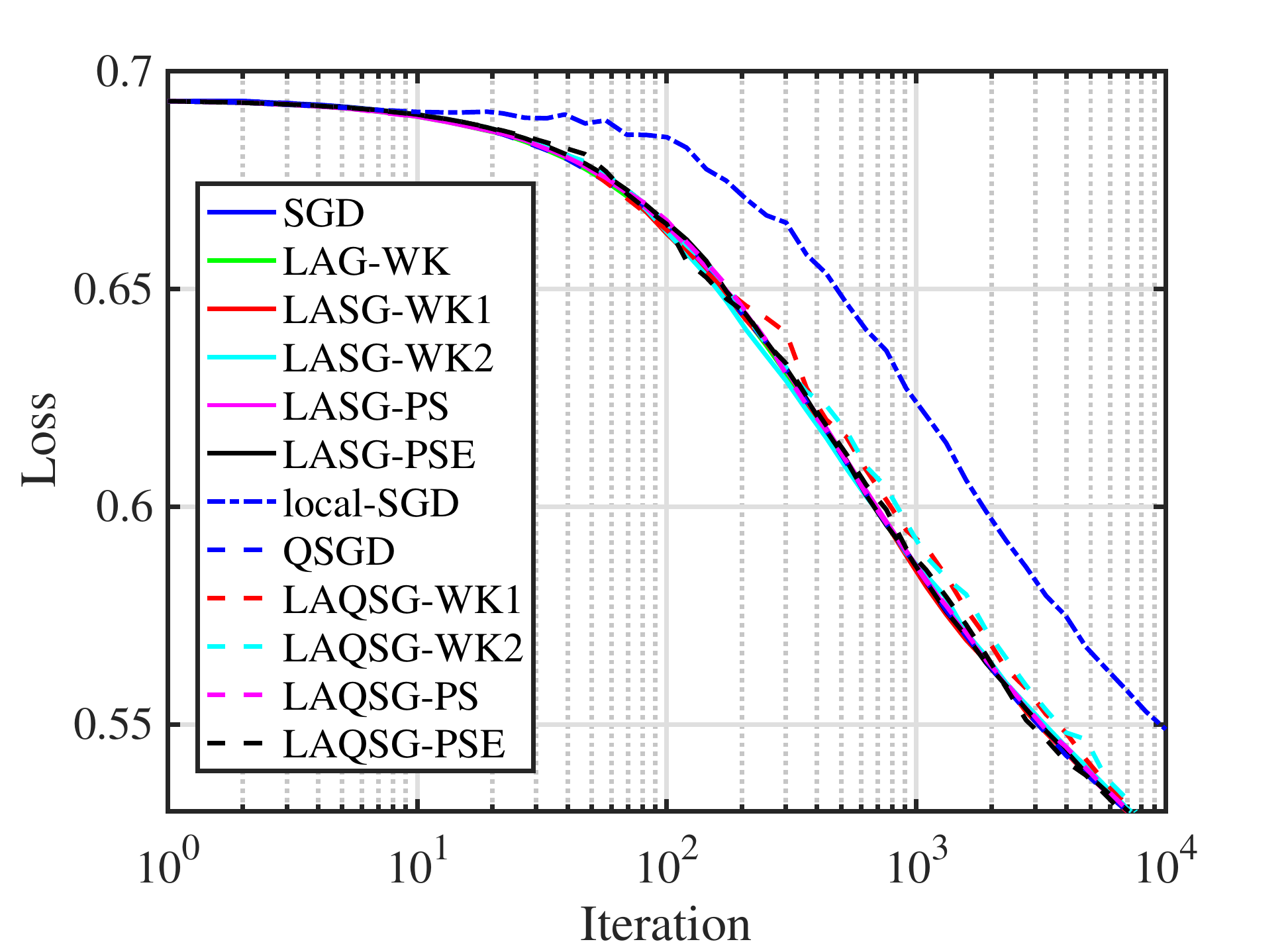}
    \hspace*{-4ex}
    \includegraphics[width=.35\textwidth]{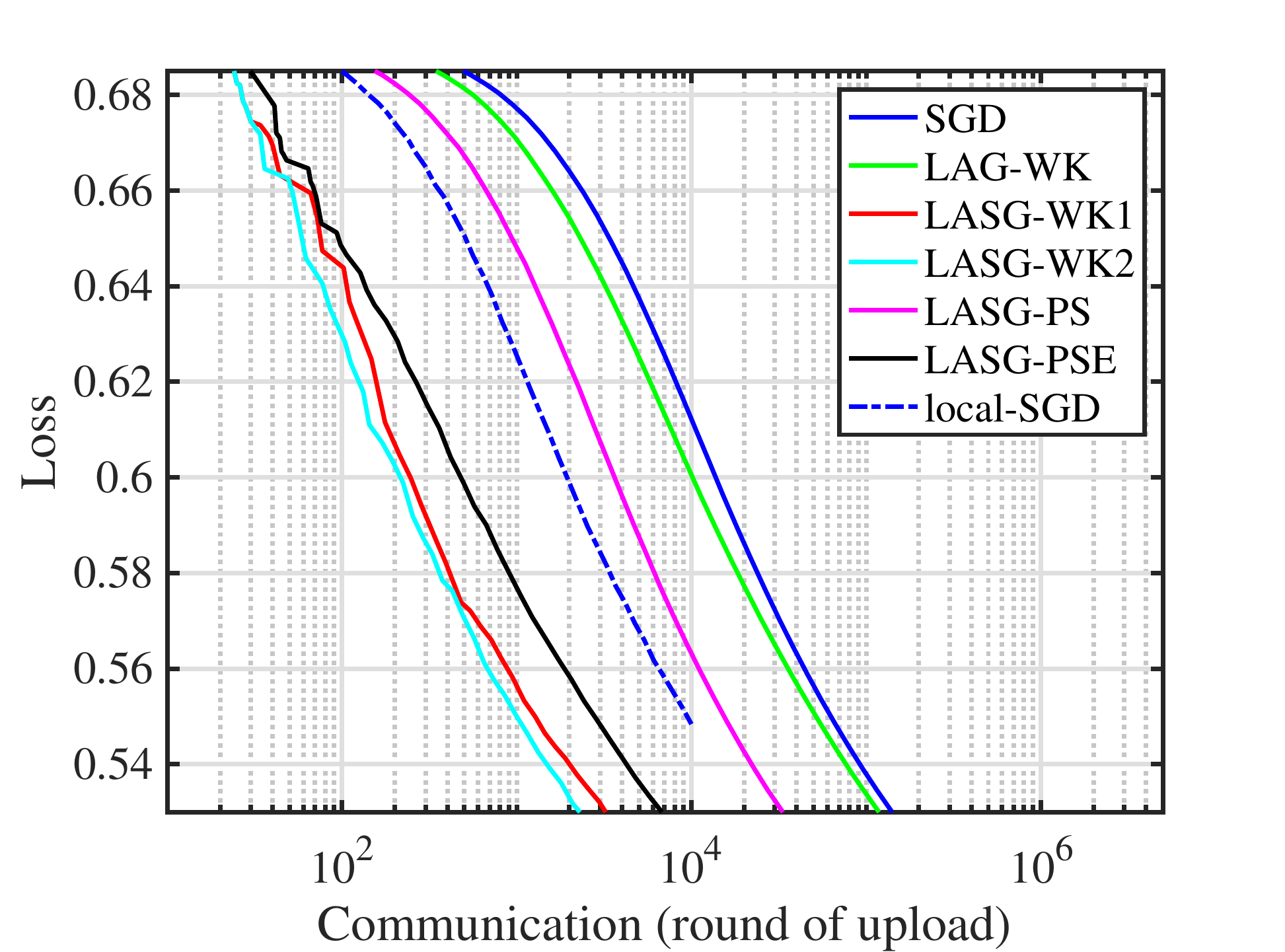}
    \hspace*{-4ex}
    \includegraphics[width=.35\textwidth]{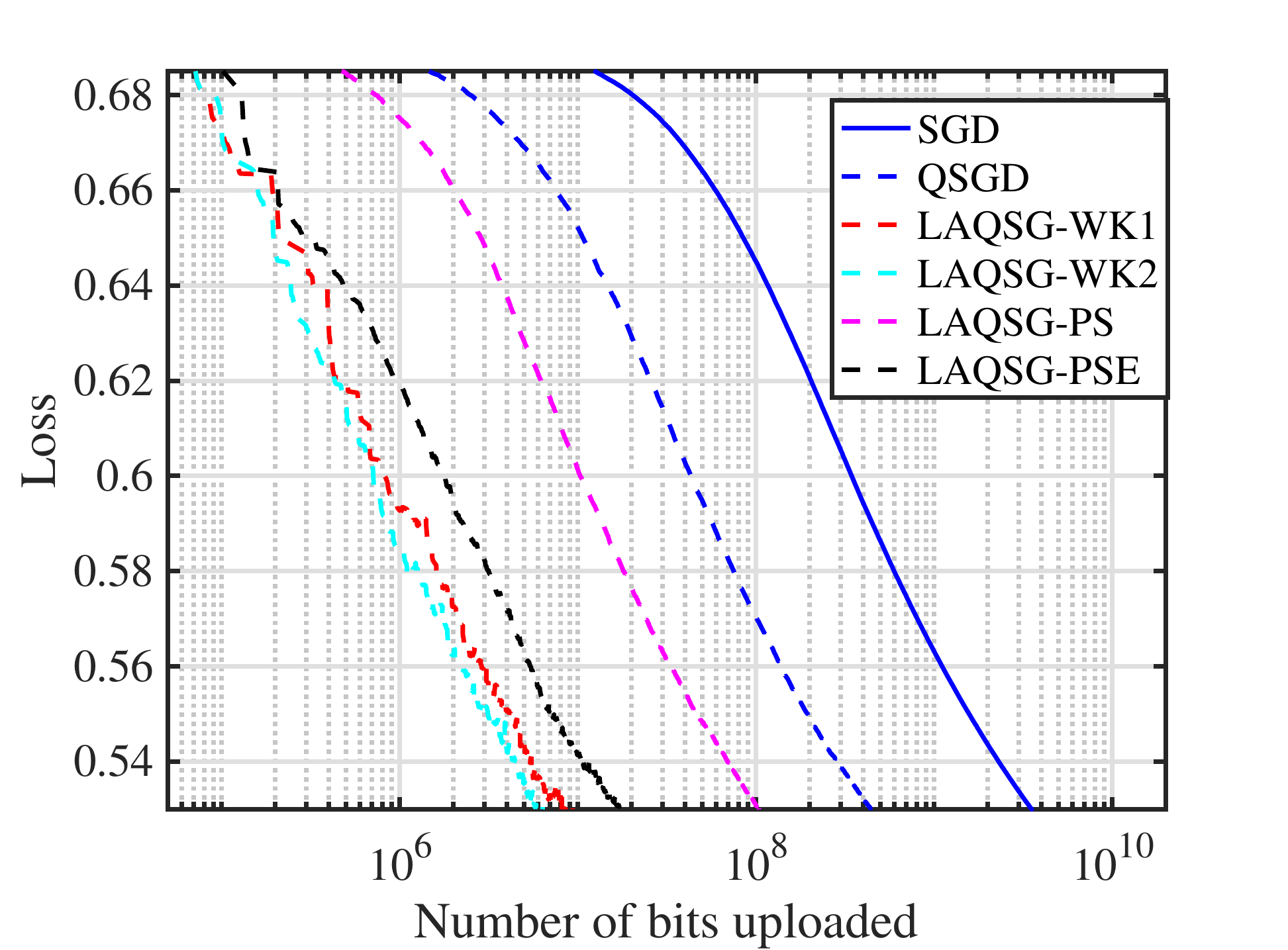}
    \caption{Logistic regression on \textit{covtype} dataset in the heterogeneous setting}
    \label{fig:covtype}
\end{figure*}

\begin{figure*}[t]
\centering
 \hspace*{-3ex}
    \includegraphics[width=.35\textwidth]{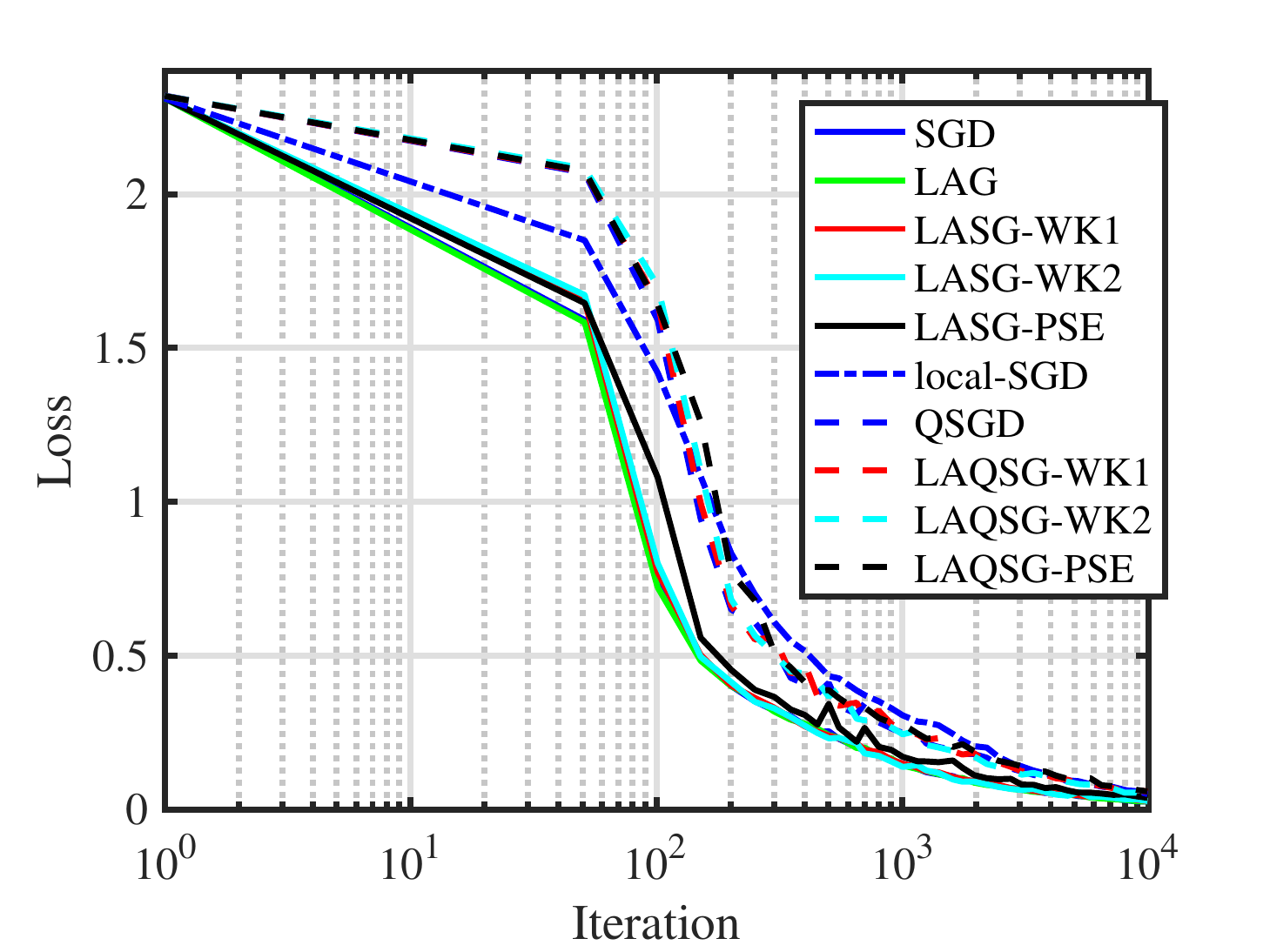}
    \hspace*{-4ex}
    \includegraphics[width=.35\textwidth]{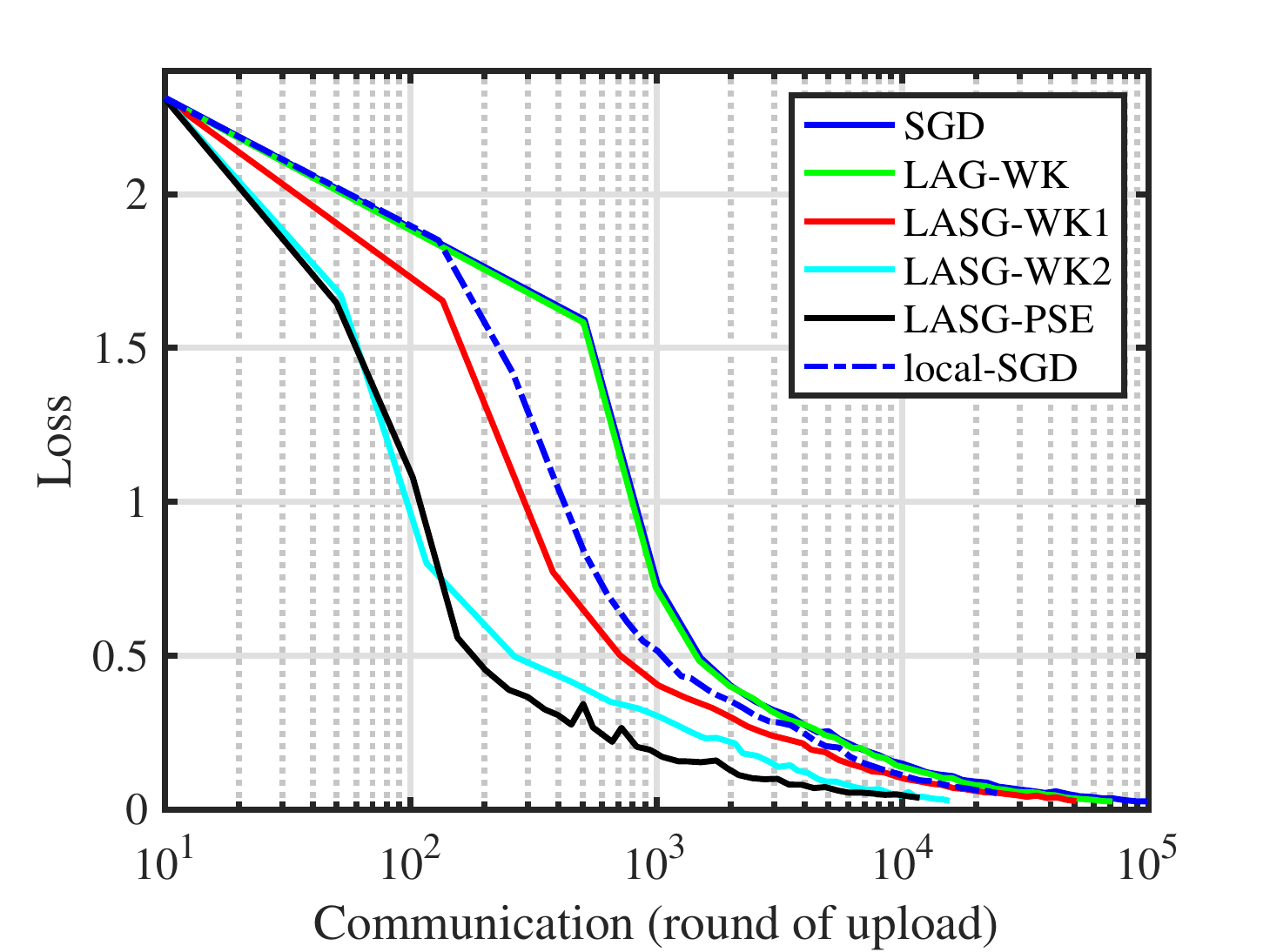}
    \hspace*{-4ex}
    \includegraphics[width=.35\textwidth]{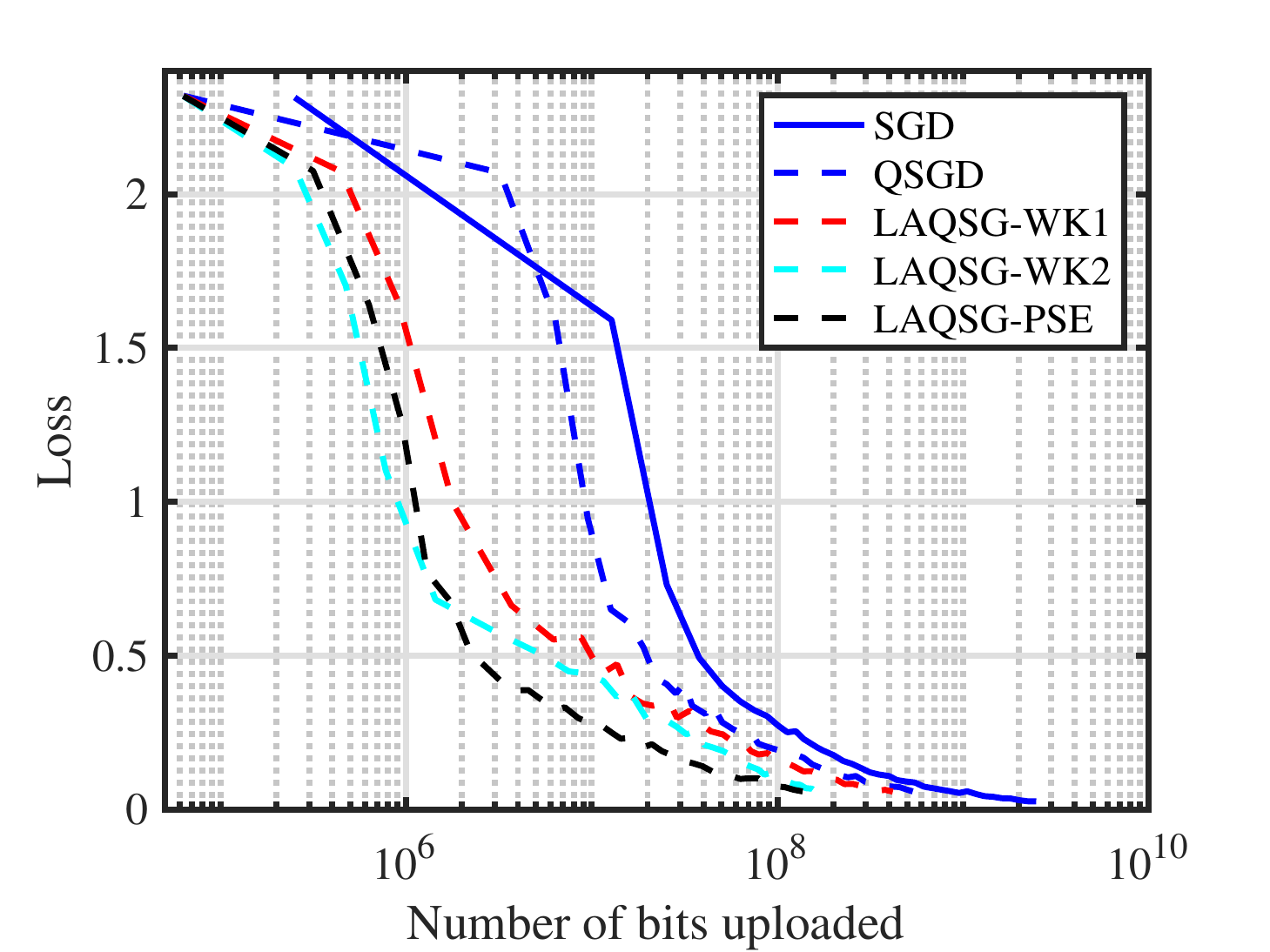}
    \caption{Training Neural network on \textit{mnist} dataset in the heterogeneous setting.}
    \label{fig:NNmnist}
\end{figure*}

\begin{figure*}[t]
\centering
 \begin{subfigure}[b]{0.5\textwidth}
\centering
\includegraphics[width=1\textwidth]{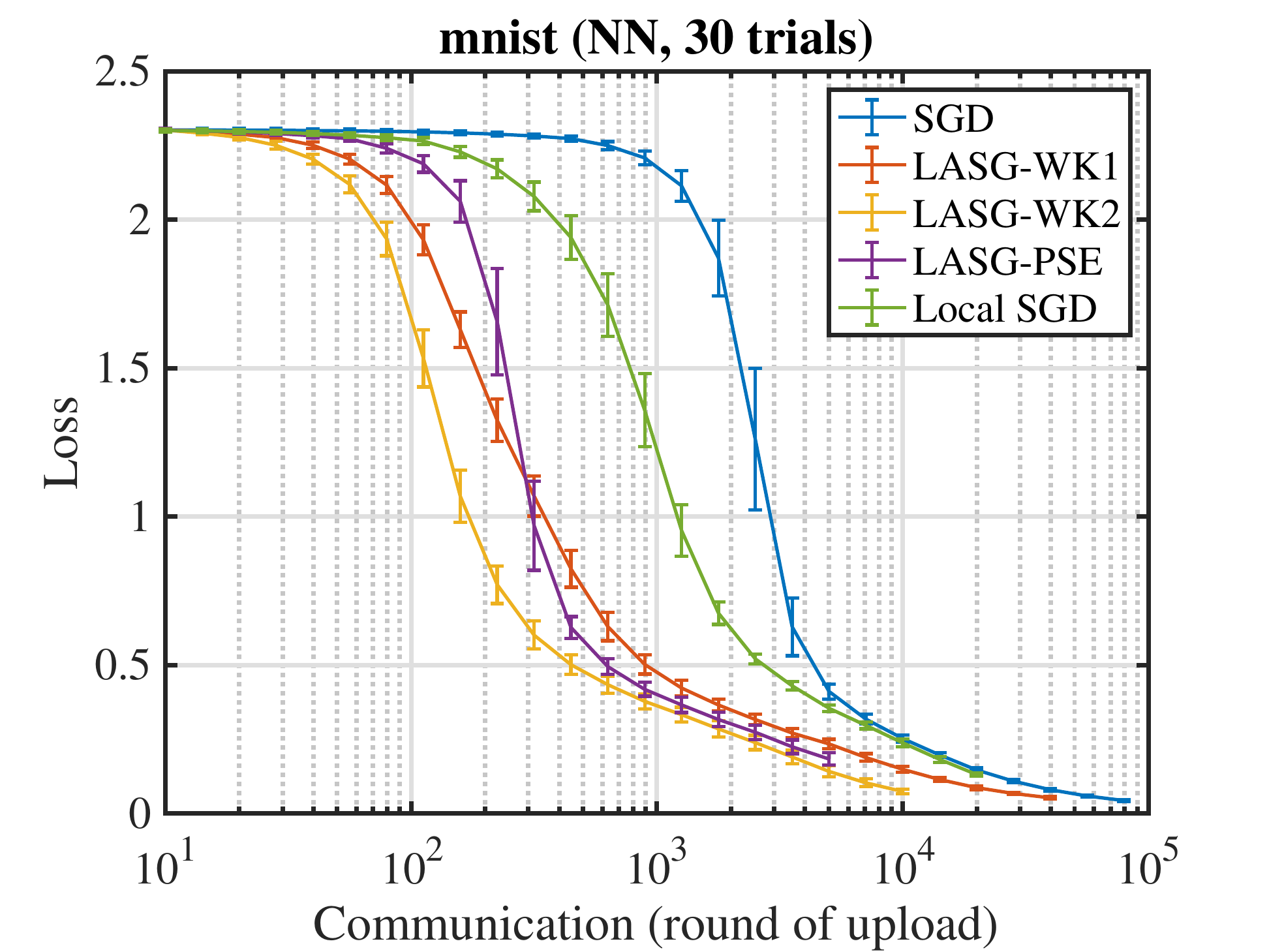}
\caption{Training loss on \textit{mnist} dataset.}
    \label{fig:mnist_train_seed}
    \end{subfigure}%
 \begin{subfigure}[b]{0.5\textwidth}
\centering
\includegraphics[width=1\textwidth]{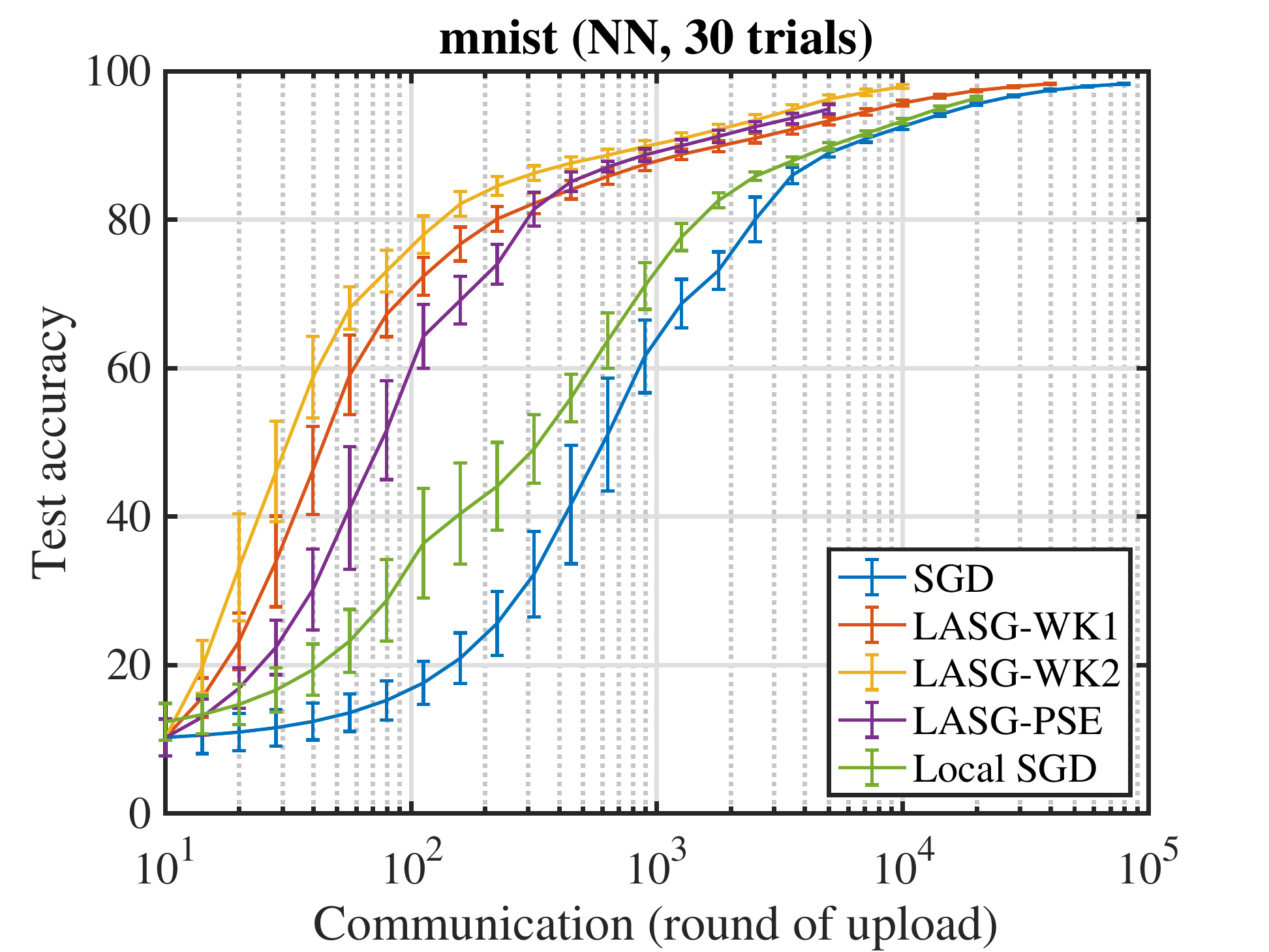}
\caption{Training loss on \textit{mnist} dataset.}
    \label{fig:mnist_test_seed}
    \end{subfigure}%
 \caption{Simulations on \textit{mnist} dataset averaged over 30 trials.}   
\end{figure*}

\begin{figure*}[t]
\centering
 \begin{subfigure}[b]{0.5\textwidth}
\centering
\includegraphics[width=1\textwidth]{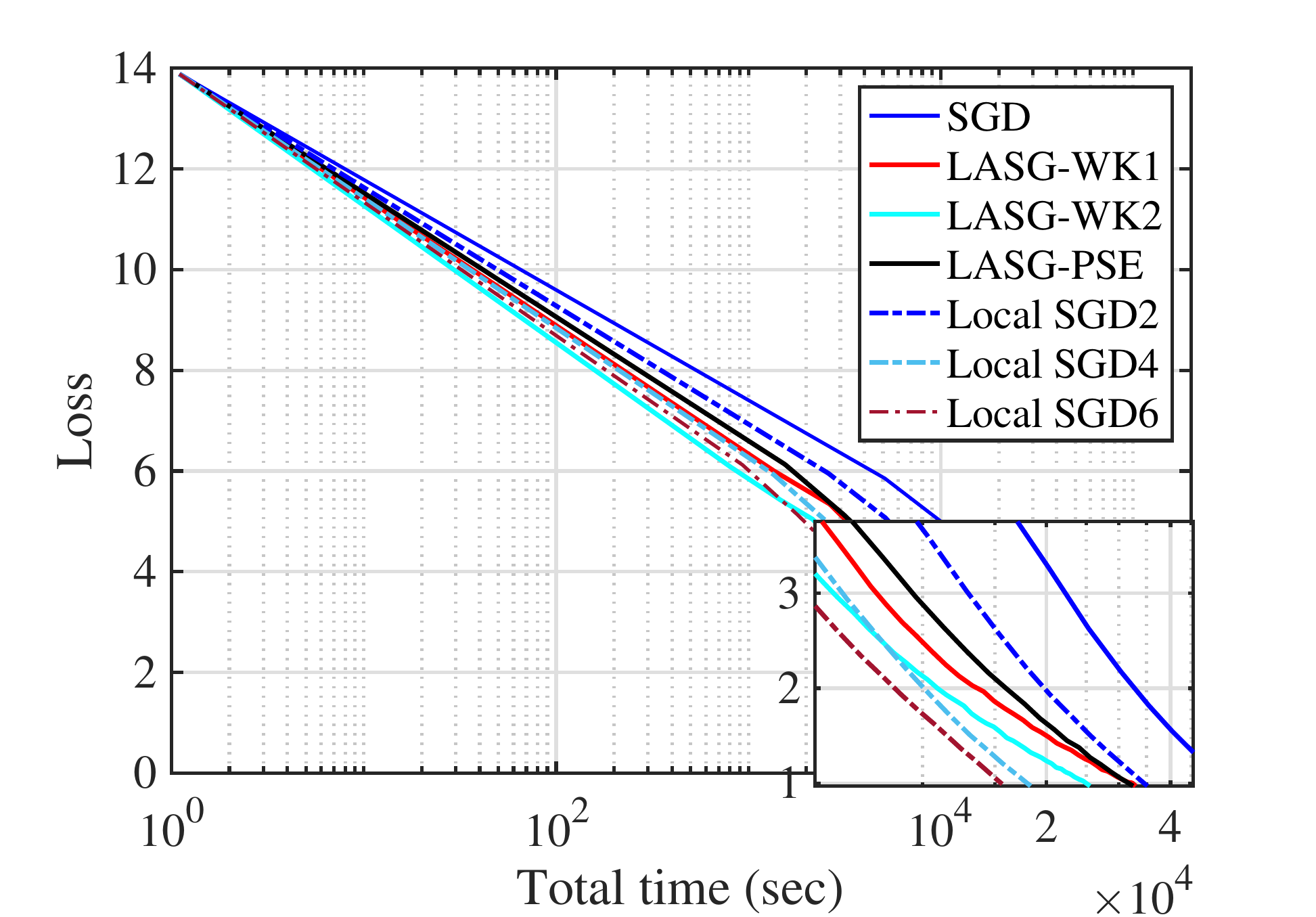}
\caption{Training loss on \textit{tiny imagenet} dataset.}
    \label{fig:imagenet_train}
    \end{subfigure}%
 \begin{subfigure}[b]{0.5\textwidth}
\centering
\includegraphics[width=1\textwidth]{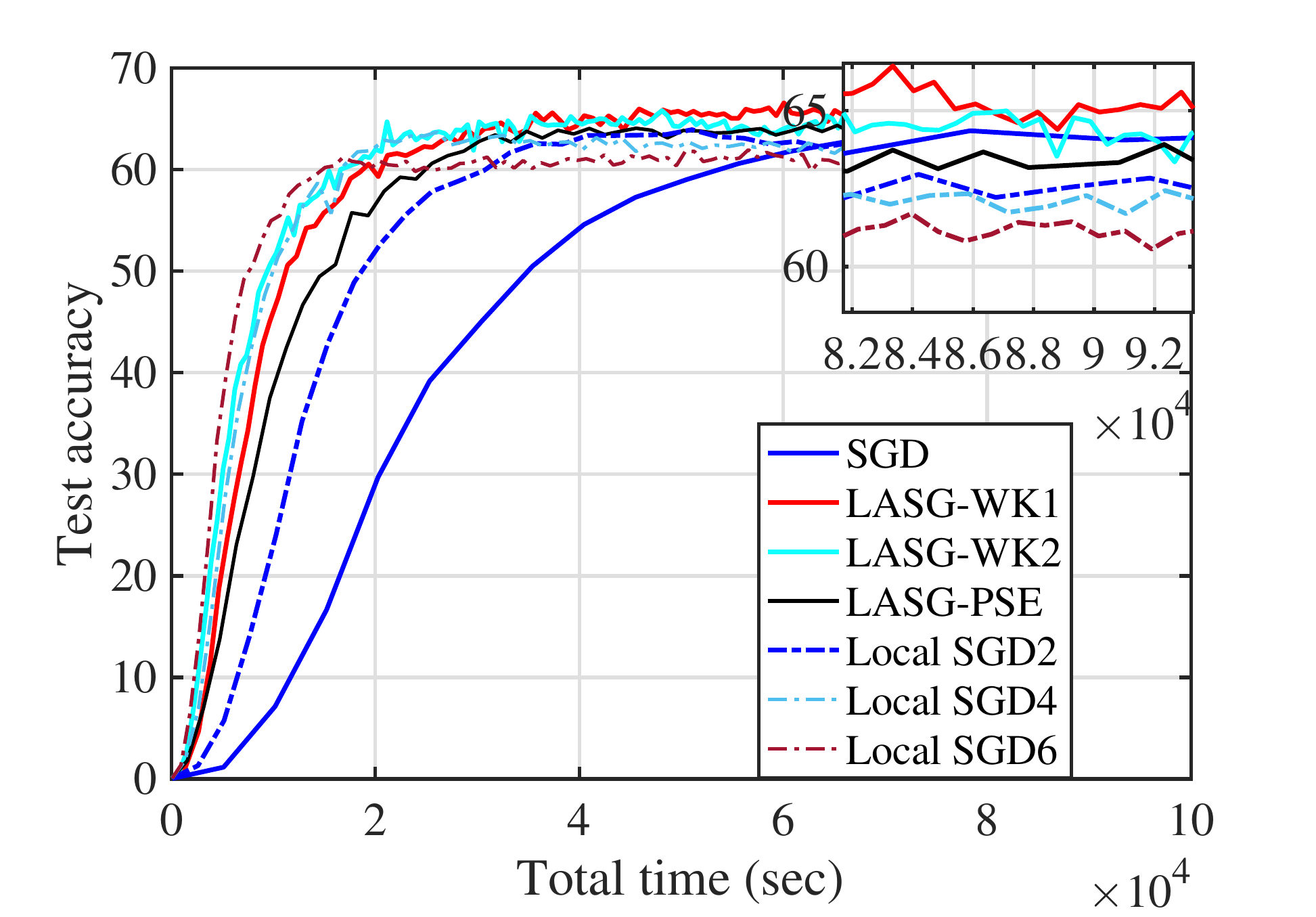}
\caption{Testing accuracy on \textit{tiny imagenet} dataset.}
    \label{fig:imagenet_test}
    \end{subfigure}%
 \caption{Simulations on \textit{tiny imagenet} dataset.}   
\end{figure*}

\begin{table*}[t]
 \vspace{-0.2cm}
\centering
\resizebox{\textwidth}{!}{
\begin{tabular}{c||c|c|c|c|c|c|c||c|c|c|c|c}
\hline
\multirow{2}{*}{Algorithm} & \multirow{2}{*}{SGD} & \multirow{2}{*}{\!\!\!LAG-WK\!\!\!} & \multicolumn{4}{|c|}{LASG} &\multirow{2}{*}{\!\!\!local SGD\!\!\!} &\multirow{2}{*}{QSGD} & \multicolumn{4}{|c}{LAQSG}\\ \cline{4-7}\cline{10-13}
                     &                      &                         & WK1 & WK2 & PS & PSE
&                           &                      & WK1 & WK2 & PS & PSE \\ \hline\hline
ijcnn1(LR)&0.4383& 0.3317 & {\bf 0.2311} & {\bf 0.2308} & 0.2736  & 0.2356   & 0.2554 & 0.4382 & {\bf 0.2307} & 0.2316    & 0.2722  & 0.2356\\ \hline
mnist(LR) &0.4214& 0.4074 & 0.1826 & {\bf 0.1281} & 0.2950  & 0.1681 & 0.2189 & 0.4228 & {\bf 0.1486} & {\bf 0.1342}    & 0.2839  & 0.1572\\ \hline
covtype(LR)&0.6772&0.6727 & 0.5604 &{\bf 0.5523}& 0.6467 & 0.5985 & 0.6288 & 0.6772 &{\bf 0.5537}&{\bf 0.5477}& 0.6456 & 0.5914 \\ \hline
mnist(NN)&0.1544&0.1457&0.1063 &{\bf 0.0510} & ---& {\bf 0.0421} & 0.1178 & 0.2492 & 0.2327 & 0.1564 &--- &{\bf 0.1046}\\ \hline
\end{tabular}}
\caption{Loss after $10^3$ and $10^4$ communication rounds for logistic regression (LR) and neural network (NN) in heterogeneous setting.}
\label{table: 1}
\end{table*}

\begin{table*}[t]
    \vspace{0.2cm}
\centering
\resizebox{\textwidth}{!}{
\begin{tabular}{c||c|c|c|c|c|c|c||c|c|c|c|c}
\hline
\multirow{2}{*}{Algorithm} & \multirow{2}{*}{SGD} & \multirow{2}{*}{\!\!\!LAG-WK\!\!\!} & \multicolumn{4}{|c|}{LASG} &\multirow{2}{*}{\!\!\!local SGD\!\!\!} &\multirow{2}{*}{QSGD} & \multicolumn{4}{|c}{LAQSG}\\ \cline{4-7}\cline{10-13}
                     &                      &                         & WK1 & WK2 & PS & PSE
&                           &                      & WK1 & WK2 & PS & PSE \\ \hline\hline
ijcnn1(LR)&0.6356& 0.3672 & 0.3050 & 0.3037 & 0.3378 & 0.3053 & 0.3226 & 0.4597 & {\bf 0.2370} & {\bf 0.2363} & 0.2798 & {\bf 0.2411}\\ \hline
mnist(LR)&0.6989& 0.6762 & 0.5081 & 0.3859 & 0.6653 & 0.4029 & 0.6125 & 0.5690 & {\bf 0.1916} & {\bf 0.1681} & 0.4215 & {\bf 0.2071} \\ \hline
covtype(LR)&0.6837&0.6806 & 0.5770 & 0.5650 & 0.6620 & 0.6178 & 0.6483 & 0.6469 & {\bf 0.5273} & {\bf 0.5277} & 0.5930 & {\bf 0.5488} \\ \hline
mnist(NN)&0.3042&0.2792&0.2157 & 0.1159 &---&{\bf 0.0814}& 0.2472 & 0.2011 & 0.1602 &{\bf 0.1104} &---&{\bf 0.0764} \\ \hline
\end{tabular}}
\caption{Loss after $10^5$, $10^6$, $10^6$ bits of uploads for LR on ijcnn1, mnist, covtype,  and $10^8$ for NN in heterogeneous setting.}
\label{table: 2}
\end{table*}

\section{Simulations}
Numerical tests have been conducted on both logistic regression and neural network models.

We benchmark LA(Q)SG with SGD, LAG-WK, local SGD and QSGD. For local SGD \citep{lin2018don}, workers perform SGD independently to update local $\{\theta_m\}_{m=1}^M$, and $\{\theta_m\}_{m=1}^M$ then are averaged over all workers every $H$ iterations. In simulations, we did a grid search for SGD learning rates.
We consider the heterogeneous setting where data with same labels are unevenly assigned to $M$ workers.

\noindent\textbf{Logistic regression on ijcnn1, MNIST and covtype.}
The data are distributed across $M=10$ workers for ijcnn1, MNIST (with digits 3, 5) and $M=20$ for Covtype. For each worker, the batch size is selected to be 0.01 of the local data size for ijcnn1, MNIST and  0.001 for Covtype. The $\ell_2$-regularization parameter is set to be $10^{-5}$. We choose stepsize $\eta=0.1$. For all LASG algorithms, $D=100$ and $c_d=0.1/\eta^2/M^2$ for $d=1,2,\ldots 10$ and $c_d=0$ for $d=11,\ldots 100$. For local-SGD, the communication period is $H=50, 10, 20$ iterations for ijcnn1, MNIST, Covtype respectively. This is optimized to save communication as much as possible without largely affecting the convergence speed. For quantization methods, we perform 4-bit stochastic quantization \citep{alistarh2017qsgd}. Numerical results are reported in Figure \ref{fig:covtype} and in Tables \ref{table: 1}, \ref{table: 2}. 
Performance averaged over multiple trails has also been reported in Figures \ref{fig:mnist_train_seed} and \ref{fig:mnist_test_seed}.
Supplementary materials have additional tests in Figures \ref{fig:ijcnn1_u}-\ref{fig:covtype_u} and Tables \ref{table: 3} and \ref{table: 4}.

\noindent\textbf{Neural network.}
We train a convolutional neural network with two convolution-ELU-maxpooling layers (ELU is a smoothed ReLU) followed by two fully-connected layers for 10 classes classification on MNIST. The data are distributed on $M=10$ workers. We choose stepsize $\eta=0.05$. Since the objective function is nonsmooth ($L_m$ is not available), LASG-PS is not considered in this test. For all LASG algorithms, we set $D=50$, $c_d=0.1/\eta^2/M^2$ for $d=1,2,\ldots 10$ and $c_d=0$ for $d=11,\ldots 50$.
For local-SGD, we set the communication period to be 4. For all quantization methods, we perform 8-bit stochastic quantization. Numerical results are reported in Figure \ref{fig:NNmnist} and listed in Tables \ref{table: 1} and \ref{table: 2}. Additional results can be found in Figures \ref{fig:NNmnist_u}-\ref{fig:NNmnist_test} and Tables \ref{table: 3}, \ref{table: 4} in the supplementary materials.

All algorithms have been tested on the popular \emph{tiny imagenet} dataset using the Resnet18 model initialized by weights pretrained on ImageNet1000; see the accuracy versus total time (communication and computation) in Figures \ref{fig:imagenet_train} and \ref{fig:imagenet_test}. For training loss, LASG-WK1 and -WK2 require much less total time than SGD and local SGD with $H=2$, but slightly more than local SGD with $H=4$ and $6$. However, as shown in Figure \ref{fig:imagenet_test}, local SGD with larger communication period sacrifices the testing accuracy by 3-4\%.

\begin{figure}[h]
\centering
\includegraphics[width=.45\textwidth]{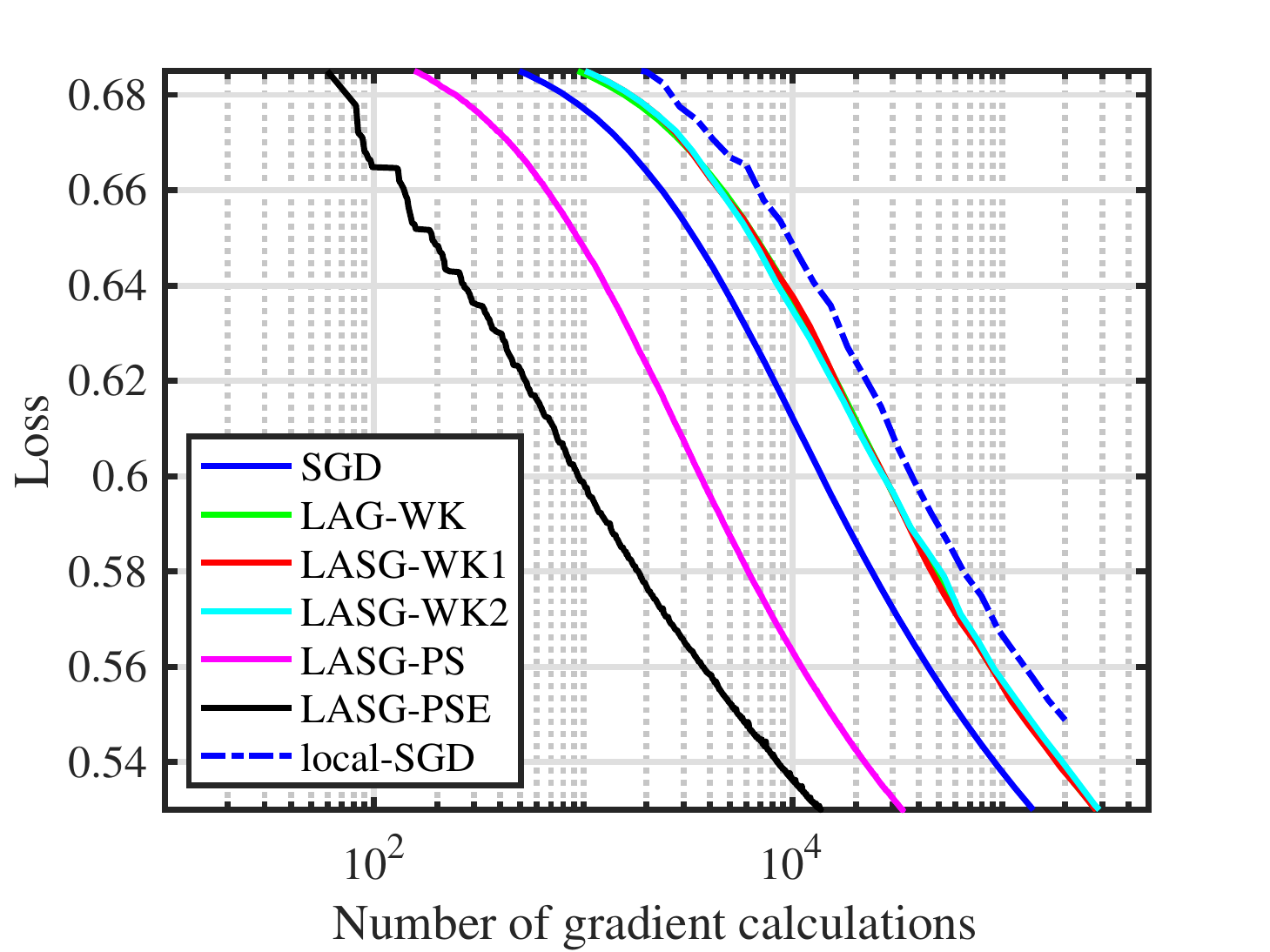}
\caption{Logistic regression on \textit{covtype} dataset in the heterogeneous setting.}
\label{fig:covtype_grad}
\end{figure}

In our numerical tests, all LASG algorithms achieve the same iteration complexity as SGD and outperform local-SGD in most cases.
Compared with SGD, LASG-WK2 and LASG-PSE reduce the number of communication rounds by around one order of magnitude for neural network training and even more for logistic regression. LASG-WK1 also reduce the communication by more than one order of magnitude for logistic regression.
Based on the results of LAG-WK and QSGD, it is evident that the selection rules \eqref{eqn:workerrule1}, \eqref{eqn:workerrule2} and \eqref{eqn:eserversule} of LASG-WK1, LASG-WK2 and LASG-PSE achieve more significant improvement in terms of saving communication and bits than the selection rule \eqref{eqn:LAGrule} of LAG-WK and stochastic quantization strategy of QSGD.

Although the performance of LASG-PS is not as impressive as other LASG algorithms in saving communication, it considerably saves local computation compared with other algorithms except LASG-PSE as shown in Figure \ref{fig:covtype_grad}. Moreover, LASG-PSE has performance gains in both saving communication and local computation.
The performance of LAQSG validates that the LASG algorithms can be easily equipped with stochastic quantization with additional benefits from quantization.

\section{Conclusions}
In this paper, we developed a class of LASG methods as communication-efficient variants of SGD.
LASG methods leverage a set of adaptive communication rules to detect and then skip less informative or redundant communication rounds between the server and workers during distributed learning.
To further reduce communication bandwidth, the quantized version of LASG is also presented.
Both LASG and their quantized version are simple to implement, and have convergence rate comparable to the original SGD. Extensions to nonsmooth and decentralized settings are also in our research agenda.

\clearpage
\onecolumn
\begin{center}
{\large \bf Supplementary materials for
``LASG: Lazily Aggregated Stochastic}
\end{center}
\begin{center}
{\large \bf  Gradients for Communication-Efficient Distributed Learning"}
\end{center}

\appendix
In this supplementary document, we first present some basic inequalities that will be used frequently in this document, and then present the missing derivations of some claims, as well as the proofs of all the lemmas and theorems in the paper, which is followed by details on our experiments. The content of this supplementary document is summarized as follows.


\section{Basics}
 Before starting the proof, we introduce our notation. First let $\htau_m^k=0$ if $m\in{\cal M}^k$; $\htau_m^k=\tau_m^k$ otherwise. Note that $\htau_m^k\in\{1,2\ldots,D\}$. With $\htau_m^k$, LASG's update can be simplified to
\begin{align}\label{eqn:LASGproof}
    \theta^{k+1}=\theta^k-\eta_k\sum\limits_{m\in{\cal M}}\ell(\theta^{k-\htau_m^k};\xi_m^{k-\htau_m^k}).
\end{align}
Define the $\sigma$-algebra $\Theta^k=\{\theta^l: \forall 1\leq l\leq k\}$. We also let $\theta^{-D},\theta^{-D+1},\ldots,\theta^{-1}=\theta^0$, which allows us to express our algorithms conveniently.
Some basic facts used in the proof are reviewed as follows.

\noindent{\bf Fact 1.} Assume that $X_1, X_2, \ldots, X_n\in\RR^p$ are independent random variables, and $EX_1=\cdots=EX_n=0$. Then
\begin{align}\label{eqn:equation1}
    \EE\Big[\|\sum\limits_{i=1}^nX_i\|^2\Big]=\sum\limits_{i=1}^n\EE\left[\|X_i\|^2\right].
\end{align}

\noindent{\bf Fact 2.} (Young's inequality) For any $\theta_1,\theta_2\in\RR^p,\varepsilon>0$,
\begin{align}\label{eqn:young}
    \dotp{\theta_1,\theta_2}\leq\frac{\|\theta_1\|^2}{2\varepsilon}+\frac{\varepsilon \|\theta_2\|^2}{2}.
\end{align}

As a consequence, we have
\begin{align}\label{eqn:young1}
    \|\theta_1+\theta_2\|^2\leq (1+\frac{1}{\varepsilon})\|\theta_1\|^2+(1+\varepsilon)\|\theta_2\|^2.
\end{align}

\noindent{\bf Fact 3.} (Cauchy-Schwarz inequality) For any $\theta_1,\theta_2,\ldots,\theta_n\in\RR^p$,
\begin{align}\label{eqn:cauchy}
    \big\|\sum\limits_{i=1}^n\theta_i\big\|^2\leq n\sum\limits_{i=1}^n\|\theta_i\|^2.
\end{align}

\noindent{\bf Fact 4.} For $k-D\leq l \leq k-\htau_m^k$,
\begin{align}\label{eqn:inequality1}
\begin{split}
    &\EE\left[\dotp{\nabla \cL(\theta^k),\nabla\ell(\theta^{l};\xi_m^k)-\nabla\ell(\theta^{l};\xi_m^{k-\htau_m^k})}\right]\\
    \overset{\text{(1a)}}{=}&\EE\left[\dotp{\nabla \cL(\theta^k)-\nabla \cL(\theta^{l}),\nabla\ell(\theta^{l};\xi_m^k)-\nabla\ell(\theta^{l};\xi_m^{k-\htau_m^k})}\right]\\
    \leq& L\EE\left[\|\theta^k-\theta^l\|\|\nabla\ell(\theta^{l};\xi_m^k)-\nabla\ell(\theta^{l};\xi_m^{k-\htau_m^k})\|\right]\\
    \overset{\text{(1b)}}{\leq}&\frac{L}{12 \sqrt{M} D\eta_k}\underbrace{\EE\left[\|\theta^k-\theta^{l}\|^2\right]}_{I_1}+\frac{6\sqrt{M} D L\eta_k}{2}\underbrace{\EE\left[\|\nabla\ell(\theta^{l};\xi_m^k)-\nabla\ell(\theta^{l};\xi_m^{k-\htau_m^k})\|^2\right]}_{I_2}\\
    \overset{\text{(1c)}}{\leq}&\frac{L}{12\sqrt{M}\eta_k}\!\sum\limits_{d=1}^{D}\!\EE\left[\|\theta^{k+1-d}\!-\!\theta^{k-d}\|^2\right]\!+\! 6\sqrt{M} D L\eta_k\sigma_m^2,
\end{split}
\end{align}
where (1a) holds due to
\begin{align*}
    &\EE\left[\dotp{\nabla\cL(\theta^{l}),\nabla\ell(\theta^l;\xi_m^k)-\nabla\ell(\theta^{l};\xi_m^{k-\htau_m^k})}\right]\\
    =&\EE\left[\EE\big[\dotp{\nabla\cL(\theta^{l}),\nabla\ell(\theta^l;\xi_m^k)-\nabla\ell(\theta^{l};\xi_m^{k-\htau_m^k})}\big|\Theta^l\big]\right]\\
    =&\EE\left[\dotp{\nabla\cL(\theta^l), \EE\big[ \nabla\ell(\theta^l;\xi_m^k)-\nabla\ell(\theta^{l};\xi_m^{k-\htau_m^k})\big|\Theta^{l}\big]}\right]\\
    =&\EE\left[\dotp{\nabla\cL(\theta^l), \nabla\cL_m(\theta^l)-\nabla\cL_m(\theta^l)}\right]\\
    =&0
\end{align*}
(1b) is a direct application of the Young's inequality \eqref{eqn:young}, and (1c) is a result of applying the Cauchy-Schwarz inequality \eqref{eqn:cauchy} to $I_1$, \eqref{eqn:equation1}, and Assumption \ref{assump:gradientestimator} to $I_2$,
\begin{align*}
    I_1=&\EE\Big[\|\sum\limits_{d=1}^{k-l}(\theta^{k+1-d}-\theta^{k-d})\|^2\Big]    \leq(k-l)\sum\limits_{d=1}^{k-l}\EE\Big[\|\theta^{k+1-d}-\theta^{k-d}\|^2\Big]\leq D\sum\limits_{d=1}^D\EE\Big[\|\theta^{k+1-d}-\theta^{k-D}\|^2\Big]\\
    I_2=&\EE\Big[\big\|\nabla\ell(\theta^l;\xi_m^k)-\nabla\ell(\theta^l;\xi_m^{k-\htau_m^k})\big\|^2\Big]\\
       =&\EE\Big[\big\|\big(\nabla\ell(\theta^l;\xi_m^k)-\nabla\cL_m(\theta^l)\big)+\big(\nabla\cL_m(\theta^l)-\nabla\ell(\theta^l;\xi_m^{k-\htau_m^k})\big)\big\|^2\Big]\\
       =&\EE\Big[\big\|\nabla\ell(\theta^l;\xi_m^k)-\nabla\cL_m(\theta^l)\big\|^2\Big]+\EE\Big[\big\|\nabla\cL_m(\theta^l)-\nabla\ell(\theta^l;\xi_m^{k-\htau_m^k})\big\|^2\Big]\\
    \leq&2\sigma_m^2.
\end{align*}

Similar to \eqref{eqn:inequality1}, it can be verified that
\begin{align}\label{eqn:inequality2}
\begin{split}
    \EE\left[\dotp{\nabla\cL(\theta^k),\nabla\cL_m(\theta^l)-\nabla\ell(\theta^l;\theta^{k-\htau_m^k})}\right]\leq \frac{L}{12\sqrt{M}\eta_k}\sum\limits_{d=1}^D\EE\left[\|\theta^{k+1-d}-\theta^{k-d}\|^2\right]+\frac{6\sqrt{M}DL\eta_k\sigma_m^2}{2}
\end{split}
\end{align}

\section{Proof of \eqref{eqn:variance-wk}, \eqref{eqn:variance-wk1} and \eqref{eqn:variance-wk2}}
The analysis in this part is analogous to that in \citep{ghadimi2013stochastic}. Before the proof, we define an auxiliary function,
\begin{align*}
    \psi_m(\theta)={\cal L}_m(\theta)-{\cal L}_m(\theta^{\star})-\dotp{\nabla{\cal L}_m(\theta^{\star}),\theta-\theta^{\star}}
\end{align*}
where $\theta^{\star}$ is a global minimizer of ${\cal L}$. And we assume that $\nabla\ell(\theta;\xi_m)$ is $\bar{L}$-Lipschitz continuous for any $\xi_m$, 
\begin{align*}
    \|\nabla\ell(\theta;\xi_m)-\nabla\ell(\theta^{\star};\xi_m)\|^2\leq 2\bar{L}(\ell(\theta;\xi_m)-\ell(\theta^{\star};\xi_m)-\dotp{\nabla\ell(\theta^{\star};\xi_m), \theta-\theta^{\star}}).
\end{align*}
Take expectation with respect to $\xi_m$ and we can obtain
\begin{align*}
    \EE_{\xi_m}[\|\nabla\ell(\theta;\xi_m)-\nabla\ell(\theta^{\star};\xi_m)\|^2]\leq2\bar{L}({\cal L}_m(\theta)-{\cal L}_m(\theta^{\star})-\dotp{\nabla{\cal L}_m(\theta^{\star}),\theta-\theta^{\star}})=2\bar{L}\psi_m(\theta).
\end{align*}
Note that $\nabla{\cal L}_m$ is also $\bar{L}$-Lipschitz continuous and thus
\begin{align*}
    \|\nabla{\cal L}_m(\theta)-\nabla{\cal L}_m(\theta^{\star})\|^2\leq2\bar{L}({\cal L}_m(\theta)-{\cal L}_m(\theta^{\star})-\dotp{\nabla{\cal L}_m(\theta^{\star}),\theta-\theta^{\star}})=2\bar{L}\psi_m(\theta).
\end{align*}

\subsection{Derivations of \eqref{eqn:variance-wk}} 
By \eqref{eqn:cauchy}, we can derive that
\begin{align*}
    \|\theta_1+\theta_2\|\leq 2\|\theta_1\|^2+2\|\theta_2\|^2\Rightarrow\|\theta_1\|^2\geq\frac{1}{2}\|\theta_1+\theta_2\|^2-\|\theta_2\|^2.
\end{align*}
As a consequence,
\begin{align*}
&\EE\Big[\big\|\nabla\ell(\theta^k;\xi_m^k)-\nabla\ell(\theta^{k-\tau_m^k};\xi_m^{k-\tau_m^k})\big\|^2\Big]\\
\geq&\frac{1}{2}\EE\Big[\big\|\big(\nabla\ell(\theta^k;\xi_m^k)-\nabla{\cal L}_m(\theta^k)\big)+\big(\nabla{\cal L}_m(\theta^{k-\tau_m^k})-\nabla\ell(\theta^{k-\tau_m^k};\xi_m^{k-\tau_m^k})\big)\big\|^2\Big]-\EE\Big[\big\|\nabla{\cal L}_m(\theta^k)-\nabla{\cal L}_m(\theta^{k-\tau_m^k})\big\|^2\Big]\\
=&\frac{1}{2}\EE\Big[\big\|\nabla\ell(\theta^k;\xi_m^k)-\nabla{\cal L}_m(\theta^k)\big\|^2\Big]+\frac{1}{2}\EE\Big[\big[\big\|\nabla\ell(\theta^{k-\tau_m^k};\xi_m^{k-\tau_m^k})-\nabla{\cal L}_m(\theta^{k-\tau_m^k})\big\|^2\big]\Big]\\
&+\underbrace{\EE\Big[\dotp{\nabla\ell(\theta^k;\xi_m^k)-\nabla{\cal L}_m(\theta^k), \nabla{\cal L}_m(\theta^{k-\tau_m^k})-\nabla\ell(\theta^{k-\tau_m^k};\xi_m^{k-\tau_m^k})}\Big]}_{F_1}-\EE\Big[\big\|\nabla{\cal L}_m(\theta^k)-\nabla{\cal L}_m(\theta^{k-\tau_m^k})\big\|^2\Big]
\end{align*}
where we used the fact that $F_1=0$ to obtain \eqref{eqn:variance-wk} since
\begin{align*}
    F_1= \EE\Big[\dotp{\underbrace{\EE\big[\nabla\ell(\theta^k;\xi_m^k)\big|\Theta^k\big]}_{=\nabla{\cal L}_m(\theta^k)}-\nabla{\cal L}_m(\theta^k), \nabla{\cal L}_m(\theta^{k-\tau_m^k})-\nabla\ell(\theta^{k-\tau_m^k};\xi_m^{k-\tau_m^k})}\Big]=0.
\end{align*}

\subsection{Derivations of \eqref{eqn:variance-wk1}}
Recall that
\begin{align*}
    \tdelta_m^k-\tdelta_m^{k-\tau_m^k}=&\big(\nabla\ell(\theta^k;\xi_m^k)-\nabla\ell(\tilde{\theta};\xi_m^k)+\nabla{\cal L}_m(\tilde{\theta})\big)-\big(\nabla\ell(\theta^{k-\tau_m^k};\xi_m^{k-\tau_m^k})-\nabla\ell(\tilde{\theta};\xi_m^{k-\tau_m^k})+\nabla{\cal L}_m(\tilde{\theta})\big)\\
    =&\underbrace{\big(\nabla\ell(\theta^k;\xi_m^k)-\nabla\ell(\tilde{\theta};\xi_m^k)+\nabla\psi_m(\tilde{\theta})\big)}_{g_m^k}-\underbrace{\big(\nabla\ell(\theta^{k-\tau_m^k};\xi_m^{k-\tau_m^k})-\nabla\ell(\tilde{\theta};\xi_m^{k-\tau_m^k})+\nabla\psi_m(\tilde{\theta})\big)}_{g_m^{k-\tau_m^k}}.
\end{align*}
And by \eqref{eqn:cauchy}, $\|\tdelta_m^k-\tdelta_m^{k-\tau_m^k}\|^2\leq 2\|g_m^k\|^2+2\|g_m^{k-\tau_m^k}\|^2$.
\begin{align*}
&    \EE[\|g_m^k\|^2]\\\leq&2\EE[\|\nabla\ell(\theta^k;\xi_m^k)-\nabla\ell(\theta^{\star};\xi_m^k)\|^2]+2\EE[\|\nabla\ell(\ttheta;\xi_m^k)-\nabla\ell(\theta^{\star};\xi_m^k)-\nabla\psi_m(\ttheta)\|^2]\\
    =&2\EE[\EE[\|\nabla\ell(\theta^k;\xi_m^k)-\nabla\ell(\theta^{\star};\xi_m^k)\|^2|\Theta^k]]+2\EE[\|\nabla\ell(\ttheta;\xi_m^k)-\nabla\ell(\theta^{\star};\xi_m^k)-\EE[\nabla\ell(\ttheta;\xi_m^k)-\nabla\ell(\theta^{\star};\xi_m^k)|\Theta^k]\|^2]\\
    \leq&4\bar{L}\EE\psi_m(\theta^k)+2\EE[\|\nabla\ell(\ttheta;\xi_m^k)-\nabla\ell(\theta^{\star};\xi_m^k)\|^2]\\
    =&4\bar{L}\EE\psi_m(\theta^k)+2\EE[\EE[\|\nabla\ell(\ttheta;\xi_m^k)-\nabla\ell(\theta^{\star};\xi_m^k)\|^2|\Theta^k]]\\
    \leq&4\bar{L}\EE\psi_m(\theta^k)+4\bar{L}\EE\psi_m(\ttheta).
\end{align*}
By nonnegativity of $\psi_m$, we have
\begin{align*}
    \EE[\|g_m^k\|^2]\leq4\bar{L}\sum\limits_{m\in{\cal M}}\EE\psi_m(\theta^k)+4\bar{L}\sum\limits_{m\in{\cal M}}\EE\psi_m(\ttheta)=4\bar{L}(\EE{\cal L}(\theta^k)-{\cal L}(\theta^{\star}))+4\bar{L}(\EE{\cal L}(\ttheta)-{\cal L}(\theta^{\star}))
\end{align*}
Similarly, we can prove $\EE[\|g_m^{k-\tau_m^k}\|^2]\leq4\bar{L}(\EE{\cal L}(\theta^{k-\tau_m^k})-{\cal L}(\theta^{\star}))+4\bar{L}(\EE{\cal L}(\ttheta)-{\cal L}(\theta^{\star}))$.
Therefore, it follows that
\begin{align*}
    \EE[\|\tdelta_m^k-\tdelta_m^{k-\tau_m^k}\|^2]\leq 8\bar{L}(\EE{\cal L}(\theta^k)-{\cal L}(\theta^{\star}))+8\bar{L}(\EE{\cal L}(\theta^{k-\tau_m^k})-{\cal L}(\theta^{\star}))+16\bar{L}(\EE{\cal L}(\ttheta)-{\cal L}(\theta^{\star})).
\end{align*}

\subsection{Derivations of \eqref{eqn:variance-wk2}}
\begin{align*}
    \nabla\ell(\theta^k;\xi_m^k)-\nabla\ell(\theta^{k-\tau_m^k};\xi_m^k)=&\big(\nabla\ell(\theta^k;\xi_m^k)-\nabla\ell(\theta^{k-\tau_m^k};\xi_m^k)+\nabla{\cal L}_m(\theta^{k-\tau_m^k})\big)-\nabla{\cal L}_m(\theta^{k-\tau_m^k})\\
    =&\big(\nabla\ell(\theta^k;\xi_m^k)-\nabla\ell(\theta^{k-\tau_m^k};\xi_m^k)+\nabla\psi_m(\theta^{k-\tau_m^k})\big)-\nabla\psi_m(\theta^{k-\tau_m^k})
\end{align*}
Similar to the proof of \eqref{eqn:variance-wk1}, we can obtain
\begin{align*}
    \EE[\|\nabla\ell(\theta^k;\xi_m^k)-\nabla\ell(\theta^{k-\tau_m^k};\xi_m^k)+\nabla\psi_m(\theta^{k-\tau_m^k})\|^2]\leq
    4\bar{L}(\EE{\cal L}(\theta^k)-{\cal L}(\theta^{\star}))+4\bar{L}(\EE{\cal L}(\theta^{k-\tau_m^k})-{\cal L}(\theta^{\star})).
\end{align*}
Combined with 
\begin{align*}
	\EE[\|\nabla\psi_m(\theta^{k-\tau_m^k})\|^2]&=\EE[\|\nabla{\cal L}(\theta^{k-\tau_m^k})-\nabla{\cal L}(\theta^{\star})\|^2]\\	
	&\leq 2\bar{L}\EE\psi(\theta^{k-\tau_m^k})\leq 2\bar{L}(\EE{\cal L}(\theta^{k-\tau_m^k})-{\cal L}(\theta^{\star}))
\end{align*}
we have
\begin{align*}
    \EE[\|\nabla\ell(\theta^k;\xi_m^k)-\nabla\ell(\theta^{k-\tau_m^k};\xi_m^k)\|^2]\leq 8\bar{L}(\EE{\cal L}(\theta^k)-{\cal L}(\theta^{\star}))+12\bar{L}(\EE{\cal L}(\theta^{k-\tau_m^k})-{\cal L}(\theta^{\star})).
\end{align*}

\subsection{Derivations of \eqref{eqn:variance-server}}
From the LASG update, we have
\begin{align*}
\theta^{k}-\theta^{k-\tau_m^k}=\sum_{d=1}^{\tau_m^k} \eta_{k-d} \sum_{m\in{\cal M}}\nabla\ell(\theta^{k-d-\tau_m^{k-d}};\xi_m^{k-d-\tau_m^{k-d}}).
\end{align*}
Then the LASG-PS condition \eqref{eqn:serverrule} implies that 
{\small
\begin{align*}
&~~~~\EE[\|\theta^k-\theta^{k-\tau_m^k}\|^2]\\
&=\EE[\|\sum\limits_{d=1}^{\tau_m^k}\sum\limits_{m\in{\cal M}}\eta_{k-d}\nabla\ell(\theta^{k-d-\tau_m^{k-d}};\xi_m^{k-d-\tau_m^{k-d}})\|^2]\\
&\leq \tau_m^k\sum\limits_{d=1}^{\tau_m^k}\eta_{k-d}^2\EE[\|\sum\limits_{m\in{\cal M}}\nabla\ell(\theta^{k-d-\tau_m^{k-d}};\xi_m^{k-d-\tau_m^{k-d}})\|^2]\\
&=\tau_m^k\sum\limits_{d=1}^{\tau_m^k}\eta_{k-d}^2\EE[\|\sum\limits_{m\in{\cal M}}\nabla\ell(\theta^{k-d-\tau_m^{k-d}};\xi_m^{k-d-\tau_m^{k-d}})-\nabla{\cal L}_m(\theta^{k-d-\tau_m^{k-d}}\|^2]\\
&~~~~~~+\tau_m^k\sum\limits_{d=1}^{\tau_m^k}\eta_{k-d}^2\EE[\|\sum\limits_{m\in{\cal M}}\nabla{\cal L}_m(\theta^{k-d-\tau_m^{k-d}})\|^2]\\
&\leq\tau_m^k\sum\limits_{d=1}^{\tau_m^k}\eta_{k-d}^2\sum\limits_{m\in{\cal M}}\sigma_m^2+\tau_m^k\sum\limits_{d=1}^{\tau_m^k}\eta_{k-d}^2\EE[\|\sum\limits_{m\in{\cal M}}\nabla{\cal L}_m(\theta^{k-d-\tau_m^{k-d}})\|^2]\\
&\leq\tau_m^k\sum\limits_{d=1}^{\tau_m^k}\sum\limits_{m\in{\cal M}}\sigma_m^2\eta_{k-d}^2
+2\tau_m^k\sum\limits_{d=1}^{\tau_m^k}\EE[\|\nabla{\cal L}(\theta^{k-d})\|^2]\eta_{k-d}^2
+2\tau_m^k\sum\limits_{d=1}^{\tau_m^k}\sum\limits_{m\in{\cal M}}L_m\EE[\|\theta^{k-d}-\theta^{k-d-\tau_m^{k-d}}\|^2]\eta_{k-d}^2\\
&\leq D\sum\limits_{d=1}^D\sum\limits_{m\in{\cal M}}\sigma_m^2\eta_{k-D}^2+2D\sum\limits_{d=1}^D\EE[\|\nabla{\cal L}(\theta^{k-d})\|^2]\eta_{k-D}^2+2D\sum\limits_{d=1}^D\sum\limits_{m\in\cal M}L_m^2\EE[\|\theta^{k-d}-\theta^{k-d-\tau_m^{k-d}}\|^2]\eta_{k-D}^2.
\end{align*}}

\section{Proof of Lemma \ref{lemma:lossdescent}}
Due to the smoothness of $\cL$ in Assumption \ref{assump:smoothness}, we have
\begin{align}\label{eqn:lossdescent1}
   \EE\cL(\theta^{k+1})=&\EE\cL(\theta^k-\eta_k\sum\limits_{m\in{\cal M}}\nabla\ell(\theta^{k-\htau_m^k};\xi_m^{k-\htau_m^k})) \nonumber\\
    \overset{\eqref{assump:smoothness}}{\leq}&\EE\cL(\theta^k)-\eta_k\underbrace{\EE\Big[\dotp{\nabla\cL(\theta^k),\sum\limits_{m\in{\cal M}}\nabla\ell(\theta^{k-\htau_m^k};\xi_m^{k-\htau_m^k})}\Big]}_{J_1}+\frac{L\eta_k^2}{2}\underbrace{\EE\Big[\big\|\sum\limits_{m\in{\cal M}}\nabla\ell(\theta^{k-\htau_m^k};\xi_m^{k-\htau_m^k})\big\|^2\Big]}_{J_2}.
\end{align}

We decompose $J_1,J_2$ as follows, 
\begin{align*}
    J_1=&\EE\Big[\dotp{\nabla\cL(\theta^k),\sum\limits_{m\in{\cal M}}\nabla\ell(\theta^k;\xi_m^k)}\Big]- \underbrace{\sum\limits_{m\in{\cal M}}\EE \Big[\dotp{\nabla\cL(\theta^k), \nabla\ell(\theta^k;\xi_m^k)-\nabla\ell(\theta^{k-\htau_m^k};\xi_m^{k-\htau_m^k})}\Big]}_{H_1}\\
       =&\EE\Big[\dotp{\nabla\cL(\theta^k),\sum\limits_{m\in{\cal M}}\EE\left[\left.\nabla\ell(\theta^k;\xi_m^k)\right|\Theta^k\right]}\Big]- H_1\\
       =&\EE\left[\|\nabla\cL(\theta^k)\|^2\right]- H_1,\\
    J_2=&\EE\Big[\big\|\sum\limits_{m\in{\cal M}}\nabla\ell(\theta^{k-\htau_m^k};\xi_m^{k-\htau_m^k})-\nabla\cL(\theta^k)+\nabla\cL(\theta^k)\big\|^2\Big]\\
       =&\underbrace{\EE\Big[\big\|\sum\limits_{m\in{\cal M}}\nabla\ell(\theta^{k-\htau_m^k};\xi_m^{k-\htau_m^k})-\nabla\cL(\theta^k)\big\|^2\Big]}_{H_2}+\EE\left[\|\nabla\cL(\theta^k)\|^2\right]\\
        &-2\EE\Big[\dotp{\nabla\cL(\theta^k), \underbrace{\sum\limits_{m\in{\cal M}}\EE\left[\left.\nabla\ell(\theta^k;\xi_m^k)\right|\Theta^k\right]}_{\nabla{\cal L}(\theta^k)=}}\Big]+2\EE\Big[\dotp{\nabla\cL(\theta^k), \sum\limits_{m\in{\cal M}}\nabla\ell(\theta^{k-\htau_m^k};\xi_m^{k-\htau_m^k})}\Big]\\
       =&H_2+\EE\left[\|\nabla\cL(\theta^k)\|^2\right]-2H_1.
\end{align*}

By taking expectation and substituting $J_1,J_2$ in \eqref{eqn:lossdescent1}, we obtain
\begin{align}\label{eqn:lossdescent2}
    \EE\cL(\theta^{k+1})\leq&\EE\cL(\theta^k)-\left(\eta_k-\frac{L\eta_k^2}{2}\right)\EE\left[\|\nabla\cL(\theta^k)\|^2\right]+\left(\frac{L\eta_k^2}{2}\right)H_2+\left(\eta_k-L\eta_k^2\right)H_1.
\end{align}

We analyze $H_1$, $H_2$ separately for different rules. First for LASG-WK1's rule \eqref{eqn:workerrule1},
\begin{align*}
        H_1\overset{\text{(2a)}}=&\sum\limits_{m\in{\cal M}}\EE\left[\dotp{\nabla\cL(\theta^k),\tilde{\delta}_m^k-\tilde{\delta}_m^{k-\htau_m^k}}\right]+\sum\limits_{m\in{\cal M}}\EE\left[\dotp{\nabla\cL(\theta^k),\nabla\ell(\ttheta;\xi^k)-\nabla\ell(\ttheta,\xi_m^{k-\htau_m^k})}\right]\nonumber\\
        \overset{\text{(2b)}}{\leq}&\frac{L\eta_k}{2}\EE\Big[\big\|\nabla\cL(\theta^k)\big\|^2\Big]+\sum\limits_{d=1}^{D}\left(\frac{c_d}{2L\eta_k}+\frac{\sqrt{M}L}{12\eta_k}\right)\EE\Big[\big\|\theta^{k+1-d}-\theta^{k-d}\big\|^2\Big]+6\sqrt{M} D L\eta_k\sum\limits_{m\in{\cal M}}\sigma_m^2.
\end{align*}
    where (2a) is due to the definition of $\delta^k$ and (2b) is obtained by \eqref{eqn:workerrule1}, \eqref{eqn:young} with $\varepsilon=\frac{M}{L\eta_k}$, and \eqref{eqn:inequality1} with $\theta^{l}=\ttheta$. Note that the definition of $\tilde{\theta}$ in Algorithm \ref{alg: LASG-WK1} implies $l=\lfloor\frac{k}{D}\rfloor\leq k-\hat{\tau}_m^k$.
Similarly, for LASG-WK2's rule \eqref{eqn:workerrule2}, we apply \eqref{eqn:young} with $\varepsilon=\frac{M}{L\eta_k}$ and \eqref{eqn:inequality1} with $l=k-\htau_m^k$.
\begin{align*}
    H_1&=\sum\limits_{m\in{\cal M}}\EE\left[\dotp{\nabla\cL(\theta^k), \nabla\ell(\theta^k;\xi_m^k)-\nabla\ell(\theta^{k-\htau_m^k};\xi_m^k)}\right]\nonumber\\
    &+\EE\left[\dotp{\nabla\cL(\theta^k), \nabla\ell(\theta^{k-\htau_m^k};\xi_m^k)-\nabla\ell(\theta^{k-\htau_m^k};\xi_m^{k-\htau_m^k})}\right]\nonumber\\
    &\leq\frac{L\eta_k}{2}\EE\left[\|\nabla\cL(\theta^k)\|^2\right]+\sum\limits_{d=1}^{D}\left(\frac{c_d}{2L\eta_k}+\frac{\sqrt{M}L}{12\eta_k}\right)\EE\left[\|\theta^{k+1-d}-\theta^{k-d}\|^2\right]+6\sqrt{M}DL\eta_k\sum\limits_{m\in{\cal M}}\sigma_m^2.
\end{align*}

For LASG-PS's rule \eqref{eqn:serverrule}, apply $\EE\left[\left.\nabla\ell(\theta^k;\xi_m^k)\right|\Theta^k\right]=\nabla\cL_m(\theta^k)$, \eqref{eqn:young} with $\varepsilon=\frac{M}{L\eta_k}$ and \eqref{eqn:inequality2} with $l=k-\htau_m^k$ to get
\begin{align*}
    H_1=&\sum\limits_{m\in{\cal M}}\EE\left[\dotp{\nabla\cL(\theta^k), \nabla\cL_m(\theta^k)-\nabla\cL_m(\theta^{k-\htau_m^k})}\right] +\EE\left[\dotp{\nabla\cL(\theta^k), \nabla\cL_m(\theta^{k-\htau_m^k})-\nabla\ell(\theta^{k-\htau_m^k};\xi_m^{k-\htau_m^k})}\right]\nonumber\\
    \leq&\frac{L\eta_k}{2}\EE\left[\|\nabla\cL(\theta^k)\|^2\right]+\sum\limits_{d=1}^{D}\left(\frac{c_d}{2L\eta_k}+\frac{\sqrt{M}L}{12\eta_k}\right)\EE\left[\|\theta^{k+1-d}-\theta^{k-d}\|^2\right]+\frac{6\sqrt{M}D L\eta_k}{2}\sum\limits_{m\in{\cal M}}\sigma_m^2.
\end{align*}

Now we deal with $H_2$. For LASG-WK1,
\begin{align*}
    H_2\overset{\eqref{eqn:cauchy}}\leq&3\EE\Big[\big\|\sum\limits_{m\in{\cal M}}\tilde{\delta}_m^k-\tilde{\delta}_m^{k-\htau_m^k}\big\|^2\Big]+3\EE\Big[\big\|\sum\limits_{m\in{\cal M}}\nabla\ell(\theta^k,\xi_m^k)-\nabla\cL(\theta^k))\big\|^2\Big]\nonumber\\
    &+3\EE\Big[\big\|\sum\limits_{m\in{\cal M}}(\nabla\ell(\ttheta;\xi_m^k)-\nabla\cL_m(\ttheta))+\sum\limits_{m\in{\cal M}}(\nabla\cL_m(\ttheta)-\nabla\ell(\ttheta;\xi_m^{k-\htau_m^k}))\big\|^2\Big]\nonumber\\
    \overset{\eqref{eqn:workerrule1},\eqref{eqn:boundedvariance},\eqref{eqn:equation1}}\leq&3\sum\limits_{d=1}^{D}c_d\EE\big[\|\theta^{k+1-d}-\theta^{k-d}\|^2\big]+9\sum\limits_{m\in{\cal M}}\sigma_m^2.
\end{align*}

For LASG-WK2,
\begin{align*}%
    H_2\overset{\eqref{eqn:cauchy}}\leq&2\EE\Big[\big\|\sum\limits_{m\in{\cal M}}\big(\nabla\ell(\theta^{k-\htau_m^k},\xi_m^{k-\htau_m^k})-\nabla\ell(\theta^{k};\xi_m^{k})\big)\big\|^2\Big]+2\EE\Big[\big\|\sum\limits_{m\in{\cal M}}\big(\nabla\ell(\theta^{k};\xi_m^{k})-\nabla\cL_m(\theta^k)\big)\big\|^2\Big]\\
    \overset{\eqref{eqn:workerrule2},\eqref{eqn:boundedvariance},\eqref{eqn:equation1}}\leq&2\sum\limits_{d=1}^{D}c_d\EE\left[\|\theta^{k+1-d}-\theta^{k-d}\|^2\right]+2\sum\limits_{m\in{\cal M}}\sigma_m^2\leq 3\sum\limits_{d=1}^{D}c_d\EE\|\theta^{k+1-d}-\theta^{k-d}\|^2+9\sum\limits_{m\in{\cal M}}\sigma_m^2.
\end{align*}

For LASG-PS,
\begin{align*}%
H_2\overset{\eqref{eqn:cauchy}}\leq&2\EE\Big[\big\|\sum\limits_{m\in{\cal M}}\big(\nabla\ell(\theta^{k-\htau_m^k},\xi_m^{k-\htau_m^k})-\nabla\cL_m(\theta^{k-\htau_m^k})\big)\Big\|^2\Big]+2\EE\Big[\big\|\sum\limits_{m\in{\cal M}}\big(\nabla\cL_m(\theta^{k-\htau_m^k})-\nabla\cL_m(\theta^k)\big)\big\|^2\Big]\\
    \overset{\eqref{eqn:serverrule},\eqref{eqn:boundedvariance},\eqref{eqn:equation1}}\leq&2\sum\limits_{d=1}^{D}c_d\EE\left[\|\theta^{k+1-d}-\theta^{k-d}\|^2\right]+2\sum\limits_{m\in{\cal M}}\sigma_m^2\leq 3\sum\limits_{d=1}^{D}c_d\EE\|\theta^{k+1-d}-\theta^{k-d}\|^2+9\sum\limits_{m\in{\cal M}}\sigma_m^2.
\end{align*}

Plug $H_1$, $H_2$ back into \eqref{eqn:lossdescent2} and get
\begin{align}\label{eqn.32}
    \EE\cL(\theta^{k+1})\leq\EE\cL(\theta^k)&-\left(\eta_k-L\eta_k^2+\frac{L^2\eta_k^3}{2}\right)\EE\left[\|\nabla\cL(\theta^k)\|^2\right]\nonumber\\
    &+\sum\limits_{d=1}^D\left(\left(\frac{3L\eta_k^2}{2}+\frac{1}{2L}-\frac{\eta_k}{2}\right)c_d+(1-L\eta_k)\frac{\sqrt{M}L}{12}\right)\EE\left[\|\theta^{k+1-d}-\theta^{k-d}\|^2\right]\nonumber\\
    &+L\eta_k^2\left(\frac{9}{2}+6\sqrt{M}D(1-L\eta_k)\right)\sum\limits_{m\in{\cal M}}\sigma_m^2.
\end{align}
Since $L\eta_k\leq 1$, \eqref{eqn.32} can be simplified as \eqref{eqn:lossdescent}.
\section{Proof of Lemma \ref{lemma:lyapunovdescent}}
By definition of $V^k$, it follows that
\begin{align*}
    V^{k+1}-V^k
    =&\EE\cL(\theta^{k+1})-\EE\cL(\theta^k)+\gamma_1\EE\left[\|\theta^{k+1}-\theta^{k}\|^2\right]\\
     &+\sum\limits_{d=1}^{D-1}(\gamma_{d+1}-\gamma_d)\EE\left[\|\theta^{k+1-d}-\theta^{k-d}\|^2\right]-\gamma_D\EE\left[\|\theta^{k+1-D}-\theta^{k-D}\|^2\right].
\end{align*}
First decompose $\EE\left[\|\theta^{k+1}-\theta^k\|^2\right]$,
\begin{align*}
    \EE\left[\|\theta^{k+1}-\theta^k\|^2\right]=&\eta_k^2\EE\Big[\|\sum\limits_{m\in{\cal M}}\nabla\ell(\theta^{k-\htau_m^k};\xi_m^{k-\htau_m^k})-\nabla\cL(\theta^k)+\nabla\cL(\theta^k)\|^2\Big]\\
    \leq&2\eta_k^2\EE\Big[\|\nabla\cL(\theta^k)\|^2\Big]+2\eta_k^2\underbrace{\EE\Big[\|\sum\limits_{m\in{\cal M}}\nabla\ell(\theta^{k-\htau_m^k};\xi_m^{k-\htau_m^k})-\nabla\cL(\theta^k)\|^2\Big]}_{=H_2}\\
    \leq&2\eta_k^2\EE\left[\|\nabla\cL(\theta^k)\|^2\right]+\sum\limits_{d=1}^D6c_d\eta_k^2\EE\left[\|\theta^{k+1-d}-\theta^{k-d}\|^2\right]+18\eta_k^2\sum\limits_{m\in{\cal M}}\sigma_m^2.
\end{align*}

By Lemma \ref{lemma:lossdescent}, it follows that
\begin{align}\label{eqn:lyapunovdescent}
\begin{split}
    V^{k+1}-V^k \leq&-\underbrace{\left(\eta_k-(L+2\gamma_1)\eta_k^2\right)}_{B_0^k}\EE\left[\|\nabla\cL(\theta^k)\|^2\right]\\
    &+\sum\limits_{d=1}^{D-1}\underbrace{\left((\eta_k+\frac{1}{2L})c_d+\frac{\sqrt{M}L}{12}+6c_d\gamma_1\eta_k^2+\gamma_{d+1}-\gamma_d\right)}_{A_d^k}\EE\left[\|\theta^{k+1-d}-\theta^{k-d}\|^2\right]\\
    &+\underbrace{\left((\eta_k+\frac{1}{2L})c_D+\frac{\sqrt{M}L}{12}+6c_D\gamma_1\eta_k^2-\gamma_D\right)}_{A_D^k}\EE\left[\|\theta^{k+1-D}-\theta^{k-D}\|^2\right]\\
    &+\underbrace{\left((\frac{9}{2}+6\sqrt{M} D)L\eta_k^2+18\gamma_1\eta_k^2\right)}_{B_1^k}\sum\limits_{m\in{\cal M}}\sigma_m^2.
\end{split}
\end{align}

\section{Proof of Theorem \ref{thm:nonconvex}}
Let $\eta_k\leq \bar{\eta}\leq \frac{1}{L}$ for $k=1,2\ldots, K$. To ensure $A_d^k\geq 0$, it is sufficient to choose $\gamma_d$ satisfying,
\begin{align*}
    &(\bar{\eta}+\frac{1}{2L})c_d+\frac{\sqrt{M}L}{12}+6c_d\gamma_1\bar{\eta}^2+\gamma_{d+1}-\gamma_d=0,\quad 0\leq d\leq D-1,\\
    &(\bar{\eta}+\frac{1}{2L})c_D+\frac{\sqrt{M}L}{12}+6c_D\gamma_1\bar{\eta}^2-\gamma_D=0.
\end{align*}

Solve the linear equations above and get
\begin{align*}
    \gamma_1=\frac{(\bar{\eta}+\frac{1}{2L})\sum\limits_{d=1}^Dc_d+\frac{\sqrt{M}DL}{12}}{1-6\sum\limits_{d=1}^Dc_d\bar{\eta}^2}.
\end{align*}

Select $c_d\leq\min\{\frac{1}{12D\bar{\eta}^2},\frac{\sqrt{M}L^2}{18}\}$ such that
\begin{align*}
    \gamma_1\leq \frac{\sqrt{M}DL}{3}.
\end{align*}

Let $\eta_k\leq\frac{1}{2L+\frac{4}{3}\sqrt{M}DL}\leq\frac{1}{2L+4\gamma_1}$ and then $B_0^k\geq\frac{\eta_k}{2}$. By taking sum $\sum\limits_{k=0}^{K-1} V^{k+1}-V^k$, it follows that
\begin{align*}
  \sum\limits_{k=0}^{K-1}\frac{\eta_k}{2}\EE\Big[\|\nabla\cL(\theta^k)\|^2\Big]\leq\cL(\theta^0)-\cL^*+\sum\limits_{k=0}^{K-1}\left(\frac{9}{2}+12\sqrt{M}D\right)L\eta_k^2\sum\limits_{m\in{\cal M}}\sigma_m^2\nonumber
\end{align*}

Define a random variable $R(K)$ with $\mathbb{P}(R(K)=j)=\frac{\eta_j}{\sum\limits_{k=0}^{K-1}\eta_k}$ and then
\begin{align*}
    \EE\Big[\|\nabla\cL(\theta^{R(K)})\|^2\Big]\leq\frac{\cL(\theta^0)-\cL^*+\sum\limits_{k=0}^{K-1}\left(\frac{9}{2}+12\sqrt{M}D\right)L\eta_k^2\sum\limits_{m\in{\cal M}}\!\sigma_m^2}{\sum\limits_{k=0}^{K-1}\frac{\eta_k}{2}}.
\end{align*}

Specially, if we choose constant stepsize
\begin{align*}
    \eta_k=\min\left\{\frac{1}{2L+\frac{4}{3}\sqrt{M}DL}, \frac{c_\eta}{\sqrt{K}}\right\}
\end{align*}
where $c_{\eta}>0$ is a constant, then
\begin{align*}
    \frac{1}{K}\sum\limits_{k=0}^{K-1}\EE\left[\|\nabla\cL(\theta^{k})\|^2\right]\leq&\frac{\cL(\theta^0)-\cL^*+K\left(\frac{9}{2}+12\sqrt{M}D\right)L\eta_0^2\sum\limits_{m\in{\cal M}}\sigma_m^2}{\frac{1}{2}K\eta_0}\\
    \leq&\frac{(4L+\frac{8}{3}\sqrt{M}DL)(\cL(\theta^0)-\cL^*)}{K}+\frac{2(\cL(\theta^0)-\cL^*)}{c_{\eta}\sqrt{K}}+\frac{c_{\eta}\left(9+24\sqrt{M}D\right)L\sum\limits_{m\in{\cal M}}\sigma_m^2}{\sqrt{K}}.
\end{align*}

\section{Proof of Theorem \ref{thm:stronglyconvex}}
By strong convexity of $\cL(\theta)$,
\begin{align}\label{eqn:strongconvexity}
    2\mu(\cL(\theta)-\cL(\theta^*))\leq\|\nabla\cL(\theta)\|^2.
\end{align}

Then \eqref{eqn:lyapunovdescent} can be rewritten as
\begin{align*}
     V^{k+1}-V^k\leq -2\mu B_0^k\EE(\cL(\theta^k)-\cL(\theta^*))+B_1^k\sum\limits_{m\in{\cal M}}\sigma_m^2+\sum\limits_{d=1}^{D}A_d^k
    \EE\Big[\|\theta^{k+1-d}-\theta^{k-d}\|^2\Big].
\end{align*}

We will choose $\gamma_d$ such that $A_d^k\leq -2\mu B_0^k\gamma_d$ for $d=1,2\ldots,D$. Then
\begin{align*}
    V^{k+1}\leq &(1-2\mu B_0^k)V^k+B_1^k\sum\limits_{m\in{\cal M}}\sigma_m^2\\
    \leq &\prod\limits_{j=0}^k(1-2\mu B_0^j)V^0+\sum\limits_{j=0}^k B_1^j\prod\limits_{i=j+1}^k(1-2\mu B_0^i)\sum\limits_{m\in{\cal M}}\sigma_m^2.
\end{align*}

If $\eta_k\leq\bar{\eta}\leq\frac{1}{L+2\gamma_1}$, then $0\leq B_0^k\leq\eta_k$. It is sufficient to choose $\gamma_d$ satisfying
\begin{align*}
    &(\bar{\eta}+\frac{1}{2L})c_d+\frac{\sqrt{M}L}{12}+6c_d\gamma_1\bar{\eta}^2+\gamma_{d+1}-\gamma_d=-2\mu\bar{\eta}\gamma_1,\quad 1\leq d\leq D-1,\\
    &(\bar{\eta}+\frac{1}{2L})c_D+\frac{\sqrt{M}L}{12}+6c_D\gamma_1\bar{\eta}^2-\gamma_D=-2\mu\bar{\eta}\gamma_1.
\end{align*}

Solve the linear equations above and get
\begin{align*}
    \gamma_1=\frac{(\bar{\eta}+\frac{1}{2L})\sum\limits_{d=1}^Dc_d+\frac{\sqrt{M}DL}{12}}{(1-6\sum\limits_{d=1}^Dc_d\bar{\eta}^2-2\mu D\bar{\eta})}.
\end{align*}

Select $c_d\leq\min\{\frac{1}{24D\bar{\eta}^2},\frac{\sqrt{M}L^2}{18}\}$ and $\bar{\eta}\leq\frac{1}{8\mu D}$. Then, we obtain
\begin{align*}
\gamma_1\leq\frac{\sqrt{M}DL}{3},
\end{align*}
and
\begin{align*}
    V^{k+1}\leq(1-\mu\eta_k)V^k+\underbrace{\left(\frac{9}{2}+12\sqrt{M}D\right)L\sum\limits_{m\in{\cal M}}\sigma_m^2}_R\eta_k^2.
\end{align*}

Let $\eta_k=\frac{2}{\mu(k+K_0)}$, where $K_0=\max\{\frac{2(L+\frac{2}{3}\sqrt{M}DL)}{\mu}, 16D\}$ to ensure $\eta_k\leq\frac{1}{L+2\gamma_1}$ and $\eta_k\leq\frac{1}{8\mu D}$. 
\begin{align*}
    V^{K}\leq&\prod\limits_{k=0}^{K-1}(1-\mu\eta_k)V^0+R\sum\limits_{k=0}^{K-1}\eta_k^2\prod\limits_{j=k+1}^{K-1}(1-\mu\eta_j)\\
    \leq&\frac{(K_0-2)(K_0-1)}{(K+K_0-2)(K+K_0-1)}V^0+\frac{R}{\mu^2}\sum\limits_{k=0}^{K-1}\frac{4}{(k+K_0)^2}\frac{(k+K_0-1)(k+K_0)}{(K+K_0-2)(K+K_0-1)}\\
    \leq&\frac{(K_0-1)^2}{(K+K_0-1)^2}V^0+\frac{4RK}{\mu^2(K+K_0-1)^2}\\
    =&\frac{(K_0-1)^2}{(K+K_0-1)^2}(\cL(\theta^K)-\cL^{\star})+\frac{4RK}{\mu^2(K+K_0-1)^2}.
\end{align*}
\section{Proof of Theorem \ref{thm:LAQSG}}
In this section we prove the convergence in Theorem \ref{thm:LAQSG}. Let $\EE_Q$ to denote the expectation with respect to the stochastic quantization $Q$. As a results of [Lemma 3.1, \cite{alistarh2017qsgd}] and Assumption \ref{assump:boundedgradient}, $b$-bit quantized gradients (1 bit for sign) have the following properties ,
\begin{align*}
    &\EE_Q\left[Q(\theta;\xi_m)\right]=\nabla\ell(\theta;\xi_m);\\
    &\EE_Q\left[\|Q(\theta;\xi_m)-\nabla\ell(\theta;\xi_m)\|^2\right]
    \leq\min\Big\{\frac{d}{(2^{b-1}-1)^2},\frac{\sqrt{d}}{2^{b-1}-1}\Big\}\|\nabla\ell(\theta;\xi_m)\|^2;\\
    &\EE_{Q,\xi_m}\left[\|Q(\theta;\xi_m)-\nabla\ell(\theta;\xi_m)\|^2\right]\leq\min\Big\{\frac{d}{(2^{b-1}-1)^2},\frac{\sqrt{d}}{2^{b-1}-1}\Big\}=:B'.
\end{align*}

Following the proof of Lemma \ref{lemma:lossdescent}, we can get
\begin{align*}
  \EE\cL(\theta^{k+1})\leq\EE\cL(\theta^k)-\left(\eta_k-\frac{L\eta_k^2}{2}\right)\EE\left[\|\nabla\cL(\theta^k)\|^2\right]+\left(\eta_k-L\eta_k^2\right)H_3+\frac{L\eta_k^2}{2}H_4,
\end{align*}
where
\begin{align*}
    H_3=&\sum\limits_{m\in{\cal M}}\EE\left[\dotp{\nabla\cL(\theta^k),\nabla\ell(\theta^k;\xi_m^k)-Q(\theta^{k-\htau_m^k};\xi_m^{k-\htau_m^k})}\right]\\
    =& H_1 +\sum\limits_{m\in{\cal M}}\EE\left[\dotp{\nabla\cL(\theta^k),\nabla\ell(\theta^{k-\htau_m^k};\xi_m^{k-\htau_m^k})-Q(\theta^{k-\htau_m^k};\xi_m^{k-\htau_m^k})}\right]\\
    \overset{\text{(3a)}}{\leq}&H_1+\frac{\sqrt{M}L}{12\eta_k}\sum\limits_{d=1}^D\EE\left[\|\theta^{k+1-d}-\theta^{k-d}\|^2\right]+\frac{6\sqrt{M}DL\eta_k}{2}B'\\
    \leq&\frac{L\eta_k}{2}\EE\left[\|\nabla\cL(\theta^k)\|^2\right]+\sum\limits_{d=1}^D\left(\frac{c_d}{2L\eta_k}+\frac{\sqrt{M}L}{6\eta_k}\right)\EE\left[\|\theta^{k+1-d}-\theta^{k-d}\|^2\right]+6\sqrt{M}DL\eta_k\sum\limits_{m\in{\cal M}}\left(\sigma_m^2+\frac{B'}{2}\right)
\end{align*}
where (3a) is obtained by an approach similar to \eqref{eqn:inequality1}, and
\begin{align*}
    H_4=&\EE\Big[\big\|\nabla\cL(\theta^k)-\sum\limits_{m\in{\cal M}}Q(\theta^{k-\htau_m^k};\xi_m^{k-\htau_m^k})\big\|^2\Big]\\
    \overset{\eqref{eqn:young1}}{\leq}&4\EE\Big[\big\|\sum\limits_{m\in{\cal M}}\nabla\ell(\theta^{k-\htau_m^k};\xi_m^{k-\htau_m^k})-Q(\theta^{k-\htau_m^k};\xi_m^{k-\htau_m^k})\big\|^2\Big]+\frac{4}{3}H_2\\
    \leq&4\sum\limits_{d=1}^Dc_d\EE\Big[\|\theta^{k+1-d}-\theta^{k-d}\|^2\Big]+12\sum\limits_{m\in{\cal M}}\left(\sigma_m^2+\frac{B'}{2}\right).
\end{align*}

Therefore,
\begin{align*}
    \EE\cL(\theta^{k+1})\leq\EE\cL(\theta^k)&-\left(\eta_k-L\eta_k^2+\frac{L^2\eta_k^3}{2}\right)\EE\left[\|\nabla\cL(\theta^k)\|^2\right]\\
    &+\sum\limits_{d=1}^D\left(\big(2L\eta_k^2+\frac{1}{2L}-\frac{\eta_k}{2}\big)c_d+(1-L\eta_k)\frac{\sqrt{M}L}{6}\right)\EE\left[\|\theta^{k+1-d}-\theta^{k-d}\|^2\right]\\
    &+L\eta_k^2\left(6+6\sqrt{M}D(1-L\eta_k)\right)\sum\limits_{m\in{\cal M}}\left(\sigma_m^2+\frac{B'}{2}\right).
\end{align*}

Similar to Lemma \ref{lemma:lyapunovdescent},
\begin{align*}
    V^{k+1}-V^k\leq&-\underbrace{\left(\eta_k-(L+2\gamma_1)\eta_k^2\right)}_{B_0^k}\EE\left[\|\nabla\cL(\theta^k)\|^2\right]\\
    &+\!\sum\limits_{d=1}^{D-1}\!\underbrace{\left((\frac{3}{2}\eta_k\!+\!\frac{1}{2L})c_d\!+\!\frac{\sqrt{M}L}{6}\!+\!8c_d\gamma_1\eta_k^2\!+\!\gamma_{d+1}\!-\!\gamma_d\right)}_{A_d^k}\EE\left[\|\theta^{k+1-d}-\theta^{k-d}\|^2\right]\\
    &+\underbrace{\left((\frac{3}{2}\eta_k+\frac{1}{2L})c_D+\frac{\sqrt{M}L}{6}+8c_D\gamma_1\eta_k^2-\gamma_D\right)}_{A_D^k}\EE\left[\|\theta^{k+1-D}-\theta^{k-D}\|^2\right]\\
    &+\!\underbrace{\left((6+6\sqrt{M} D)L\eta_k^2\!+\!24\gamma_1\eta_k^2\right)}_{B_1^k}\!\sum\limits_{m\in{\cal M}}\!\left(\sigma_m^2\!+\!\frac{B'}{2}\right).
\end{align*}

Let $\bar{\eta}\geq\eta_k$ and choose $\gamma_d$ such that
\begin{align*}
    &(\frac{3}{2}\bar{\eta}+\frac{1}{2L})c_d+\frac{\sqrt{M}L}{6}+8c_d\gamma_1\bar{\eta}^2+\gamma_{d+1}-\gamma_d=0,\quad 1\leq d\leq D-1\\
    &(\frac{3}{2}\bar{\eta}+\frac{1}{2L})c_D+\frac{\sqrt{M}L}{6}+8c_D\gamma_1\bar{\eta}^2-\gamma_D=0.
\end{align*}

Solve the linear equations above and get
\begin{align*}
    \gamma_1=\frac{(\frac{3}{2}\bar{\eta}+\frac{1}{2L})\sum\limits_{d=1}^Dc_d+\frac{\sqrt{M}DL}{6}}{1-8\sum\limits_{d=1}^Dc_d\bar{\eta}^2}.
\end{align*}

Select $c_d\leq\min\{\frac{1}{16D\bar{\eta}^2},\frac{\sqrt{M}L^2}{24}\}$ such that $\gamma_1\leq \frac{\sqrt{M}DL}{2}$.
Set $\eta_k\leq\frac{1}{2L+2\sqrt{M}DL}\leq\frac{1}{2L+4\gamma_1}$. By taking sum $\sum\limits_{k=0}^{K-1} V^{k+1}-V^k$ it follows that
\begin{align*}
    &\sum\limits_{k=0}^{K-1}\frac{\eta_k}{2}\EE\Big[\|\nabla\cL(\theta^k)\|^2\Big]\leq\cL(\theta^0)-\cL^{\star}+\sum\limits_{k=0}^{K-1}(6+18\sqrt{M}D)\eta_k^2\sum\limits_{m\in{\cal M}}\left(\sigma_m^2\!+\frac{B'}{2}\right).
\end{align*}

Define a random variable $R(K)$ with $\mathbb{P}(R(K)=j)=\frac{\eta_j}{\sum\limits_{k=0}^{K-1}\eta_k}$ and then
\begin{align*}
    \EE\Big[\|\nabla\cL(\theta^{R(K)})\|^2\Big]\leq\frac{\cL(\theta^0)-\cL^{\star}+\sum\limits_{k=0}^{K-1}\left(6+18\sqrt{M}D\right)L\eta_k^2\sum\limits_{m\in{\cal M}}\left(\sigma_m^2+\frac{B'}{2}\right)}{\sum\limits_{k=0}^{K-1}\frac{\eta_k}{2}}.
\end{align*}
If we choose constant stepsize
\begin{align*}
    \eta_k=\min\left\{\frac{1}{2L+2\sqrt{M}DL}, \frac{c_\eta}{\sqrt{K}}\right\},
\end{align*}
then
{\small
\begin{align*}
    \frac{1}{K}\sum\limits_{k=0}^{K-1}\EE\Big[\|\nabla\cL(\theta^{k})\|^2\Big]\leq&\frac{\cL(\theta^0)-\cL^{\star}+K\left(6+18\sqrt{M}D\right)L\eta_0^2\sum\limits_{m\in{\cal M}}\sigma_m^2}{\frac{1}{2}K\eta_0}\\
    \leq&\frac{4(L+\sqrt{M}DL)(\cL(\theta^0)-\cL^{\star})}{K}+\frac{2(\cL(\theta^0)-\cL^{\star})}{c_{\eta}\sqrt{K}}+\frac{c_{\eta}\left(12+36\sqrt{M}D\right)L\sum\limits_{m\in{\cal M}}\left(\sigma_m^2+\frac{B'}{2}\right)}{\sqrt{K}}.
\end{align*}}

\section{Additional numerical results}
The additional numerical results in this section include both homogeneous and heterogeneous setting.\\  
\textbf{Homogeneous}: Data samples are shuffled and uniformly partitioned to $M$ workers.\\
 \textbf{Heterogeneous}: Data samples with same labels are unevenly partitioned and assigned to $M$ workers.


\begin{table*}[h!]
\centering
\resizebox{\textwidth}{!}{
\begin{tabular}{|c|c|c|c|c|c|c|c|c|c|c|c|c|}
\hline
\multirow{2}{*}{Alg} & \multirow{2}{*}{SGD} & \multirow{2}{*}{\!\!\!LAG-WK\!\!\!} & \multicolumn{4}{|c|}{LASG} &\multirow{2}{*}{\!\!\!local SGD\!\!\!} &\multirow{2}{*}{QSGD} & \multicolumn{4}{|c|}{LAQSG}\\ \cline{4-7}\cline{10-13}
                     &                      &                         & WK1 & WK2 & PS & PSE
&                           &                      & WK1 & WK2 & PS & PSE \\ \hline
ijcnn1(LR)& 0.4276 & 0.3352 & 0.2252 & 0.2254 & 0.2675 &{\bf 0.2246}& 0.2460 & 0.4278 & 0.2253 &{\bf 0.2247}& 0.2660 &{\bf 0.2250}\\ \hline
mnist(LR) & 0.4219 & 0.4205 & 0.1894 &{\bf 0.1372}& 0.2940 & 0.1690 & 0.1979 & 0.4213 &{\bf 0.1615}&{\bf 0.1319}& 0.2948 & 0.1662 \\ \hline
covtype(LR)&0.6761 & 0.6759 &{\bf 0.5725}&{\bf 0.5653}& 0.6477 & 0.5952 & 0.5885 & 0.6764 & 0.5856 &{\bf 0.5652}& 0.6486 & 0.5915 \\ \hline
mnist(NN) & 0.1612 & 0.1508 & 0.1068 &{\bf 0.0395}& --- &{\bf 0.0416}&{\bf 0.0868}& 0.2388 & 0.2963 & 0.1545 & --- & 0.0950\\ \hline
\end{tabular}}
\caption{Objective value after 1000 and 10000 communication rounds for logistic regression (LR) and neural network (NN) respectively in the homogeneous setting.}
\label{table: 3}
\end{table*}

\begin{table*}[h!]
\centering
\resizebox{\textwidth}{!}{
\begin{tabular}{|c|c|c|c|c|c|c|c|c|c|c|c|c|}
\hline
\multirow{2}{*}{Alg} & \multirow{2}{*}{SGD} & \multirow{2}{*}{\!\!\!LAG-WK\!\!\!} & \multicolumn{4}{|c|}{LASG} &\multirow{2}{*}{\!\!\!local SGD\!\!\!} &\multirow{2}{*}{QSGD} & \multicolumn{4}{|c|}{LAQSG}\\ \cline{4-7}\cline{10-13}
                     &                      &                         & WK1 & WK2 & PS & PSE
&                           &                      & WK1 & WK2 & PS & PSE \\ \hline
ijcnn1(LR)& 0.6331 & 0.3605 & 0.2911 & 0.2920 & 0.3257 & 0.2930 & 0.3085 & 0.4498 &{\bf 0.2297}&{\bf 0.2296}& 0.2726 &{\bf 0.2296}\\ \hline
mnist(LR) & 0.6875 & 0.6813 & 0.5375 & 0.4597 & 0.6760 & 0.3334 & 0.5640 & 0.5625 &{\bf 0.2837}&{\bf 0.1710}& 0.4385 &{\bf 0.2165}\\ \hline
covtype(LR)&0.6831 & 0.6826 & 0.5923 & 0.5832 & 0.6635 & 0.6148 & 0.6081 & 0.6436 &{\bf 0.5434}&{\bf 0.5346}& 0.5977 &{\bf 0.5519}\\ \hline
mnist(NN) & 0.2965 & 0.3029 & 0.2142 &{\bf 0.0902}& --- &{\bf 0.0746}& 0.1751 & 0.1857 & 0.1623 & 0.1111 &--- &{\bf 0.0663} \\ \hline
\end{tabular}}
\caption{Objective value after 1e5, 1e6, 1e6 bits of uploads for logistic regression on ijcnn1, mnist, covtype,  and 1e8 for neural network in the homogeneous setting.}
\label{table: 4}
\end{table*}

\begin{figure*}[ht]
    \centering
    \hspace*{-2ex}
    \includegraphics[width=.33\textwidth]{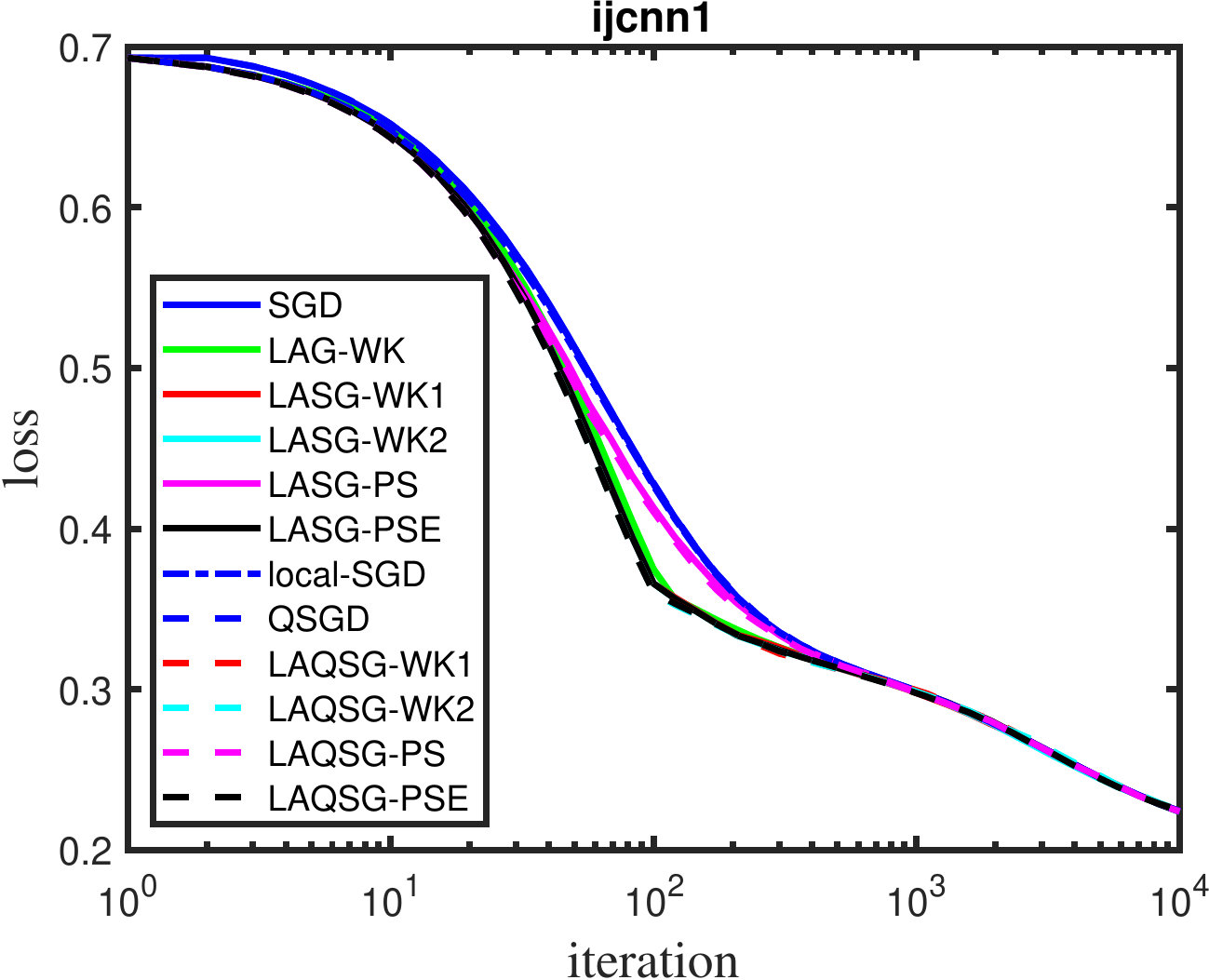}
    \hspace*{-1ex}
    \includegraphics[width=.33\textwidth]{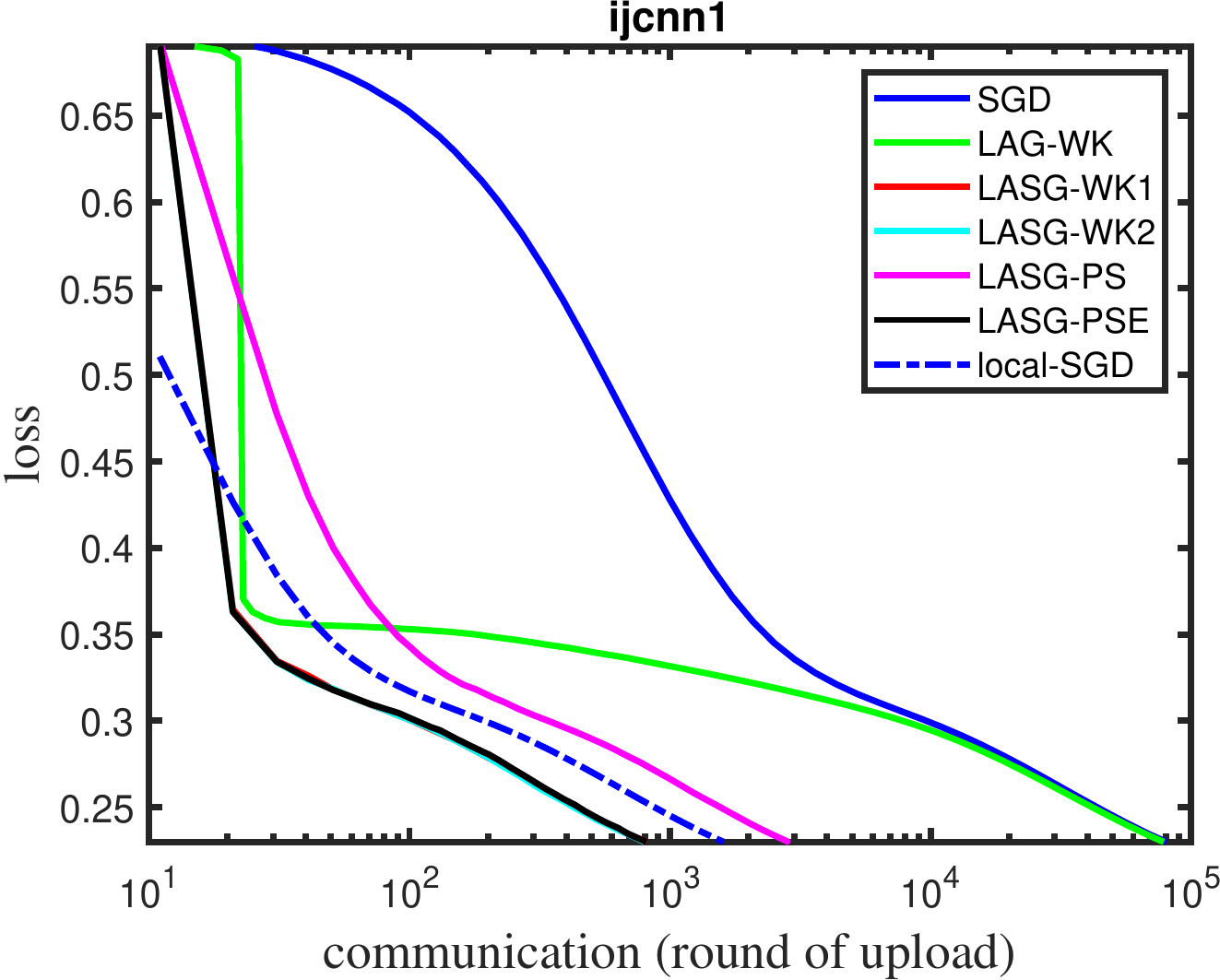}
    \hspace*{-1ex}
    \includegraphics[width=.33\textwidth]{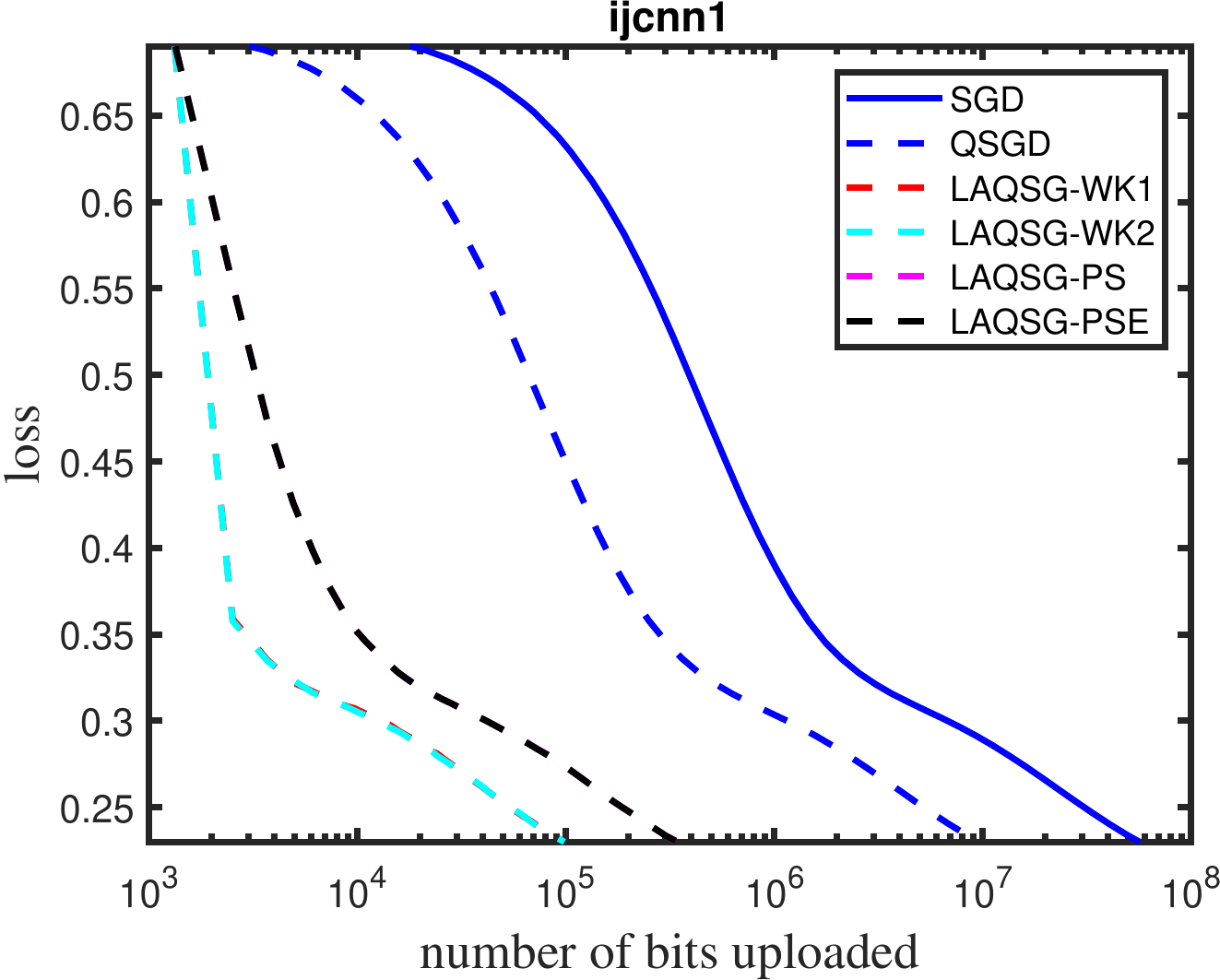}\\
    \caption{Logistic regression on \textit{ijcnn1} in the homogeneous setting.}
    \label{fig:ijcnn1_u}
\end{figure*}

\begin{figure*}[ht]
    \centering
    \hspace*{-2ex}
    \includegraphics[width=.33\textwidth]{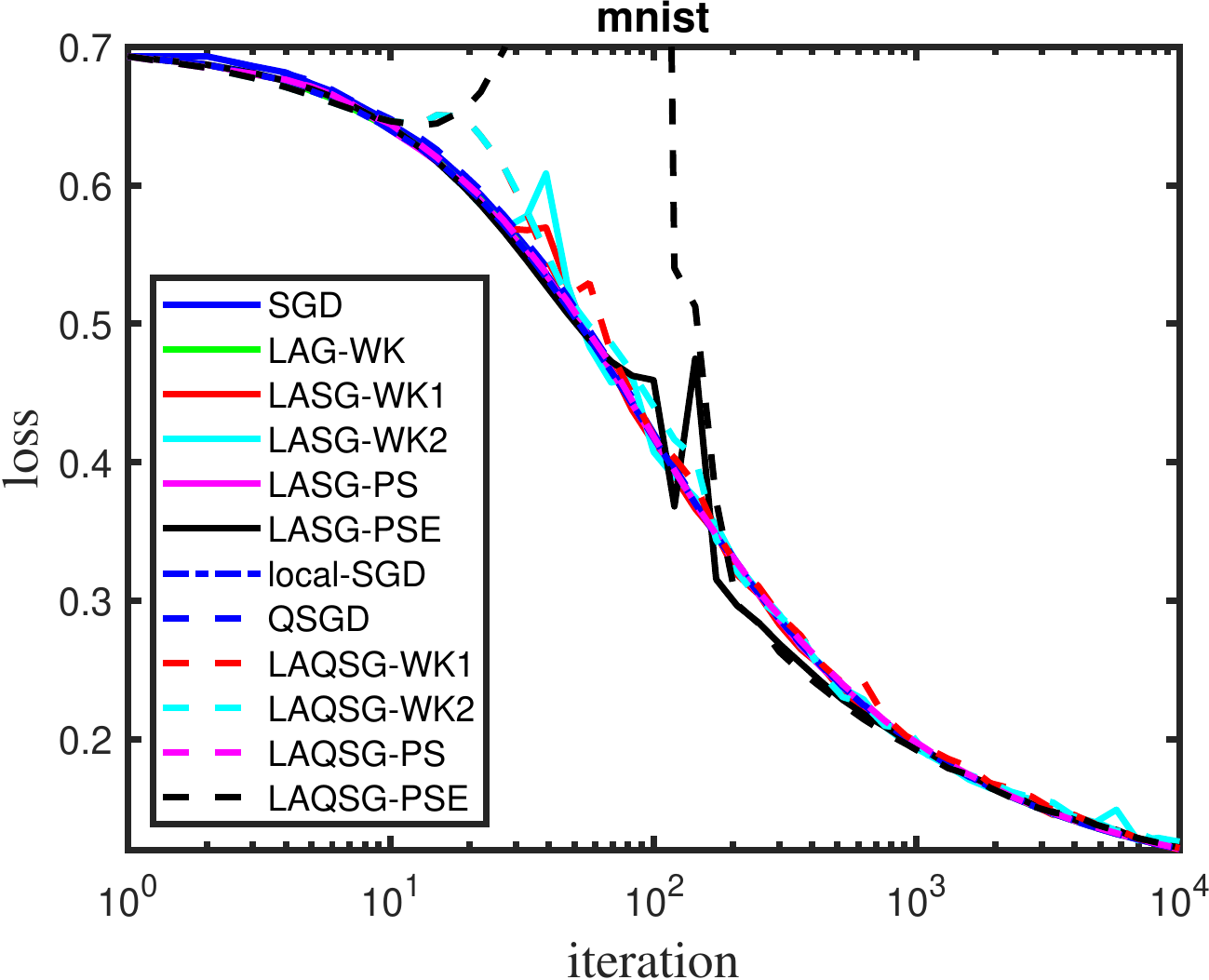}
    \hspace*{-0.5ex}
    \includegraphics[width=.33\textwidth]{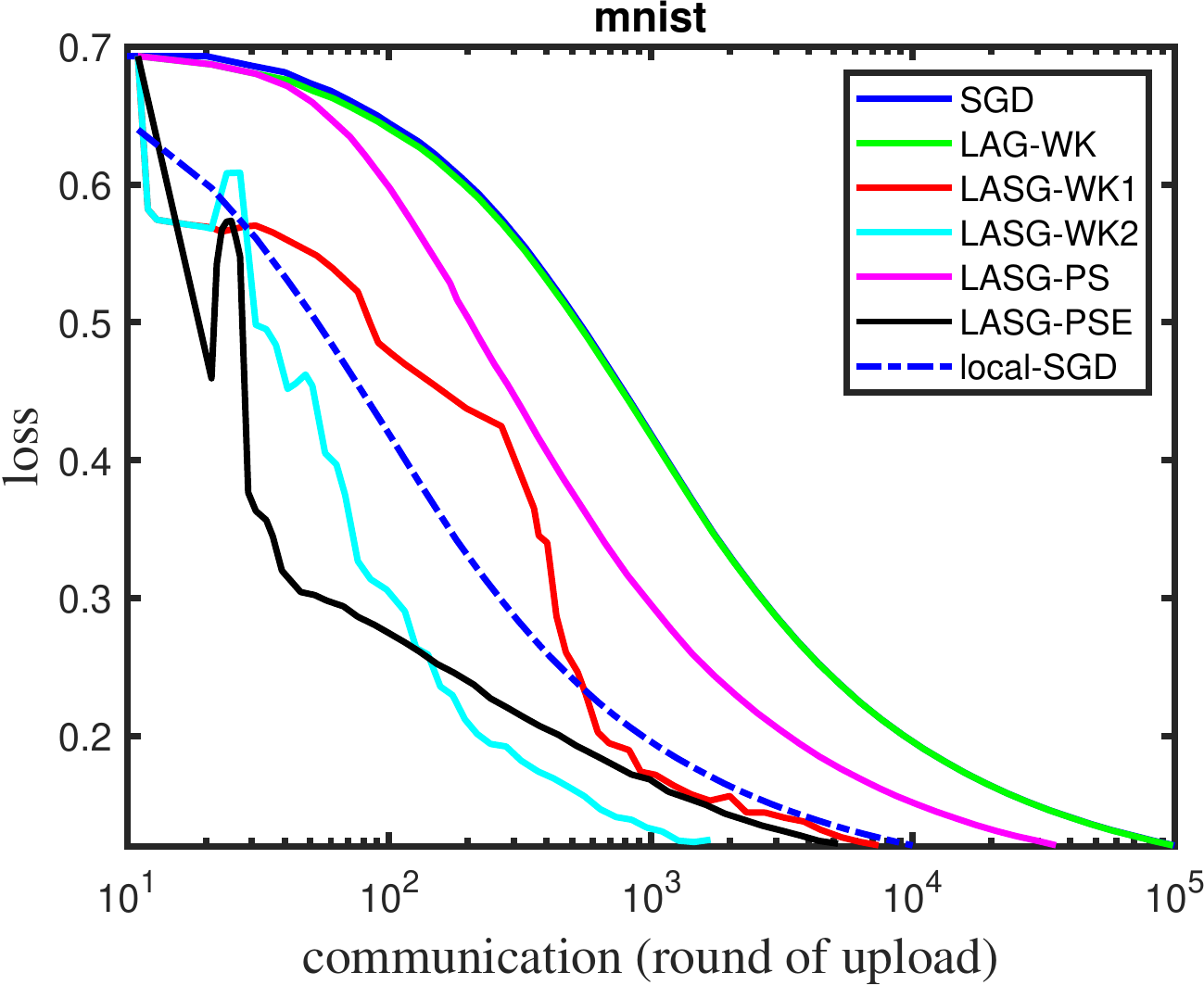}
    \hspace*{-0.5ex}
    \includegraphics[width=.33\textwidth]{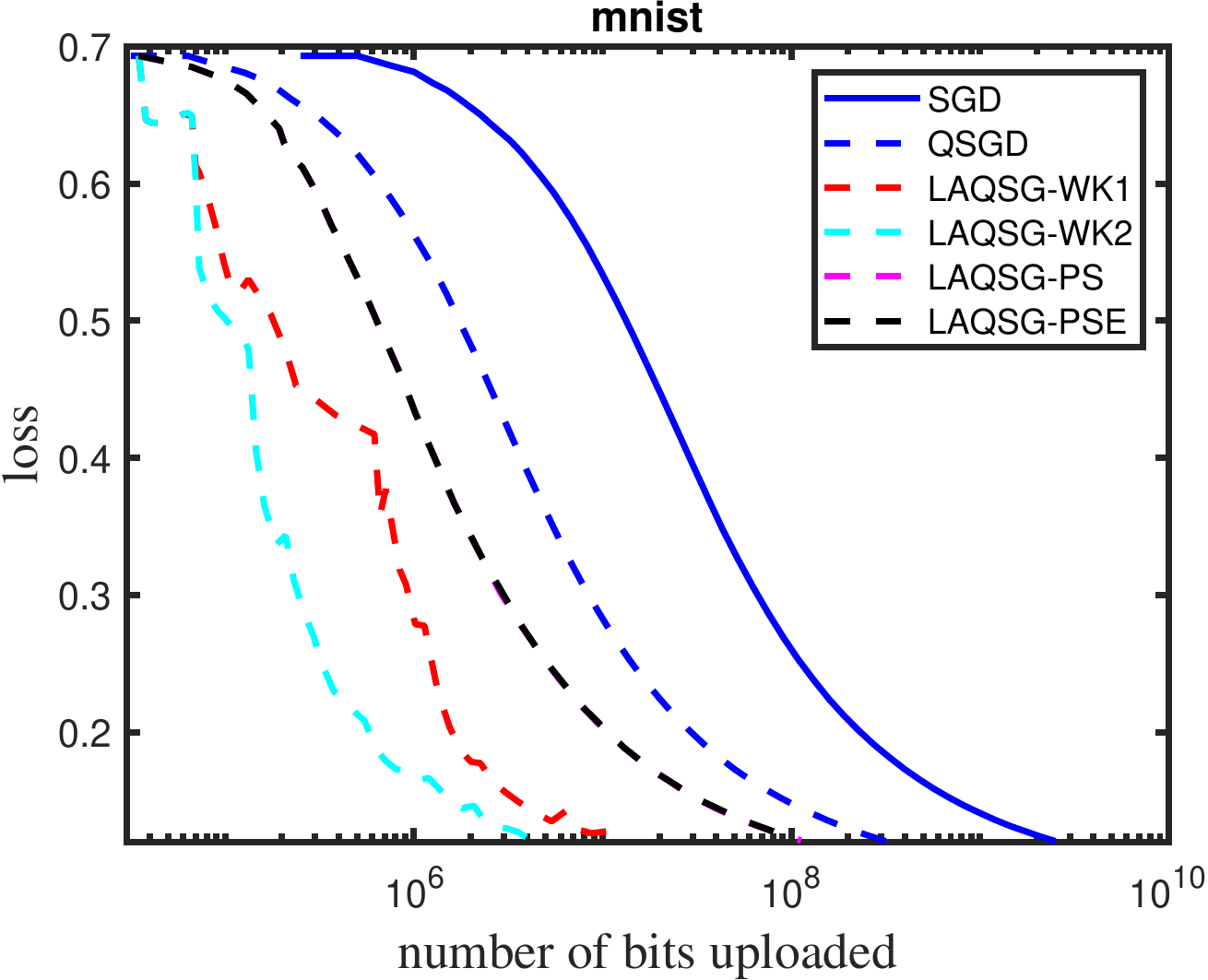}
    \caption{Logistic regression on \textit{mnist} digits 3 and 5 in the homogeneous setting}
    \label{fig:mnist_u}
\end{figure*}

\begin{figure*}[ht]
    \centering
    \hspace*{-2ex}
    \includegraphics[width=.33\textwidth]{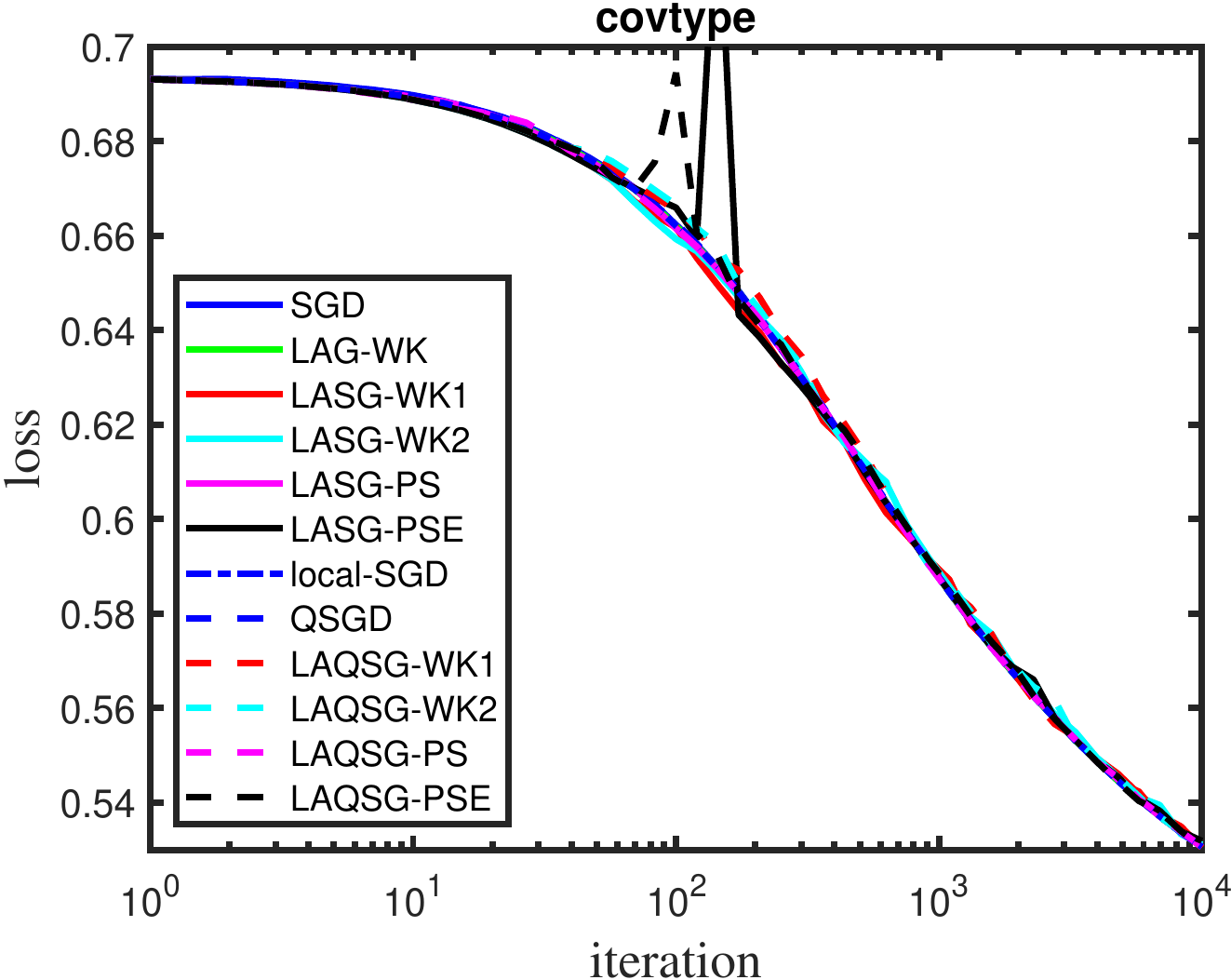}
    \hspace*{-1ex}
    \includegraphics[width=.33\textwidth]{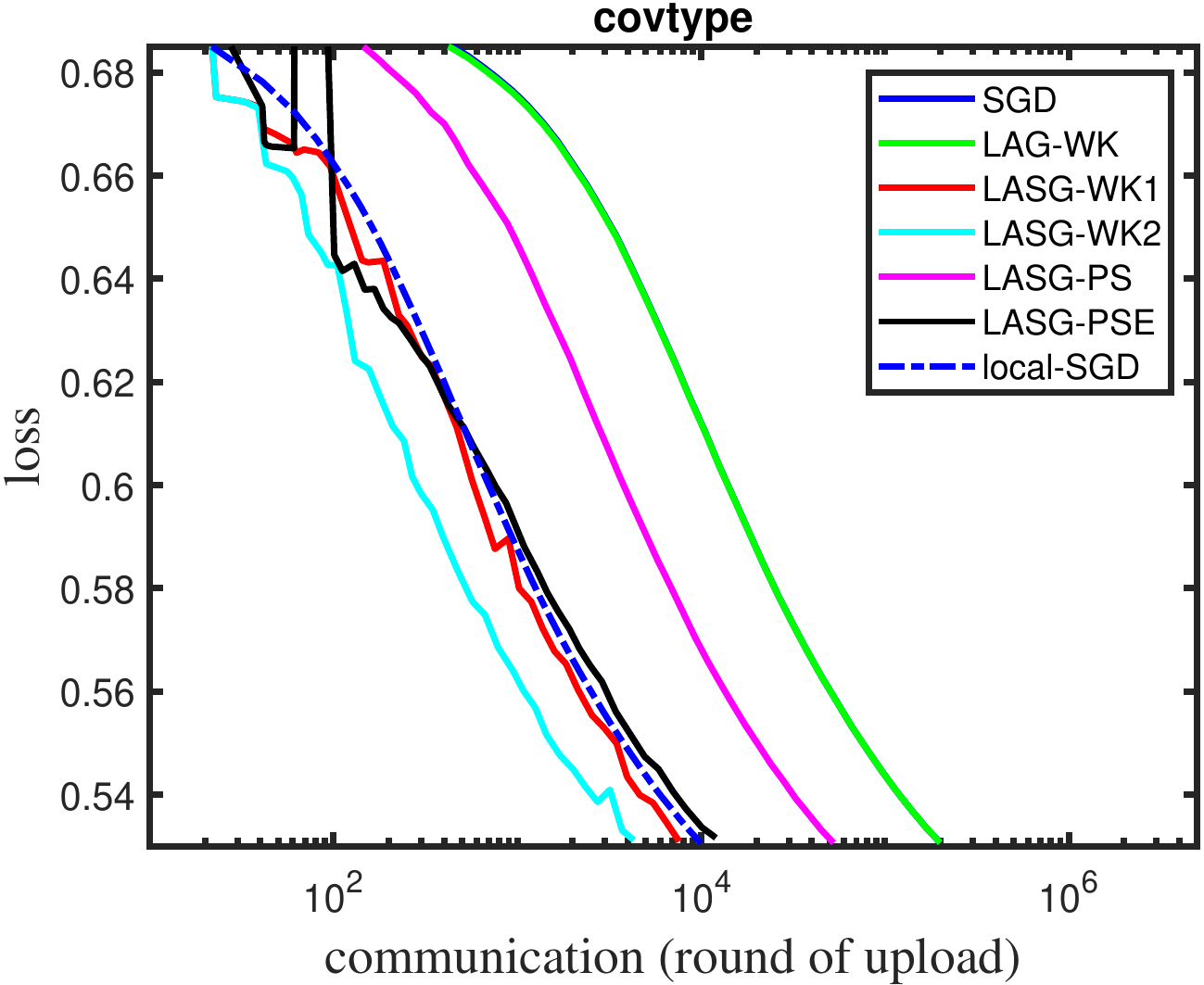}
    \hspace*{-1ex}
    \includegraphics[width=.33\textwidth]{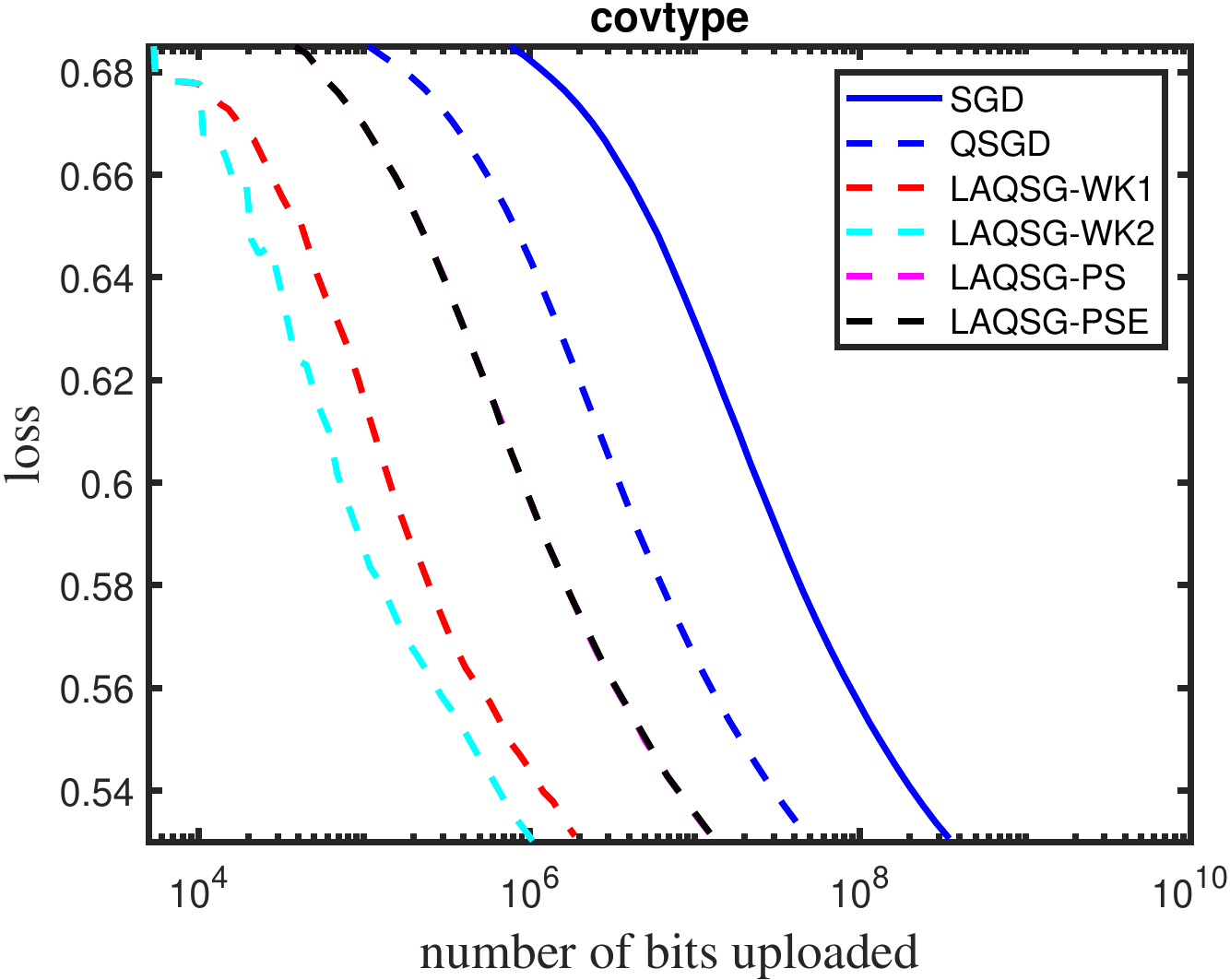}\\
    \caption{Logistic regression on dataset \textit{covtype} in the homogeneous setting}
    \label{fig:covtype_u}
\end{figure*}

\begin{figure*}[ht]
    \centering
    \hspace*{-2ex}
    \includegraphics[width=.33\textwidth]{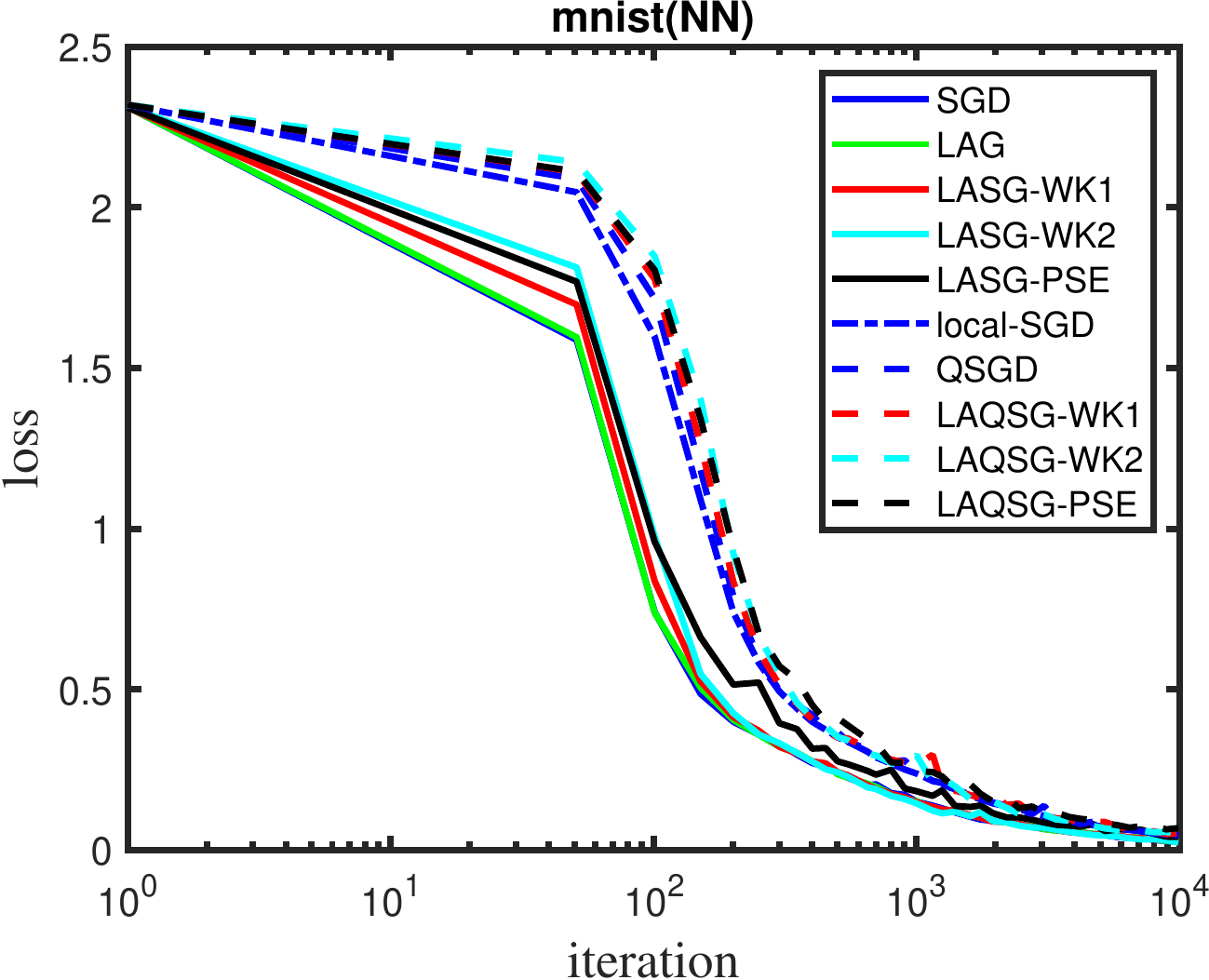}
    \hspace*{-1ex}
    \includegraphics[width=.33\textwidth]{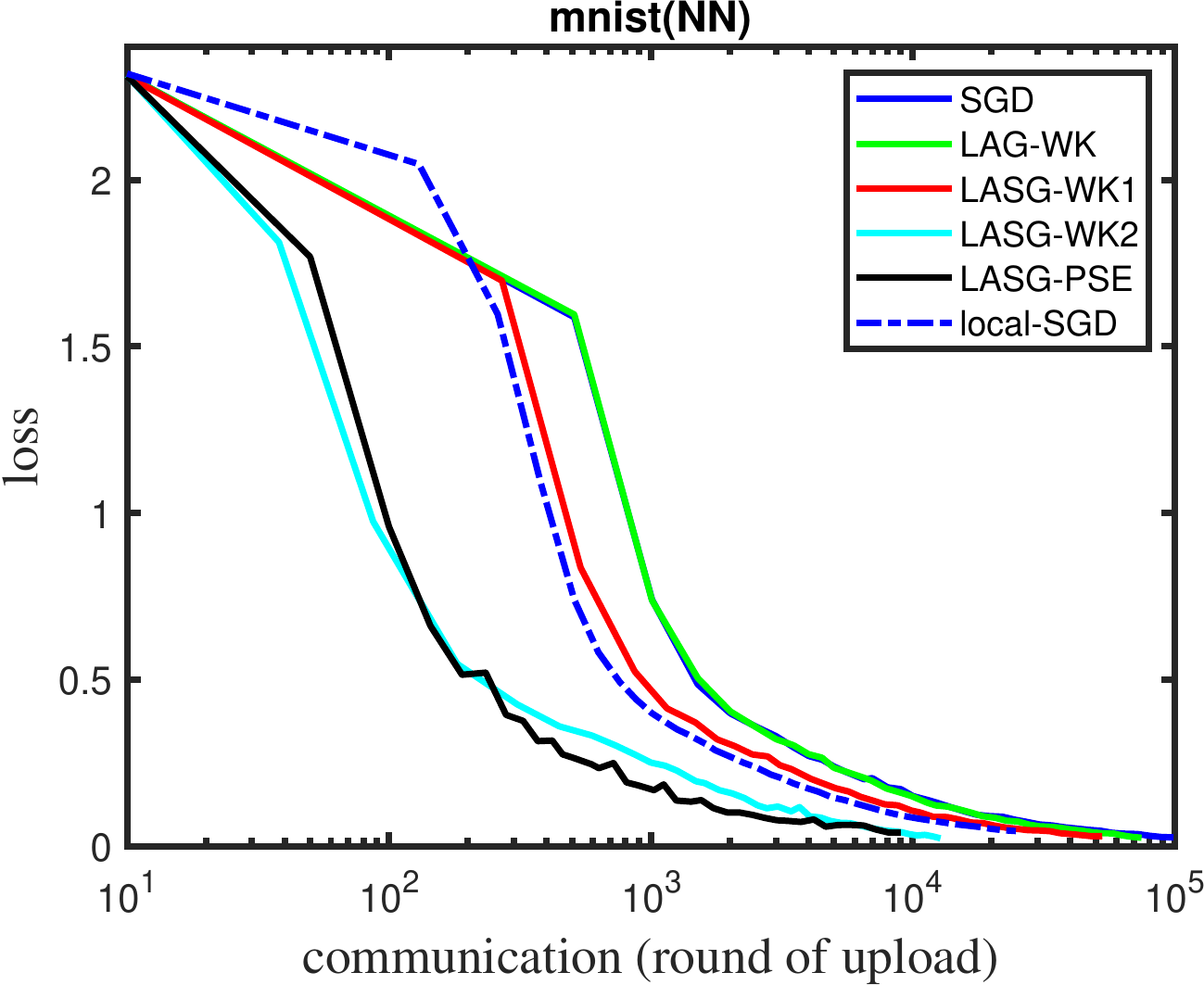}
    \hspace*{-1ex}
    \includegraphics[width=.33\textwidth]{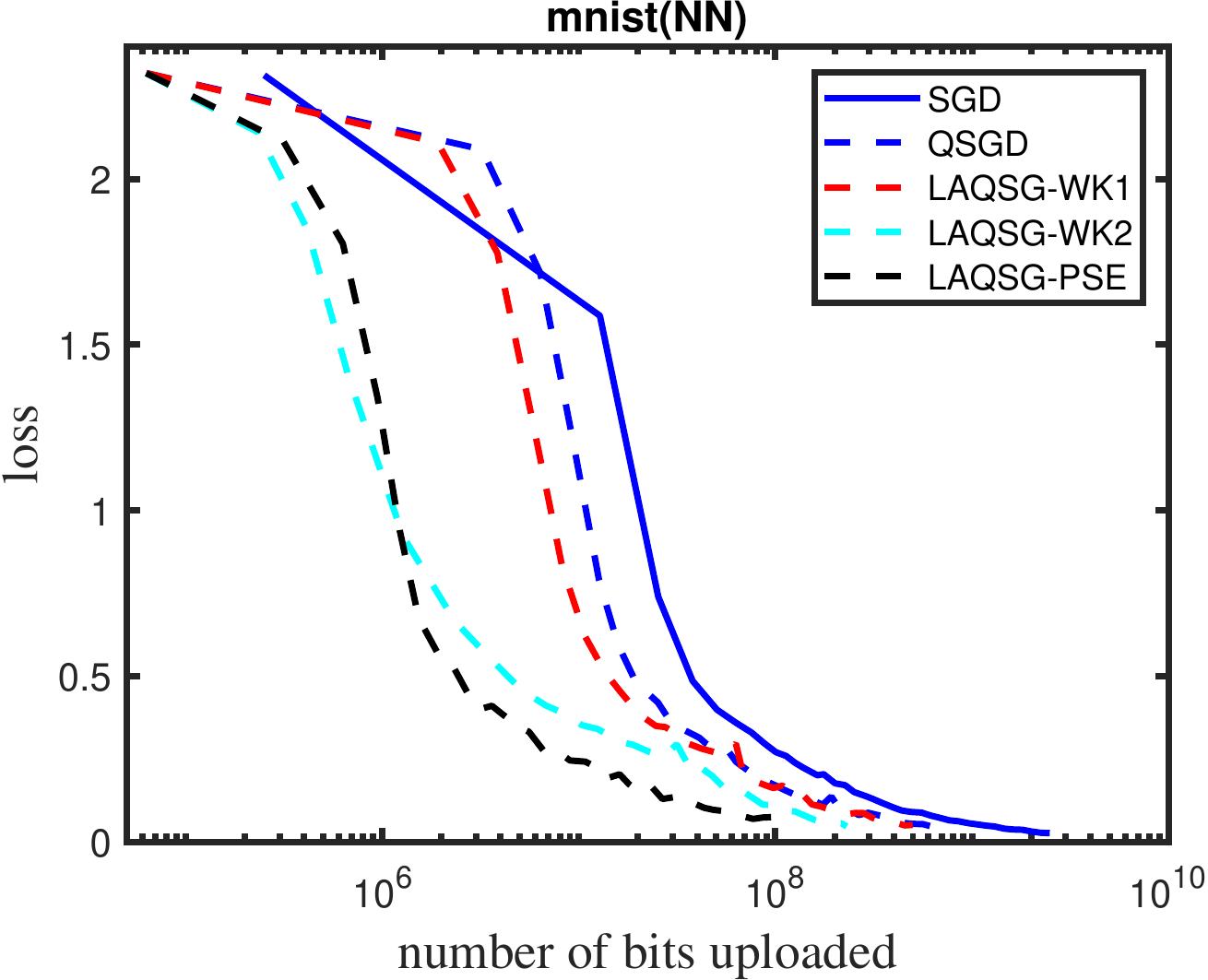}
    \caption{Neural network on dataset \textit{mnist} in the homogeneous setting}
    \label{fig:NNmnist_u}
\end{figure*}

\begin{figure*}
    \centering
    \hspace*{-2ex}
    \includegraphics[width=.33\textwidth]{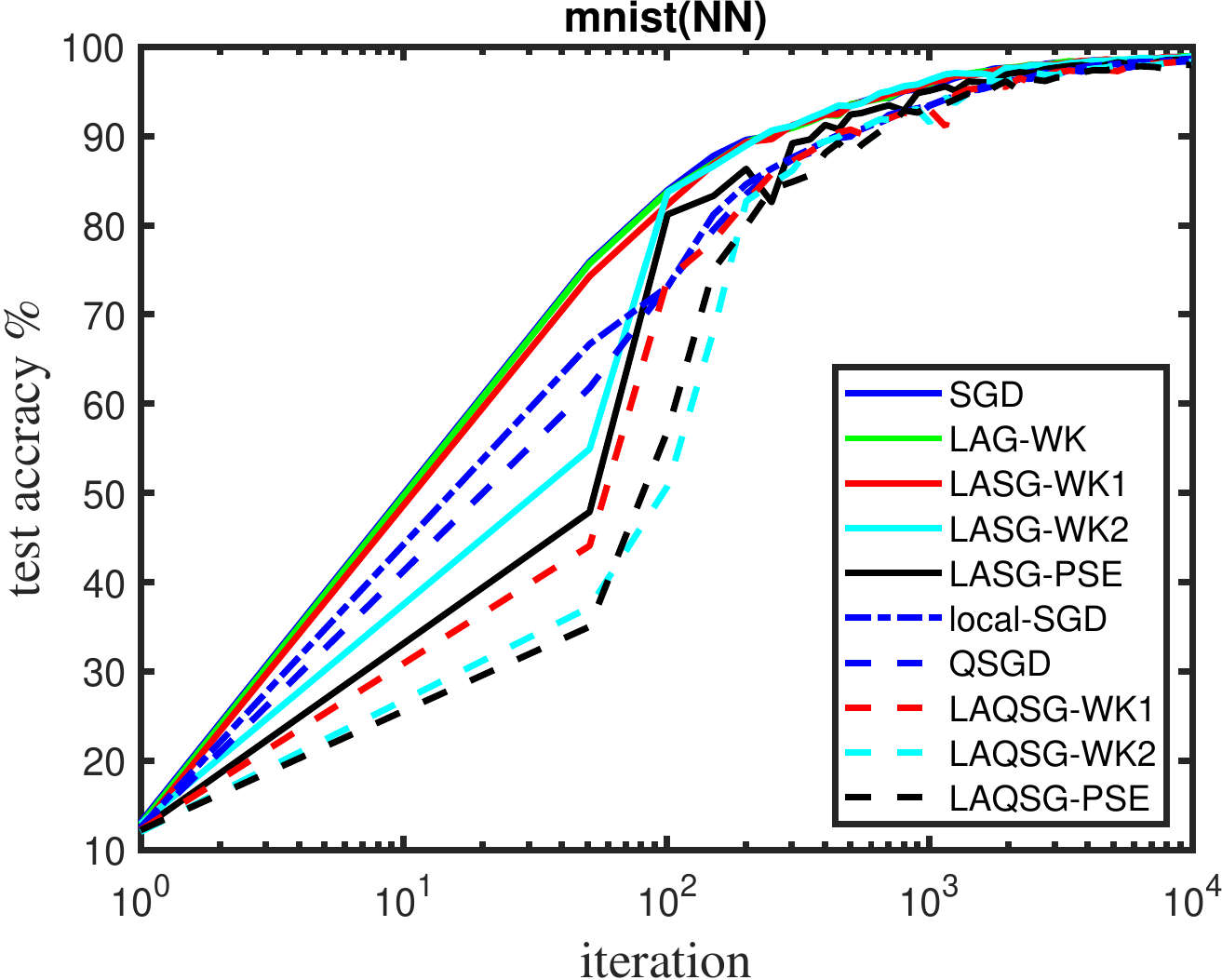}
    \hspace*{-1ex}
    \includegraphics[width=.33\textwidth]{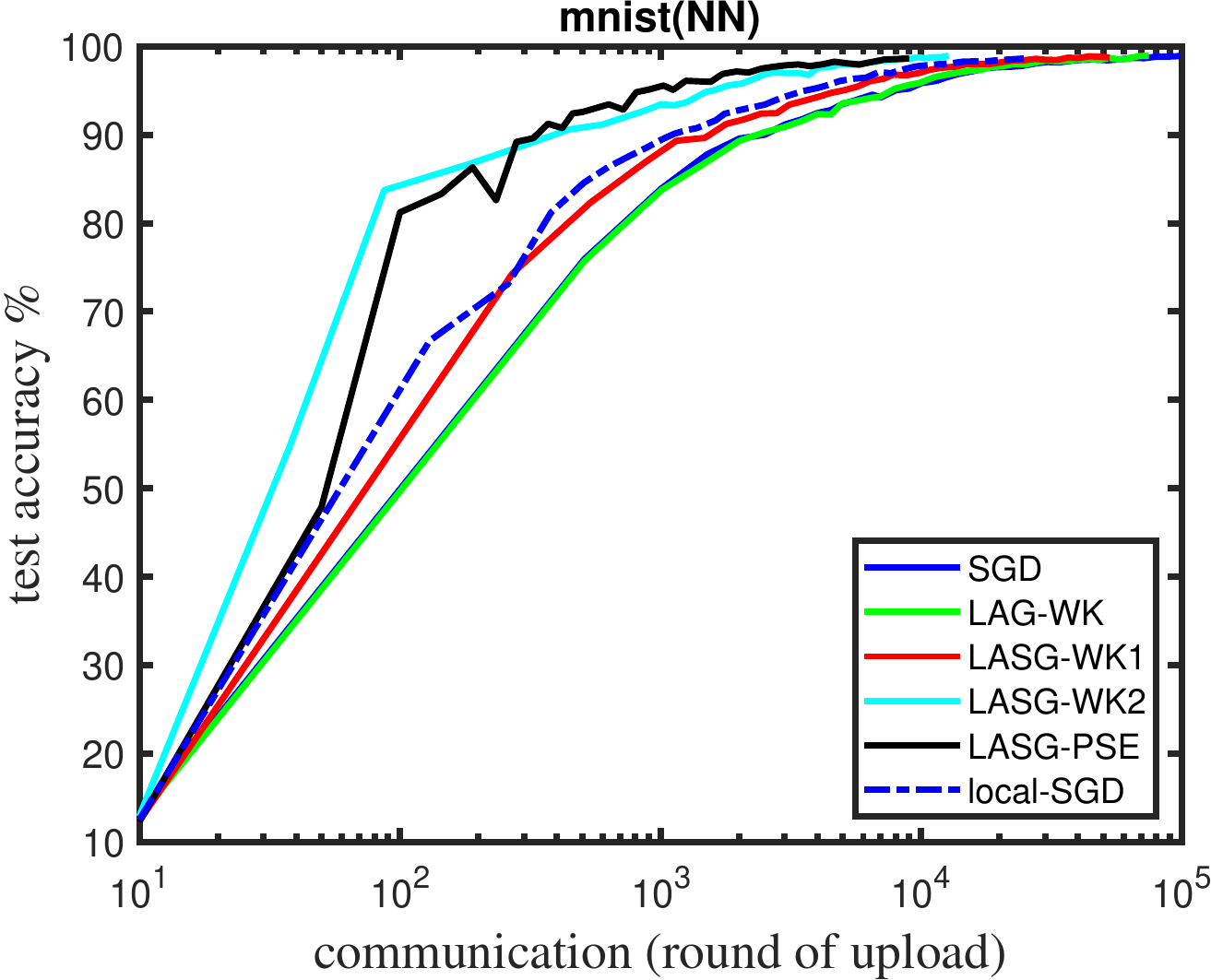}
    \hspace*{-1ex}
    \includegraphics[width=.33\textwidth]{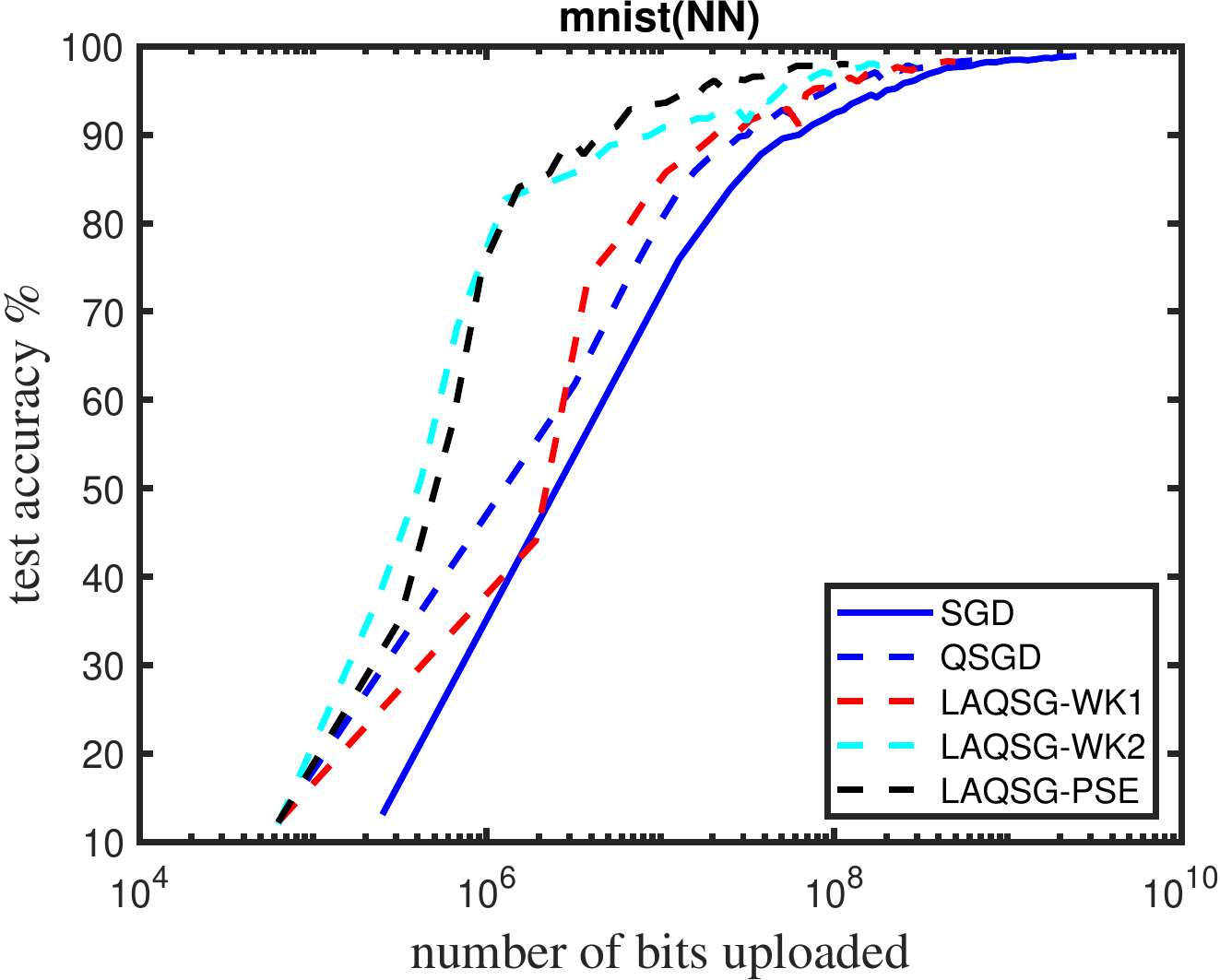}
    \caption{Neural network on dataset \textit{mnist} in the homogeneous setting}
    \label{fig:NNmnist_test_u}
\end{figure*}

\begin{figure*}
    \centering
    \hspace*{-2ex}
    \includegraphics[width=.33\textwidth]{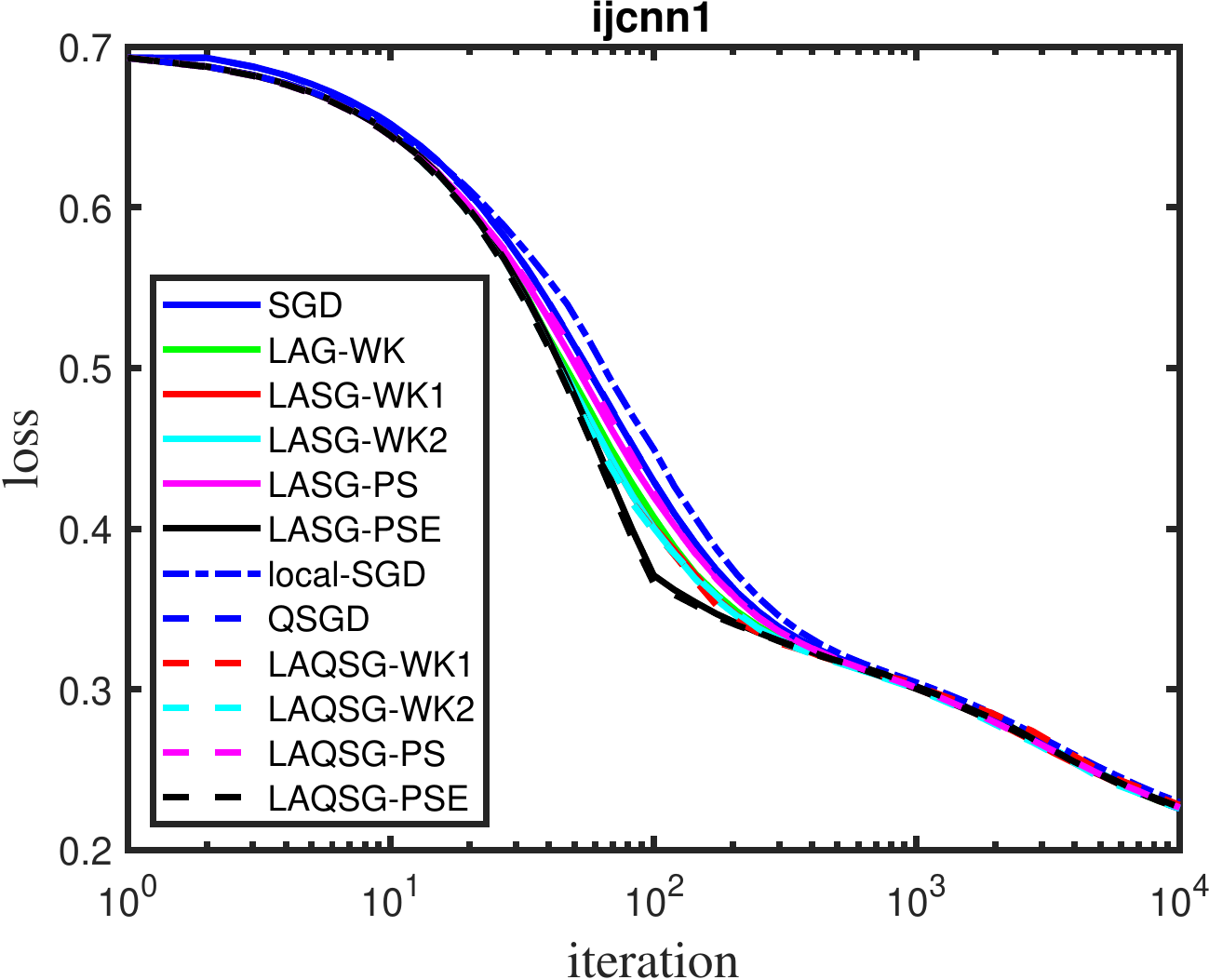}
    \hspace*{-1ex}
    \includegraphics[width=.33\textwidth]{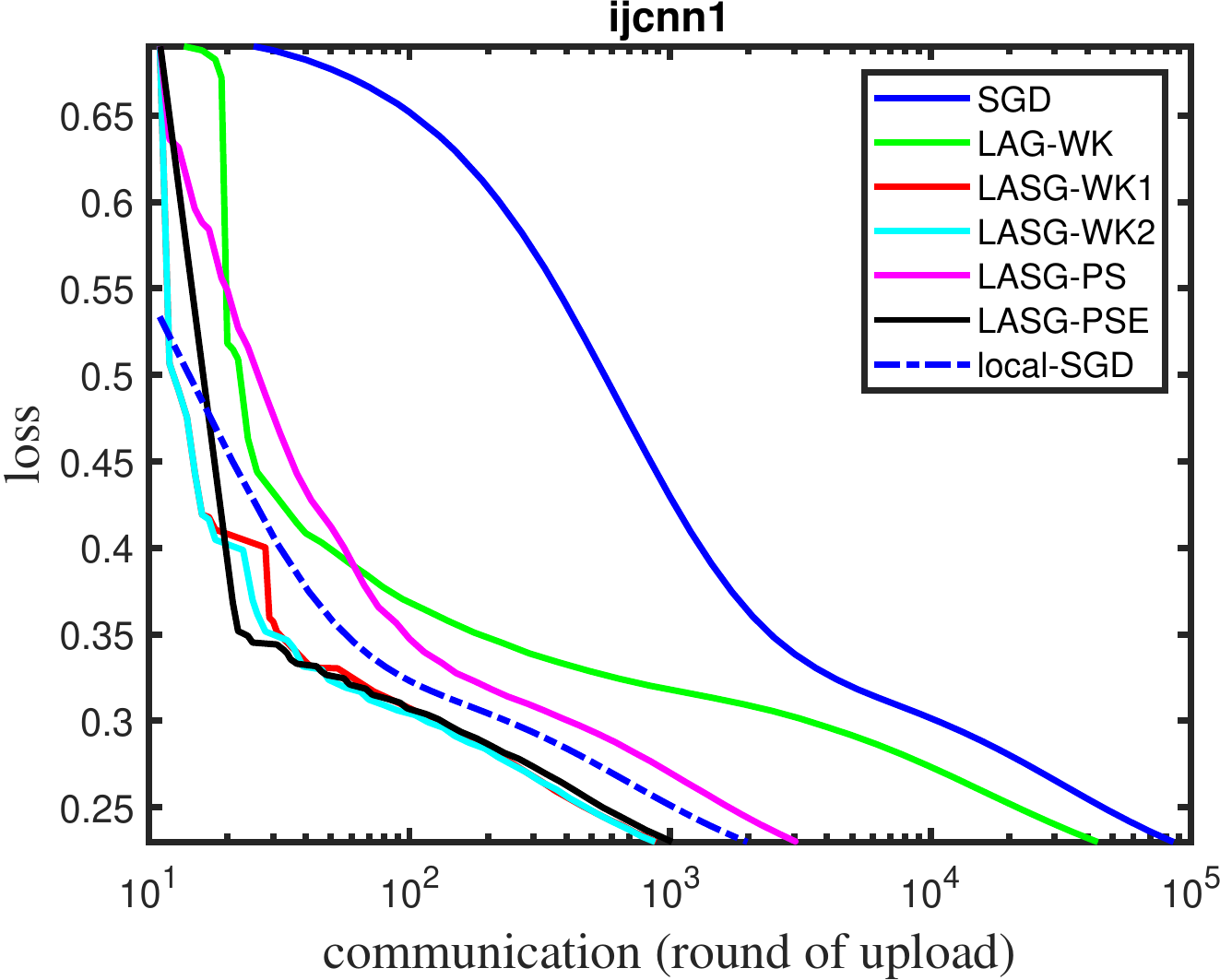}
    \hspace*{-1ex}
    \includegraphics[width=.33\textwidth]{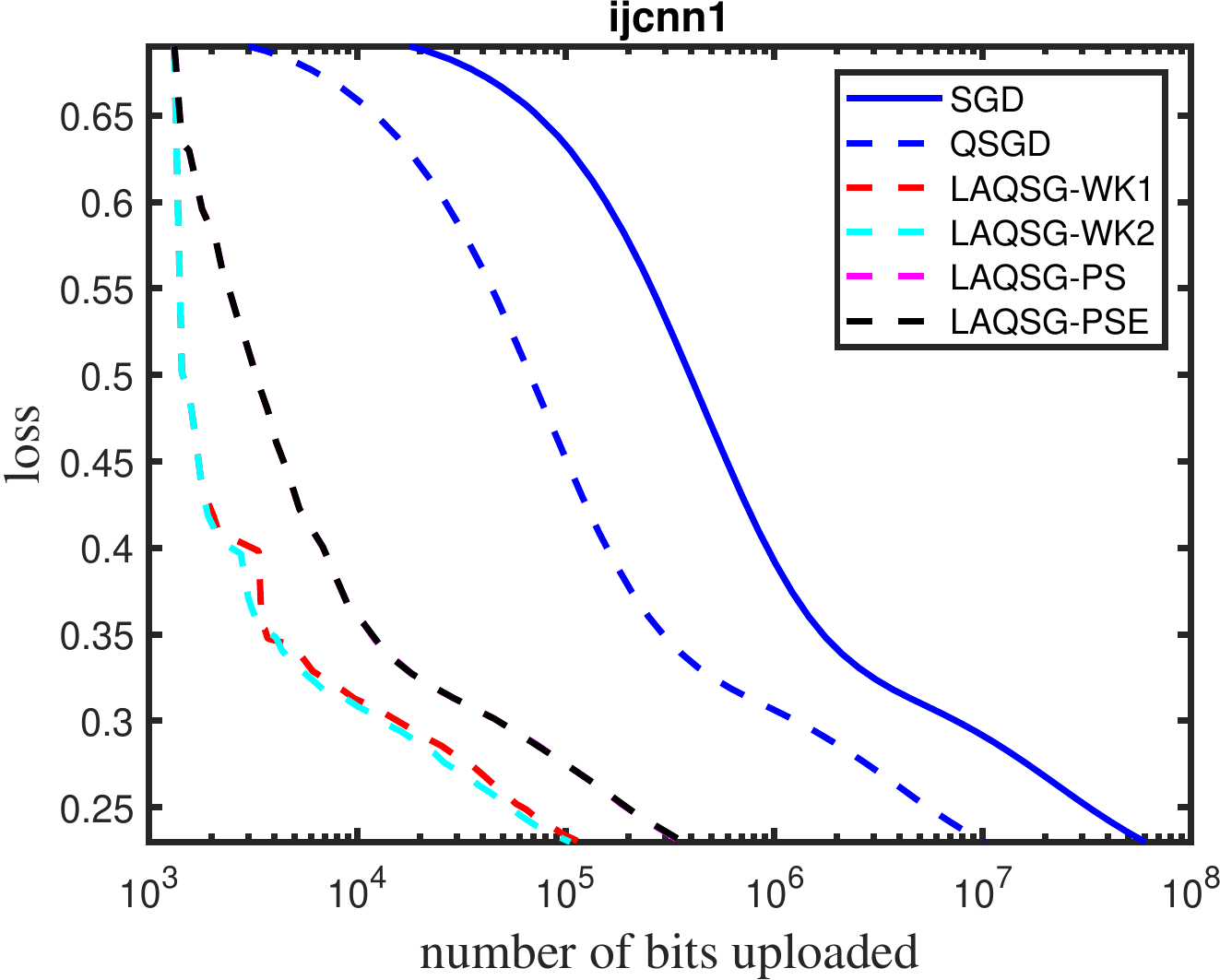}\\
    \caption{Logistic regression on \textit{ijcnn1} in the heterogeneous setting}
    \label{fig:ijcnn1}
\end{figure*}

\begin{figure*}[ht]
    \centering
    \hspace*{-2ex}
    \includegraphics[width=.33\textwidth]{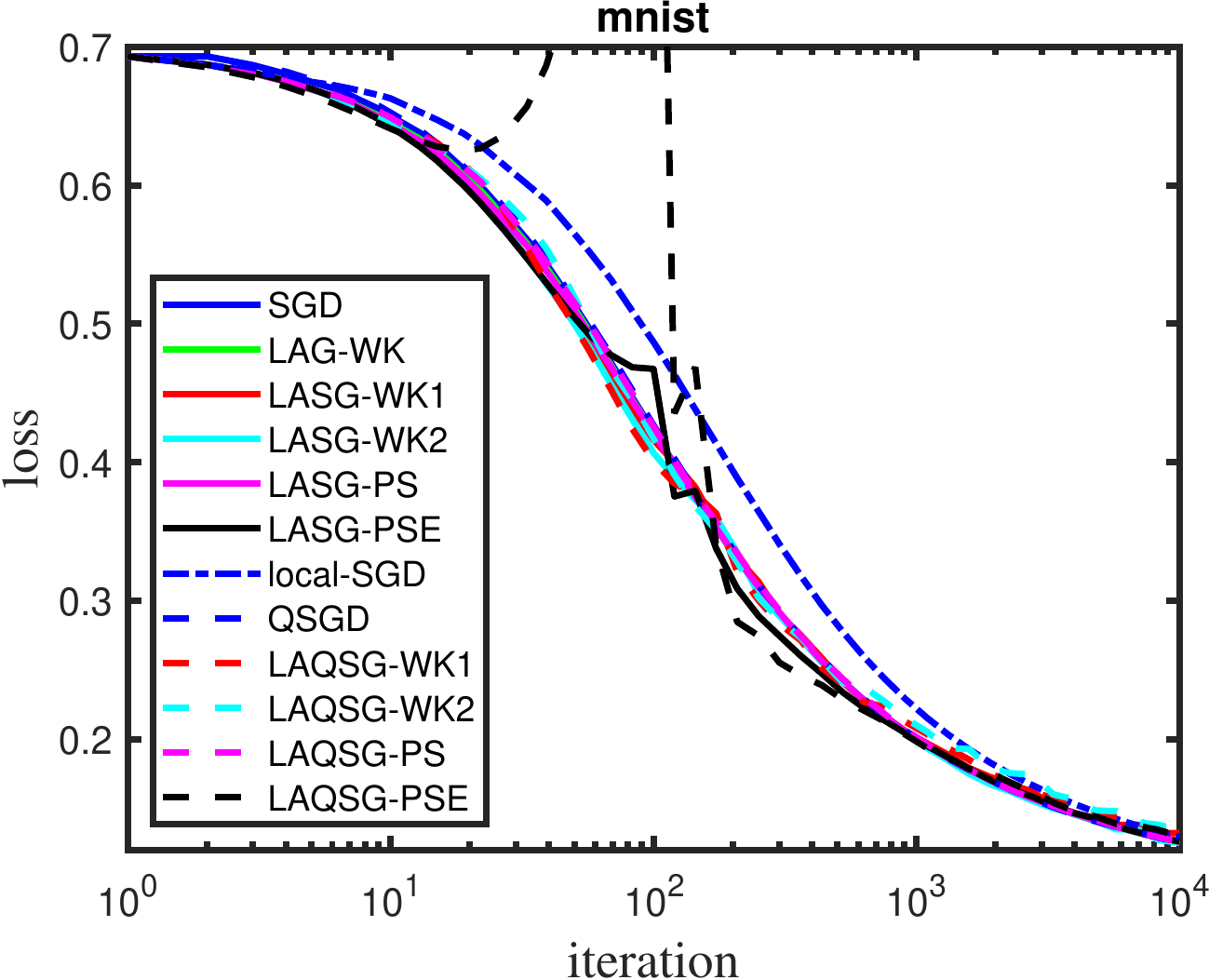}
    \hspace*{-0.5ex}
    \includegraphics[width=.33\textwidth]{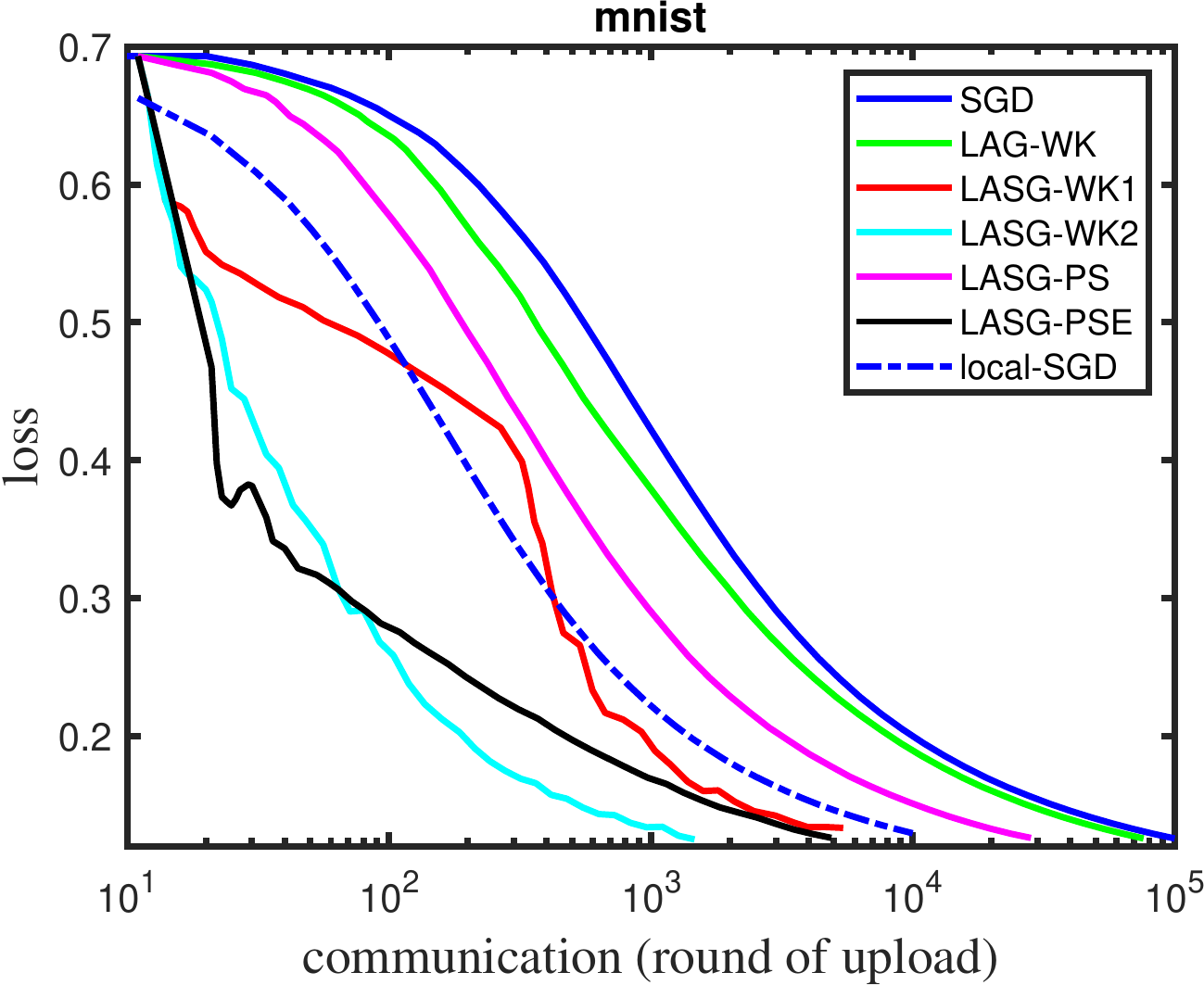}
    \hspace*{-0.5ex}
    \includegraphics[width=.33\textwidth]{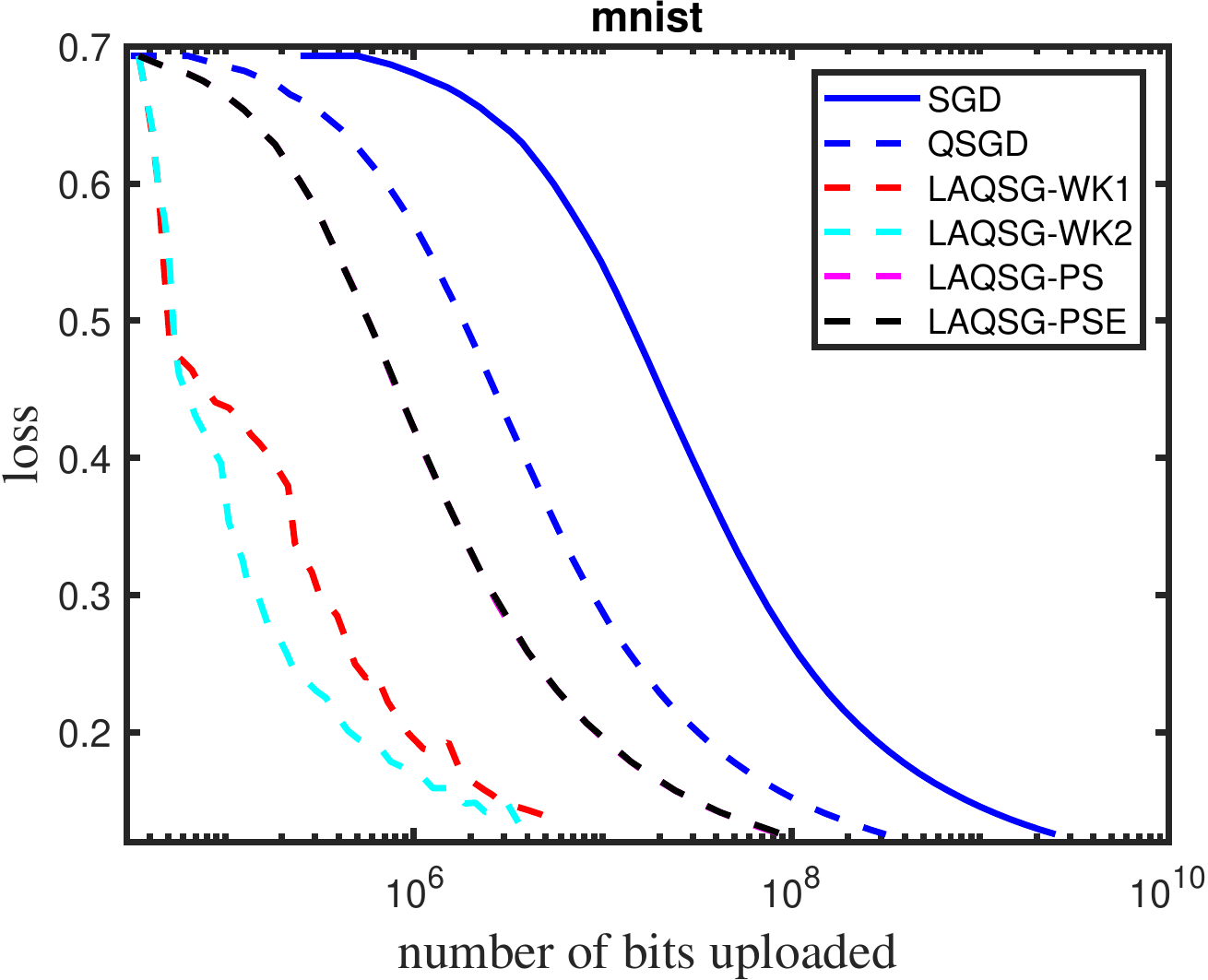}
    \caption{Logistic regression on \textit{mnist} digits 3 and 5 in the heterogeneous setting}
    \label{fig:mnist}
\end{figure*}


\begin{figure*}[ht]
    \centering
    \hspace*{-2ex}
    \includegraphics[width=.33\textwidth]{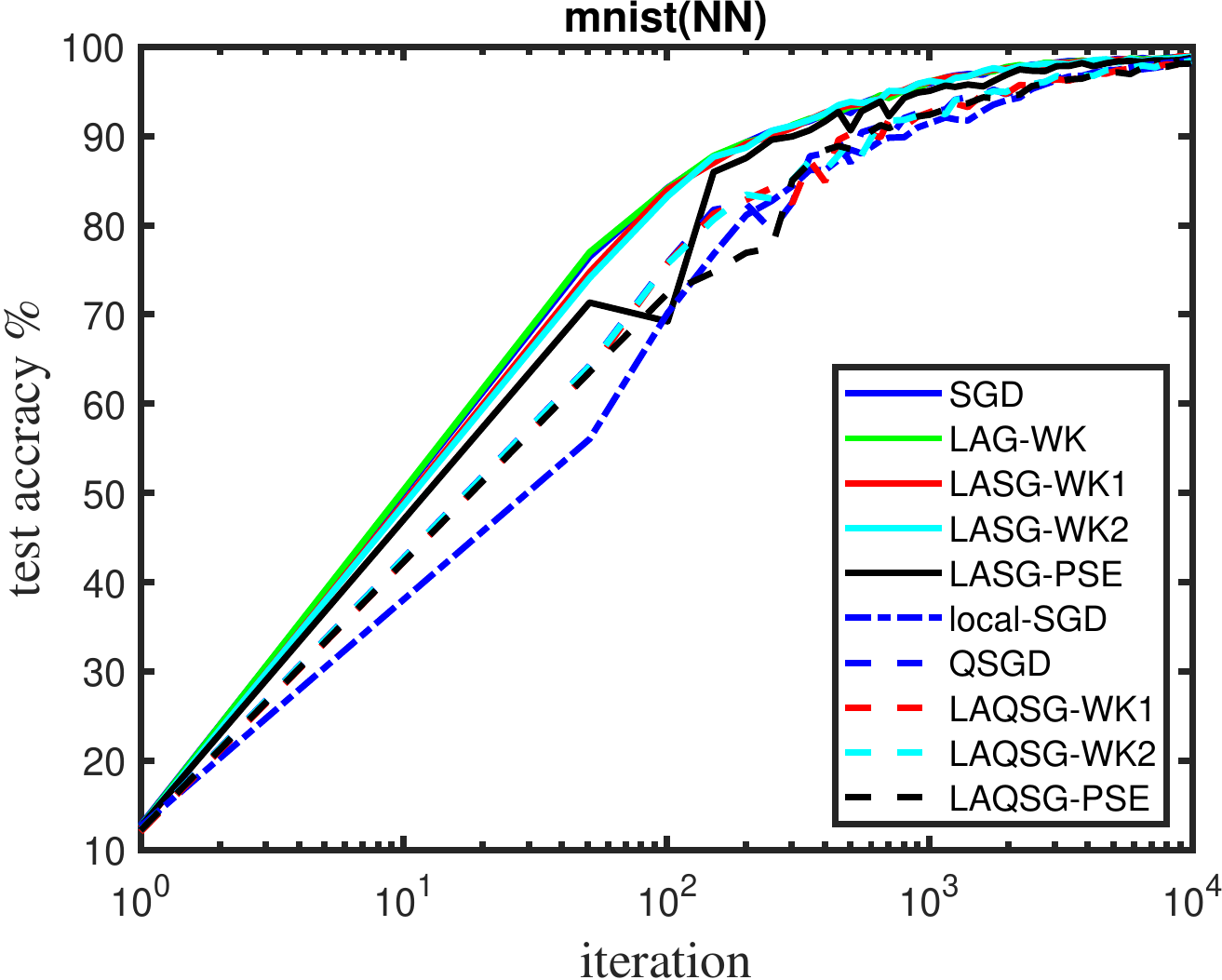}
    \hspace*{-1ex}
    \includegraphics[width=.33\textwidth]{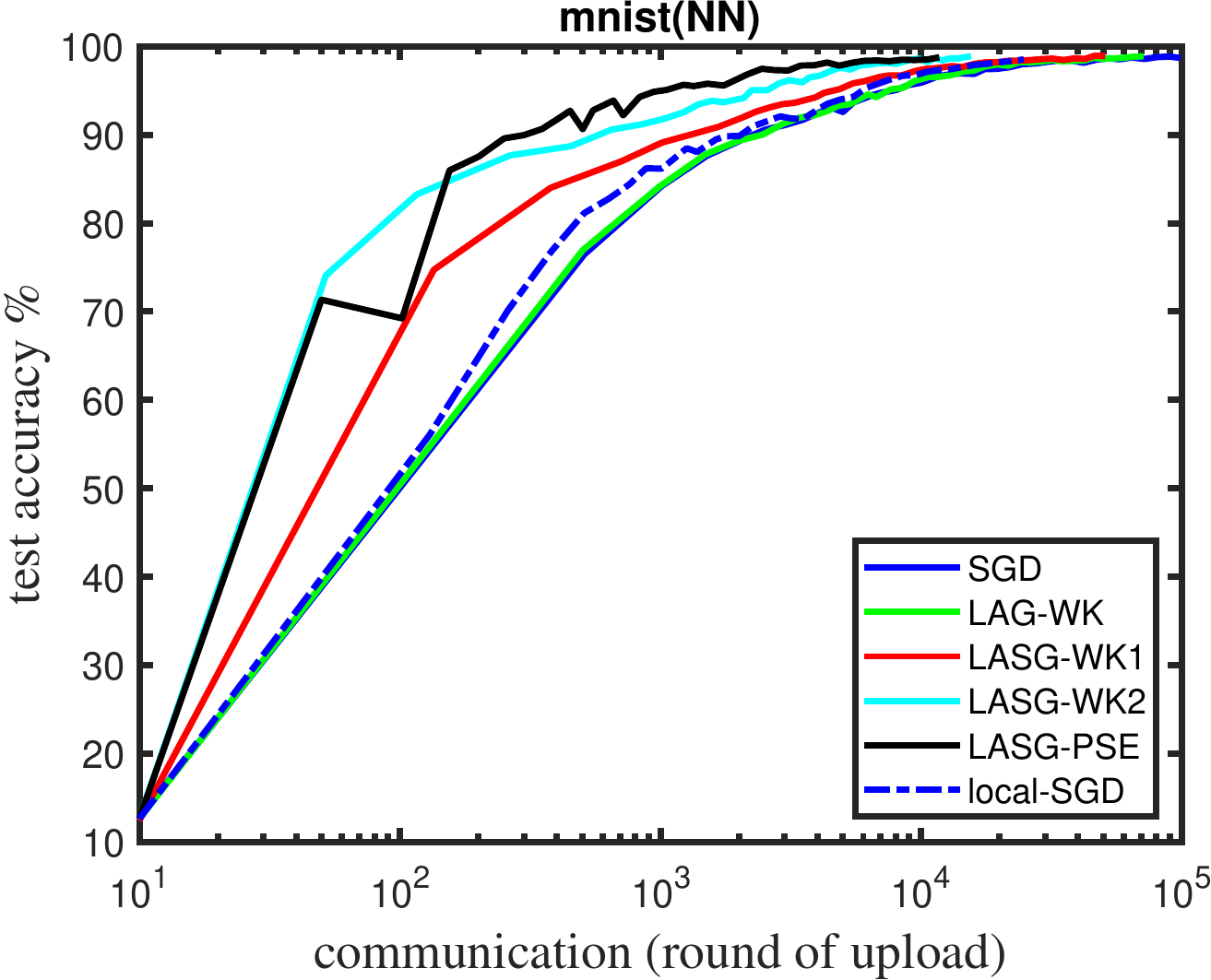}
    \hspace*{-1ex}
    \includegraphics[width=.33\textwidth]{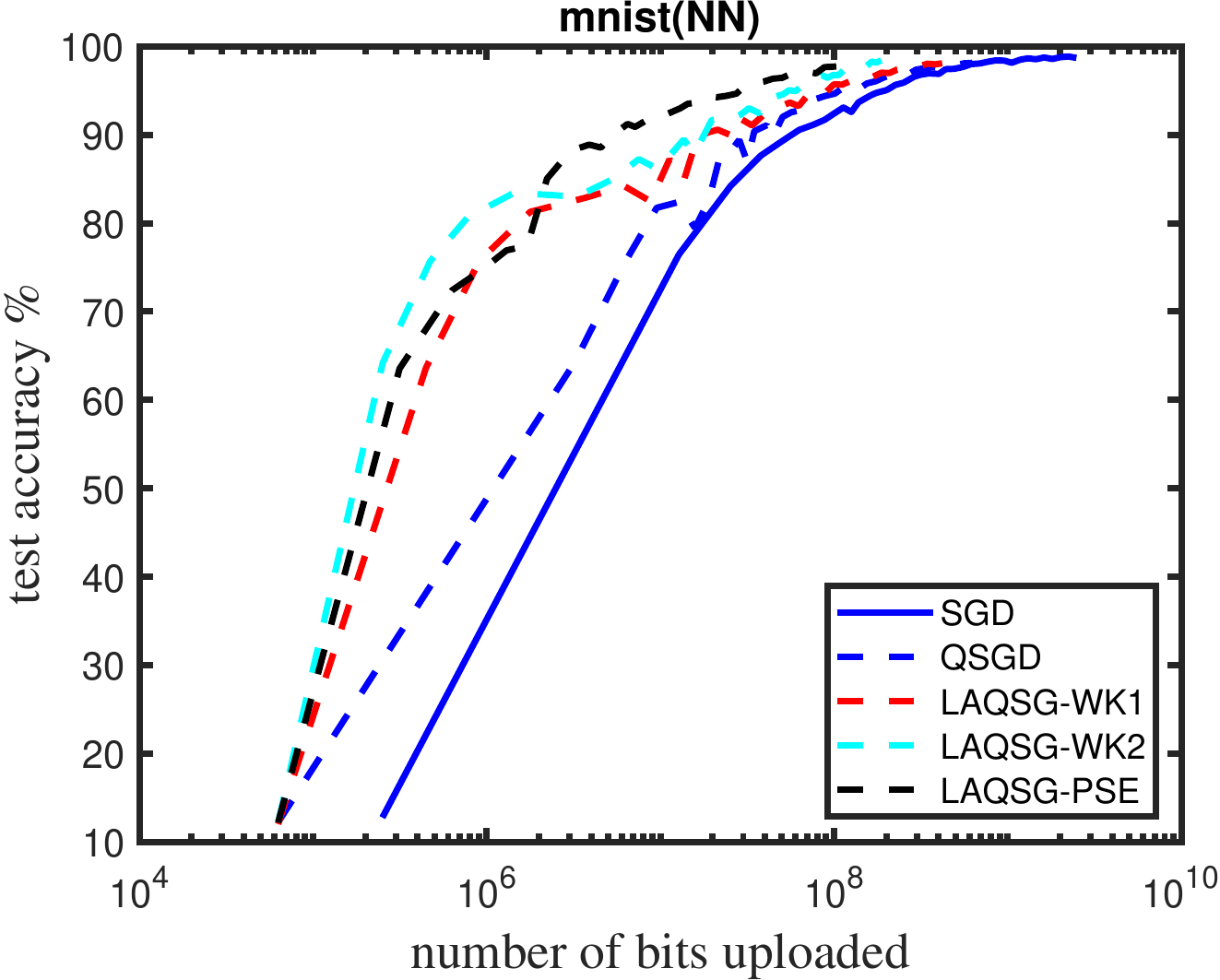}
    \caption{Neural network on dataset \textit{mnist} in the heterogeneous setting}
    \label{fig:NNmnist_test}
\end{figure*}


\end{document}